\input amstex
\documentstyle{amsppt}
\magnification1200
\tolerance=500
\def\n#1{\Bbb #1}

\def\fr{\hbox{fr}}

\def\Gal{\hbox{Gal }}

\def\Hom{\hbox{Hom}}

\def\Res{\hbox{Res }}

\def\det{\hbox{det }}

\def\Supp{\hbox{ Supp }}
\def\Proj{\hbox{Proj }}
\def\Sing{\hbox{Sing}}
\def\Spec{\hbox{ Spec }}

\def\s{\sigma}

\def\e11{E_{11}}

\def\ga{\goth A}

\def\ve{\varepsilon}
\def\vf{\varphi}
\def\de{\delta}

\def\ga{\gamma}
\def\la{\lambda}
\def\La{\Lambda}

\def\be{\beta}
\def\al{\alpha}
\def\sg{\hbox{ sgn}}

\def\suchthat{\hbox{ such that }}
\topmatter
\title
Resultantal varieties related to zeroes of L-functions of Carlitz modules
\endtitle
\author
 A. Grishkov, D. Logachev
\endauthor
\NoRunningHeads
\address
Departamento de Matem\'atica e estatistica
Universidade de S\~ao Paulo, Brasil; DM, ICE, Universidade Federal do Amazonas, Manaus, Brasil
\endaddress
\thanks Thanks: The authors are grateful to FAPESP, S\~ao Paulo, Brazil for a finansial support (process No. 2013/10596-8). The second author is grateful to Max-Planck-Institut f\"ur Mathematik, Bonn; Hausdorff Institut f\"ur Mathematik, Bonn; to Sinnou David, R\'egis de la Br\`eteche, Lo\"{\i}c Merel, Institut de Math\'ematiques de Jussieu, UPMC, Paris, for invitations in academic years 2012--13, 2013--14. The authors are grateful to a reader who advised to include the Section 7 to the paper and gave the idea of the proof of Propositions 7.2, 7.3, and to Suvrit Sra, Gjergji Zaimi, Christian Stump, Christian Krattenthaler who contributed to the proof of III.6. 
\endthanks
\abstract We show that there exists a connection between two types of objects: some kind of resultantal varieties over $\n C$, from one side, and varieties of twists of the tensor powers of the Carlitz module such that the order of 0 of its $L$-functions at infinity is a constant, from another side. Obtained results are only a starting point of a general theory. We can expect that it will be possible to prove that the order of 0 of these $L$-functions at 1 (i.e. the analytic rank of a twist) is not bounded --- this is the function field case analog of the famous conjecture on non-boundedness of rank of twists of an elliptic curve over $\n Q$. The paper contains a calculation of a non-trivial polynomial determinant. 
\endabstract 
\endtopmatter
\document 
CONTENTS\footnotetext{E-mail:  logachev94{\@}gmail.com}
\medskip
0. Introduction. 
\medskip
Part I.  $L$-functions of twisted Carlitz modules (characteristic $p$ case). 
\medskip
1. Definitions of $M$ and of $L(M,T)$. 
\medskip
2. Action of $GL_2(\n F_q)$ on the set of $P$, $\goth C^n_P$ and on $L(\goth C^n_P,T)$.
\medskip
3. Matrix $\goth M(P,n,k)$: explicit formula for $L(\goth C_P^n,T)$.
\medskip
4. Distinguished coset of rank $\ge 1$ in the group of twists. 
\medskip
5. Conjugateness of $\goth M(M)$ and $\goth M(\gamma(M))$ for $\gamma\in GL_2(\n F_q)$, and similar results.
\medskip
6. Numerical results and conjectures.
\medskip
7. Results obtained without application of the Lefschetz trace formula.
\medskip
Part II. Resultantal varieties (characteristic 0 case, the rank at infinity). 
\medskip
8. Case of any $q$.
\medskip
9. Case $q=2$.
\medskip
Part III. Calculation of a determinant. 
\medskip
Appendix: some auxiliary results and remarks. 
\medskip
{\bf 0. Introduction.} 
\medskip
{\bf A. General information.} The present paper contains 3 parts that can be read independently. Part I is inspired by the following problem. Let $E$ be an elliptic curve over $\n Q$. Its (analytic) rank is the order of 0 of $L(E,s)$ at $s=1$. There is 
\medskip
{\bf Problem 0.1.} Are the ranks of twists of $E$ bounded? (Conjecturally not). 
\medskip
We consider the function field case analog of this problem. Let $q$ be a power of a prime $p$. An analog of an elliptic curve is a Drinfeld module of rank 2. Nevertheless, there exist simpler but non-trivial objects --- tensor powers $\goth C^n$ of the Carlitz module $\goth C=\goth C_q$ over $\n F_q$. A twist $\goth C^n_P$ of $\goth C^n$ depends on a polynomial $P=\sum_{i=0}^m a_i\theta^i\in \n F_q[\theta]$, where $a_i\in \n F_q$.  Let $L(\goth C^n_P, T)$ be its $L$-function. The (analytic) rank of $\goth C^n_P$ at 1 is the order of 0 of $L(\goth C^n_P, T)$ at $T=1$, it is denoted by $r_1(\goth C^n_P)$.  There is 
\medskip
{\bf Problem 0.2.} Are $r_1(\goth C^n_P)$ bounded? ($q$, $n$ are fixed, $P$ varies). 
\medskip
A version of the Lefschetz trace formula gives us an explicit expression of $L(\goth C^n_P, T)$: it is the characteristic polynomial of a matrix $\goth M(P,m,*)$ attached to $P$ whose entries belong to $\n F_q[a_0,\dots,a_m][t]$  (Theorem 3.3), $t$ is an independent variable. This means that the set of $(a_0,\dots,a_m)\in A^{m+1}(\n F_q)$ such that $r_1(\goth C^n_P)\ge i$ is the set of $\n F_q$-points of an affine algebraic variety (denoted by $X_1(q,n,m,i)$) over $\bar \n F_q$, i.e. the set of zeroes of some (non-homogeneous) polynomials $D_{**}\in\n F_p[a_0,\dots, a_m]$. 
\medskip
At the first glance, there are too many $D_{**}$'s. If they were independent then the answer to Problem 0.2 would be yes. This looks especially likely for $n>1$ (see (6.18) for the numerical results for the value of rank for the case $q=3$, $n=2$). 
\medskip
It turns out that it is much easier to consider the behaviour (as $P$ varies) of $L(\goth C^n_P, T)$ not at $T=1$ but at $T=\infty$. The corresponding analytic rank is denoted by $r_\infty(\goth C^n_P)$ (or $r_\infty(P)=r_\infty(a_0,\dots,a_m)$ ). We have $r_\infty(P)=r_\infty(\lambda P)$, hence the set of $(a_0:...:a_m)\in P^m(\n F_q)$ such that $r_\infty(\goth C^n_P)\ge i$ is the set of $\n F_q$-points of an algebraic variety (denoted by $X_\infty(q,n,m,i)$) over $\bar \n F_q$. It is the set of zeroes of some other (homogeneous) polynomials $H_{**}\in\n F_p[a_0,\dots, a_m]$. 
\medskip
We have some results showing that $H_{**}$ are highly dependent. Since $D_{**}$ are linear combinations of $H_{**}$ with integer coefficients, we can expect that further research will give us a solution to Problem 0.2. Maybe even for $n>1$ the $r_1(\goth C^n_P)$ would be unbounded! 
\medskip
There exists a natural lift of entries of $\goth M(P,m,*)$ in characteristic 0. This is the subject of Part II (Sections 8, 9). Namely, we use the matrix $\goth M(P,m,*)$ as an initial object (i.e. we do not take into consideration that it comes from the theory of $L$-functions of twisted Carlitz modules) in characteristic 0, we consider its characteristic polynomial and varieties $X_\infty(q,n,m,i)\subset P^m(\n C)$ as varieties of zeroes of its coefficients. It turns out that $X_\infty(q,n,m,i)$ are tightly related to varieties that appear in the classical problem of description of the condition that $q$ polynomials in 1 variable\footnotemark \footnotetext{These polynomials are obtained by "splitting" of $P$, see 8.28.}, of degree $\le \frac mq$, have $\ge i$ common roots (see, for example, [P], [GKZ] and other papers), moreover, they shed new light to this theory. Section 8 contains results that are valid for any $q$, and Section 9 contains much more detailed results for the case $q=2$, that generalizes the theory of the resultant of 2 polynomials. See Introduction to Part II for more details. 

Both Sections 8, 9 can be read independently on the other parts of the paper. Moreover, Section 9 gives a self-contained definition and a list of (mostly conjectural) properties of varieties $X_\infty(2,n,m,i)$, see 9.7. Their study is a problem in the style of nineteen century (classical italian-style algebraic geometry). The most interesting (conjectural) property of these varieties is the following Conjecture 9.3: Supp $X_\infty(2,n,m,i)$ do not depend on $n$, although they depend on $n$ as schemes. For $n=0$ we get a description of the above resultantal variety as a principal irreducible component of some natural complete intersection (of hypersurfaces of zeroes of coefficients of the characteristic polynomial of $\goth M(P,m,*)$ ). It is important that for $q=2$ the matrix $\goth M(P,m,*)$ is the Sylvester matrix of two polynomials $P_{[0]}$, $P_{[1]}$ with interchanged rows. Row interchange is essential: for the Sylvester matrix itself (apparently) we have no such a nice theory. Moreover, in our case $\deg(P_{[0]})-\deg(P_{[1]})=0,\pm1$. We have no analog of this theory for the case of the difference of degrees $\ne 0,\pm1$. 
\medskip
It turns out that for $q=2$ there is a simple expression for det $\goth M(P,m,*)$. Its proof is the subject of Part III. Particularly, it implies immediately Conjecture 9.3 for $i=1$. The proof is a direct, elementary combinatorial calculation, rather tedious and long, it is completely independent on the first two parts of the paper. 
\medskip
The present paper opens a way to future research. It will be necessary to prove results of 9.7, to answer questions 9.7.14, to get analogs of 9.7 for $q>2$. 
\medskip
The same questions arise for varieties $X_1(q,n,m,i)$. First, we should find their dimension in characteristic 0, second, in characteristic $p$. See Remark 8.2.4 for the simplest calculation that has to be made. This will open a way to solve Problem 0.2. Further, we should get analogs of these results for the case of Carlitz modules over any curve over $\n F_q$, not necessarily $P^1(\n F_q)$, for Drinfeld modules of rank $>1$ etc. Most likely there will be analogs of the calculation of Part III for these objects. For example, many researchers study resultants of polynomials in several variables. To get their analogs in the present theory, we should find / define Drinfeld module-type objects depending on polynomials in several variables. 
\medskip
Because of the natural action of $GL_2(\n F_q)$ on $P^m$ and concordance of the $L$-function with this action (see Section 2), these problems should be solved in the quotient space (ind-stack) of $\underset {\longrightarrow}\to{\hbox {lim}}\  P^{m+1}(\bar \n F_q)$ by the action of $GL_2(\n F_q)$. 
\medskip
Further, it would be interesting to check whether these quotients of $X_*(q,n,m,i)$ are good from the point of view of coding theory (see [TVN] for the introduction to coding theory).
\medskip
Conjecturally, all irreducible components of $X_\infty(q,n,m,i)$ are rational varieties. We can expect that if we consider twists of a Drinfeld module $\vf$ instead of twists of $\goth C^n$, then irreducible components of the corresponding variety $ X_{\vf,\infty}(q,n,m,i)$ are non-rational varieties. Particularly, they can be elliptic curves over $\bar \n F_q$. If so, we can associate to $\vf$ these elliptic curves over $\bar \n F_q$ which looks a non-trivial construction. 
\medskip
The notion of ordinary determinantal variety can be greatly generalized, see, for example, [FP], [GKZ]; we can expect that there exist similar generalizations of the varieties $X_*(q,n,m,i)$. A natural question: what is the projective dual of $X_*(q,n,m,i)$? Have we an analog of a simple formula $D(r,m,n)^\vee=D(m-r,m,n)$ ([GKZ], 1.4.11), where $D(r,m,n)$ is the determinantal variety of matrices of rank $\le r$ in $P^{mn-1}$ --- the space of all $m\times n$ matrices? 
\medskip
{\bf B. More details. } The idea to consider twists $\goth C^n_P$ is inspired by the analogy with the number field case. If $E$ is an elliptic curve over $\n Q$ and $L(E,s)$ its $L$-function, then the order of 0 of $L(E,s)$ at $s=1$ --- the center of the symmetry of the functional equation for $L(E,s)$ --- is called the analytic rank of $E$. It is an important invariant of $E$, it enters in the statement of the Birch and Swinnerton-Dyer conjecture. Let $E_D$ be the twist of $E$ by a quadratic field $\n Q(\sqrt{D})$. The sign of the functional equation for $L(E_D,s)$ defines the parity of the analytic rank of $E_D$, it depends on $\chi(D)$ where $\chi$ is a quadratic character, see for example [Sh] for the exact statement. The set of all twists of $E$ is an abelian group $G=\Hom(\Gal(\n Q), \n Z/2)$ (if $Aut(E)\ne \n Z/4, \ \n Z/6$), and the set of even twists is a subgroup $G_0\subset G$ of index 2, hence the set of odd twists is the coset $G-G_0$. Conjecturally, for almost all $D$ the rank takes the minimal possible value, i.e. 0 for the even case and 1 for the odd case. Nevertheless, (rare) jumps, i.e. values of rank $\ge 2$, occur; conjecturally, for any $r$ there exists infinitely many $D$ such that the rank of $E_D$ is $r$, although the (conjectural) asymptotics of these $D$ is not known. 

Let us consider the function field case. Let $M$ be an Anderson t-motive over $\n F_q(\theta)$. Let $L(M,s)$, $s\in S_\infty:=\n C_\infty^*\times \n Z_p$, be the $L$-function of $M$ --- this function or its versions for $\tau$-sheaves and crystals are defined for example in [B05], Definition 15, [B02], Definition 2.19, the original definition of its first version is due to [G], 3.4.2a. We shall consider its version $L(M,T)\in (\n F_q[t])[[T]]$ which is defined for example in [G], 3.2.15, or [TW], formula (2.1), or in [L]; this version is called a na\"{\i}ve $L$-function in [B12]. 
\medskip
{\bf Remark 0.3.} In most earlier cited papers the $L$-function $L(M,s)$, first version, is considered as a function in variable $s$. In terms of $\goth C^n$, this is the same as to consider $n$ as a variable: dilatation of $s$ to an integer $n$ corresponds to the tensor multiplication of $M$ by $\goth C^n$: see for example [B05], Proposition 9, or (7.1) of the present paper for the trivial $M$. Unlike this approach, we consider $n$ fixed, and we consider $T$ as a variable. 

Earlier some results on vanishing of $L(\goth C^n,T)$ at $T=1$ were obtained in [T], [L].  
\medskip
{\bf Remark 0.4.} Since there is no functional equation for $L(M,T)$, \footnotemark \footnotetext{Formulas 2.4 can be considered as an analog of the functional equation, but they are too elementary, and no one of them gives the notion of its center of symmetry.} the choice of $T=1$ does not seem too natural. Nevertheless, there is no essential difference between investigation of the zero at $T=1$ and at $T=c$ for $c\in \n F_q^*$, see Remark 4.3 for details, and the choices of $T=1$ and $T=\infty$ are the simplest possible ones. For example, [L], Proposition 2.1, p. 2604 can be considered as an analog of the strong form of the Birch and Swinnerton-Dyer conjecture for $L(M,T)$ at $T=1$. Really, in Part II we start to consider the case of zero at $T=$ any point. 
\medskip
Therefore, we consider twists $\goth C^n_P$ and their $L$-funcions $L(\goth C^n_P,T)$. Twists $\goth C^n_{P_1}$, $\goth C^n_{P_2}$ are isomorphic over $\n F_q(\theta)$ iff $P_1/P_2 \in \n F_q(\theta)^{*(q-1)}$, hence there is no twists if $q=2$ (it is necessary to emphasize that conversely the characteristic 0 case analog of the present theory is the most interesting for $q=2$, see II, Section 9, and III). 

As it was mentioned above, the set of $P=\sum_{i=0}^m a_i\theta^i$ such that $r_1(\goth C^n_P)\ge i$, resp. $r_\infty(\goth C^n_P)\ge i$, is the set of $\n F_q$-points of an affine algebraic variety $X_1(q,n,m,i)$, resp. projective  algebraic variety $X_\infty(q,n,m,i)$ over $\bar \n F_q$ (we use the same notations for $X_*(q,n,m,i)$, $(*=1, \infty)$ and their lifts to $\bar \n Q$). 
\medskip
{\bf Remark 0.5.} \footnotemark \footnotetext{The authors are grateful to a reader who indicated them the subject of the present remark.} We have: only the sets $X_*(q,n,m,i)(\n F_q)$ are canonically defined (as sets of twists having $r_*\ge i$), but not $X_*(q,n,m,i)$ as varieties or schemes. Really, in principle $X_*(q,n,m,i)(\n F_q)$ can be defined as the sets of $\n F_q$-zeroes of other systems of polynomials (not of $D_{**}$, $H_{**}$), hence the sets of $\bar \n F_q$-zeroes of these polynomials can be another. The same holds for $X_*(q,n,m,i)\subset P^m(\bar \n Q)$ --- they depend on lifts of coefficients of $D_{**}$, $H_{**}$ to $\n Z$ which clearly are not canonical. For example, since $X_*(q,n,m,i)(\n F_q)\subset A^{m+1}(\n F_q)$ or $\subset P^m(\n F_q)$ are finite sets, it is meaningless to consider even its dimension.

Really, the situation is not too bad, because the sets $X_*(q,n,m,i)(\n F_q)$ form a concordant system (an ind-variety) as $m\to\infty$: $$X_1(q,n,m-(q-1),i)(\n F_q)=X_1(q,n,m,i)(\n F_q)\cap A^{m+1-(q-1)}$$ where $A^{m+1-(q-1)}\subset A^{m+1}$ is the subspace \{the last $q-1$ coordinates are 0\}. See Remark 4.2 for the similar formula for $X_\infty$. Most likely it will be possible to define the notion of codimension, irreducible components, their degrees etc. of these ind-varieties, this is a subject of future research. Probably the codimension of $X_1(q,n,m,i)(\n F_q)$ in $A^{m+1}$ (conjecturally, it does not depend on $m$) will be defined in terms of $$\underset{m\to\infty}\to{\hbox{lim}} \ \ \hbox{log}_q(\ \#X_1(q,n,m,i)(\n F_q)\ /\ \#A^{m+1}(\n F_q)\ )$$ Degrees of irreducible components can be defined in terms of counting of quantities of elements in $X_1(q,n,m,i)(\n F_q)$ crossed with linear subspaces of complementary dimension, etc. For example, a strong evidence for codimension conjecture 8.7 for $q=3$, case $i>1$ comes from counting of $\#X_\infty(3,1,m,i)(\n F_q)$ for large $m$ and $i$, see Table 6.15. 
\medskip
Further, it is necessary to emphasize that $D_{**}$, $H_{**}$ are the most natural systems of polynomials such that $X_*(q,n,m,i)(\n F_q)$ are their $\n F_q$-zeroes. 
\medskip
Finally, $H_{**}$ of Section 8 are the most natural lifts to $\n Z$ of $H_{**}$ in characteristic $p$ (we would be very wondered if for other lifts of $H_{**}$ the analog of Theorem 8.6, Conjectures 8.7, 9.7 were true!). For example, if $q$ is not a power of $p$ then there is no $X_*(q,n,m,i)(\n F_q)$ in finite characteristic at all, and if $q=2$  they are trivial, while in both these cases $X_\infty(q,n,m,i)$ are non-trivial objects. Moreover, in characteristic 2 there is no difference between signs + and $-$ (i.e. $-1$ can be lifted to 1), while if we change randomly the signs in the lift of (3.1) in characteristic 0 then the results of Section 9 are not (most likely) true. 
\medskip
Moreover, we have 
\medskip
{\bf Conjecture 0.6.} $X_*(q,n,m,i)$ over $\bar \n F_q$ and $X_*(q,n,m,i)$ over $\bar \n Q$ have the same properties, i.e. their dimensions, quantities of irreducible components, degrees, singularities etc. coincide. 
\medskip
This coincidence (first of all, coincidence of dimensions) does not hold if we choose other lifts of $D_{**}$, $H_{**}$.
\medskip
The paper is organized as follows. Section 1 contains definitions of Anderson t-motives, their $L$-functions, tensor powers of Carlitz modules and their twists. We consider in Section 2 the action of the group $GL_2(\n F_q)$ (really, of some larger monoid) on the set of polynomials $P$. If two polynomials $P_1$, $P_2$ belong to the same orbit of this action then there exists a relation between $L(\goth C^n_{P_1},T)$ and $L(\goth C^n_{P_2},T)$. We give the list of these relations. Section 3 is devoted to the proof of Theorem 3.3 giving us the explicit form of $L(\goth C^n_P,T)$. It is based on the Lefschetz trace formula. We define a matrix $\goth M(P,*)$ that plays a key role in the subsequent theory. We show in Section 4 that there exists a coset of index $(q-1)^2$ in the group of twists of $\goth C^n$ (which is $\n F_q(\theta)^*/\n F_q(\theta)^{*(q-1)}$) such that $r$ of the  corresponding twists is $\ge 1$, see Propositions 4.4, 4.5. This coset can be considered as an analog of the coset of index 2 of odd rank of the number field case, although this analogy is far to be complete: $r$ of these twists is $\ge 1$, but not necessarily odd. For the case of non-twisted Carlitz modules existence of these (so-called trivial) zeros of $L(\goth C^n,T)$ was already known, see for example [B12], p. 51, and [T]. 

In principle, relations between $L(\goth C^n_{P_1},T)$ and $L(\goth C^n_{P_2},T)$ existing when $P_1$, $P_2$ belong to the same orbit of $GL_2(\n F_q)$, should not come from relations between $\goth M(P_1,*)$, $\goth M(P_2,*)$. We show in Section 5 that this is not the case --- really, these relations between $\goth M(P_1,*)$, $\goth M(P_2,*)$ exist (for one case we cannot prove existence of these relations, although they should exist). 

Section 6 contains results of computations of $r_1$, $r_\infty$ for polynomials of low degrees. We show that if $D_{**}$ were independent then $r_1$ would be bounded. It turns out that shift-stable polynomials (see 6.9) show higher values of $r_1$. For the case $q=3$, $n=1$ the highest value of $r_1$ is 3 for polynomials of degrees $\le 15$ and 6 for shift-stable polynomials of degrees $\le 39$. We consider problems of correlation between values of $r_1$ and $r_\infty$. 

Some properties of $L(\goth C^n_P,T)$ can be proved without use of the Lefschetz trace formula. These proofs are given in Section 7. 

In Section 8 we consider the varieties $X_\infty(m,i)=X_\infty(q,n,m,i)$ over $\n C$. We show that they contain two subvarieties 
$X_l(m,i)$, $X_r(m,i)$ corresponding to the zeroes of the left and right multiplication by some matrices. Theorem 8.6, Conjecture 8.7 give us formulas for their codimensions, showing that polynomials $H_{**}$ defining $X_\infty(m,i)$ are highly dependent. Examples show that for $q=3$, $n=1$, $i=1$ we have: $X_l(m,1)$, $X_r(m,1)$ are irreducible components of $X_\infty(m,1)$ of the same codimension. Section 9 contains a large list of conjectural properties of $X_\infty(2,n,m,i)$ ($n=0,1$), some proofs and examples. Particularly, for $q=2$ we have $X_l(2,n,m,i)=X_r(2,n,m,i)$, it is the principal irreducible component of $X_\infty(2,n,m,i)$. A conjectural description of other components is given. For more details and directions of further research see the Introduction to Part II.
\medskip
Part III is devoted to the calculation of det $\goth M(P,m,*)$.
\medskip
{\bf Part I.  $L$-functions of twisted Carlitz modules (characteristic $p$ case). }
\medskip
{\bf 1. Definitions of $M$ and of $L(M,T)$.} Standard reference for t-motives is [A86], we use its notations. 
The Anderson ring $\n F_q(\theta)[ t, \tau]$ is the ring of 
non-commutative polynomials over $\n F_q(\theta)$ satisfying the following relations:
$$t\theta=\theta t, \ t\tau = \tau t, \ \tau \theta = \theta^q \tau$$
We need the following (less general than in [A86]) version of the definition of Anderson t-motives $M$ over $ \n F_q(\theta)$: 
\medskip
{\bf Definition.} An Anderson t-motive $M$ is a $\n F_q(\theta)[ t, \tau]$-module such that 
\medskip
(1.1) $M$ considered as a $\n F_q(\theta)[ t]$-module is free of finite dimension $r$; \footnotemark \footnotetext{The number $r$ is called the (ordinary) rank of $M$. It should not be confused with the analytic rank of $M$ at 1 sometimes denoted by $r$ as well. Throughout the present paper we consider only the case of $M$ of ordinary rank 1.}
\medskip
(1.2) $M$ considered as a $\n F_q(\theta)[\tau]$-module is free of finite dimension $n$;  
\medskip
(1.3) $\exists k>0$ such that $(t-\theta)^k M/\tau M=0$.
\medskip
Equivalently, we can consider $M$ as a free finite-dimensional $\n F_q(\theta)[ t]$-module endowed with a map $\tau:M\to M$ satisfying $\tau(\theta m)=\theta^q\tau(m)$, $\tau(t m)=t\tau(m)$ such that conditions equivalent to (1.2), (1.3) hold. 

$L(M,T)$ is defined for example in [L], upper half of page 2603 ($\fr$, $\tau$, $\vf$ of the present paper are respectively $\tau$, $u$, $\tau$ of [L]. Sorry.) Its explicit definition is the following. Let $Q\in M_{r\times r}(\n F_q(\theta)[ t])$ be the matrix of multiplication by $\tau$ in a $\n F_q(\theta)[ t]$-basis of $M$. Let $\goth P$ be an irreducible polynomial in $ \n F_q[\theta]$. $M$ is called good at $\goth P$ if there exists a $\n F_q(\theta)[ t]$-basis of $M$ such that all entries of $Q$ are integer at $\goth P$. The set of bad primes is denoted by $S$. 

We need the following notation. For $a\in (\n F_q[\theta]/\goth P)[ t]$, $a=\sum c_i t^i$ where $c_i\in \n F_q[\theta]/\goth P$, we denote $a^{(k)}:=\sum c_i^{q^k} t^i$, for a matrix $A=(a_{ij})\in M_{r\times r}((\n F_q[\theta]/\goth P)[ t])$ \ $A^{(k)}:=(a_{ij}^{(k)})$ and $A^{[k]}:=A^{(k-1)}\cdot\dots\cdot A^{(1)}\cdot A$.

The local $\goth P$-factor $L_\goth P(M,T)$ is defined as follows $(\goth P\not\in S)$. Let $d$ be the degree of $\goth P$ and $\tilde Q\in  M_{r\times r}((\n F_q[\theta]/\goth P)[ t])$ the reduction of $Q$ at $\goth P$. We have:
$$L_\goth P(M,T):= \det(I_r-\tilde Q^{[d]}T^d)^{-1}\in \n F_q[t][[T^d]]$$ (because obviously det $(I_r-\tilde Q^{[d]}T)\in \n F_q[t,T]$ and does not depend on a $\n F_q(\theta)[ t]$-basis of $M$); 
$$L_S(M,T):=\prod_{\goth P\not\in S} L_\goth P(M,T)\in \n F_q[t][[T]]$$

The Carlitz module $\goth C$ is an Anderson t-motive having $r=n=1$. Let $\{e\}=\{e_1\}$ be the only element of a basis of $M$ over $\n F_q(\theta)[ \tau]$. $\goth C$ is given by the equation $te=\theta e +\tau e$, i.e. its $Q$ is $t-\theta$. Further, $\goth C^n$ --- the $n$-th tensor power of $\goth C$ --- has the ordinary rank $r=1$, dimension $n$, its ($1\times 1$)-matrix $Q$ is $(t-\theta)^n$. For $P=\sum_{i=0}^m a_i\theta^i\in \n F_q[\theta]$ we denote by $\goth C^n_P$ its $P$-twist whose $Q$ is $P(t-\theta)^n$.  Two such twists $\goth C^n_{P_1}$, $\goth C^n_{P_2}$ are isomorphic over $\n F_q(\theta)$ iff $P_{1}/P_{2}\in \n F_q(\theta)^{*(q-1)}$. Since below we shall change $t$, sometimes we shall denote its local (resp. global) $L$-function $L_\goth P(\goth C^n_P,T)$ (resp. $L(\goth C^n_P,T)$) by $L_\goth P(\goth C^n_P,t,T)$ (resp. $L(\goth C^n_P,t,T)$). See 7.1 for the relation between $L(\goth C^n_P,T)$ and the Goss' $L$-function of the ring $\n F_q[t]$. 
\medskip
{\bf 2. Action of $GL_2(\n F_q)$ on the set of $P$, $\goth C^n_P$ and on $L(\goth C^n_P,T)$.} We fix $q$ and $n$. Really, we consider not only the action of $GL_2(\n F_q)$, but of its direct product with the additive monoid $\n Z^+=\{n|n \ge 0\}$ and with the multiplicative monoid $(\n F_q[\theta])^{q-1}$. The group $GL_2(\n F_q)$ acts tautologically on $\n F_q^2$ and hence on $S^m(\n F_q^2)$. We identify $S^m(\n F_q^2)$ and the set of $P$ of degree $\le m$, hence we get the action of $GL_2(\n F_q)$ on the set of these $P$ and on the set of $\goth C^n_P$. The following lemma gives us explicitly the action of 4 types of generating elements $\mu_{d}:=\left(\matrix 1&0\\ d&1 \endmatrix \right)$, $\nu_{c}:=\left(\matrix c&0\\ 0&1 \endmatrix \right)$, $\iota:=\left(\matrix 0&1\\ 1&0 \endmatrix \right)$, $\tau_{c}:=\left(\matrix c&0\\ 0&c \endmatrix \right)$ of  $GL_2(\n F_q)$ on the set of polynomials $P=\sum_{i=0}^m a_i \theta^i$, where $a_i \in \n F_q$ and $m\equiv -n$ mod $q-1$.
\medskip
{\bf Lemma 2.1.} 1. $\mu_{d}(P)= \sum_{i=0}^m a_i (\theta+d)^i$.
\medskip
2. $\nu_{c}(P)= \sum_{i=0}^m a_i (c\theta)^i$.
\medskip
3. $\iota(P)=\sum_{i=0}^m a_{m-i} \theta^i$.
\medskip
4. $\tau_{c}(P)= c^{-n}P$. $\square$
\medskip
The action of $k\in \n Z^+$ is defined by the following formula: 
\medskip
{\bf Definition 2.1.5.} $\s_k(\goth C^n_P):=\goth C^{q^k\cdot n}_P$.
\medskip
For $\mu_d$, $\nu_c$, $\iota$ we apply the same definitions to $\goth P$.
\medskip
{\bf Remark.} 1. We do not require $a_m=0$ in 2.1.3. Hence, we can choose different values of $m$ satisfying $m\equiv -n$ mod $q-1$, and different $\iota(P)$ are not equal as polynomials, but they are well-defined as an element of $\n F_q(\theta)^{*(q-1)}$. Hence, $\iota$ is well-defined on the set of twists. 
\medskip
2. If $n>1$ then we can consider the action of a slightly larger group denoted by $GL_2(\n F_q)_{(n)}$, namely $GL_2(\n F_q)_{(n)}\subset GL_2(\bar \n F_q)$ is generated by $GL_2(\n F_q)$ and $\tau_{c}:=\left(\matrix c&0\\ 0&c \endmatrix \right)$ where $c\in \bar \n F_q$, $c^n\in \n F_q$. The action of $\tau_{c}$ is given by the formula 2.1.4.
\medskip
It is obvious that the action of $GL_2(\n F_q)_{(n)}\times \n Z^+$ is concordant even with the local $L$-factors of $\goth C^n_P$: 
\medskip
{\bf Lemma 2.2.}  1. $L_{\mu_{d}(\goth P)}(\goth C^n_{\mu_{d}(P)},t-d,T)= L_\goth P(\goth C^n_P,t,T)$; 
\medskip
2. $L_{\nu_{c^{-1}}(\goth P)}(\goth C^n_{\nu_{c}(P)},c^{-1}t,c^nT)= L_\goth P(\goth C^n_P,t,T)$;
\medskip
3. $L_{\iota(\goth P)}(\goth C^n_{\iota(P)},t^{-1},(-t)^nT)= L_\goth P(\goth C^n_P,t,T)$ ($\goth P\ne \theta$);
\medskip
4. $L_\goth P(\goth C^n_{\tau_c(P)},t,c^{n}T)= L_\goth P(\goth C^n_P,t,T)$, where $c\in \bar \n F_q$, $c^n\in \n F_q$;
\medskip
5. $L_\goth P(\s_k(\goth C^{n}_P),t,T)=L_\goth P(\goth C^{n}_P,t^{q^k},T)$.  $\square$
\medskip
{\bf Remark 2.2.6.} The case $\goth P= \theta$ in 2.2.3 corresponds to the point $0\in P^1(\n F_q)$. We have $\iota(0)=\infty$, and 2.2.3 remains true for this case, see below the proof of Theorem 3.3.
\medskip
{\bf Remark 2.2.7.} Formula 2.2.5 shows that the investigation of $\goth C^{q^k\cdot n}_P$ can be reduced to the investigation of $\goth C^{n}_P$. See also 5.5 where it is shown that the matrix $\goth M(P,q^kn,\goth k)$ (see 3.1 below) can be expressed in terms of $\goth M(P,n,\goth k)$. For the case $P=1$ this subject is developed in [B12], Section 10: it is shown that it is possible to find an analog of $\goth M(P,n,\goth k)$ whose size $\goth k$ is expressed in terms of the digits of the $q$-digit expansion of $n$, which gives a much lower bound for $\goth k$ than the bound $\goth k\ge\frac{m+n}{q-1}$ used below. Clearly it is possible to get analogous results for $\goth C^n_P$ for any $P$. 
\medskip
For completeness, we mention also the following observation. If $P_2=P_1P^{q-1}$ then $\goth C^{n}_{P_2} =\goth C^{n}_{P_1}$, and we have
\medskip
{\bf Observation 2.3.} $L_\goth P(\goth C^{n}_{P_2}),t,T)=L_\goth P(\goth C^{n}_{P_1}),t,T)$ if $\goth P \not\vert P$.  $\square$
\medskip
{\bf Corollary 2.4.} 1. If $P_2=\mu_{d}(P_1)$ then $L(\goth C^n_{P_2},t-d,T)=L(\goth C^n_{P_1},t,T)$. 
\medskip
2.  If $P_2=\nu_{c}(P_1)$ then $L(\goth C^n_{P_2},c^{-1}t,c^nT)=L(\goth C^n_{P_1},t,T)$. 
\medskip
3.  If $P_2=\iota(P_1)$ and $S=\theta$ then $L_S(\goth C^n_{P_2},t^{-1},(-t)^nT)=L_S(\goth C^n_{P_1},t,T)$.\medskip
4.  If $P_2=\tau_{c}(P_1)$ then $L(\goth C^n_{P_2},t,c^{n}T)=L(\goth C^n_{P_1},t,T)$, where $c\in \bar \n F_q$, $c^n\in \n F_q$. 
\medskip
5. $L(\s_k(\goth C^{n}_P),t,T)=L(\goth C^{n}_P,t^{q^k},T)$.  $\square$
\medskip
{\bf 3. Matrix $\goth M(P,n,k)$: explicit formula for $L(\goth C_P^n,T)$.} Let $P=\sum_{i=0}^m a_i \theta^i$ be as above and $M=\goth C^n_P$. We denote by $\goth M(P,n,k)=\goth M(P,t,n,k)$ (here $k$ is sufficiently large) the matrix in $M_{k\times k}(\n F_q[t])$ defined by the formula 

$$\goth M(P,n,k)_{i,j}=\sum_{l=0}^n (-1)^l\binom{n}{l}a_{iq-j-l}\ t^{n-l}\eqno{(3.1)}$$ (here $a_*=0$ if $*\not\in [0,\dots, m]$). Particularly, for $n=1$ we have $\goth M(P,1,k)_{i,j}=a_{iq-j}t-a_{iq-j-1}$ and
\medskip
\noindent
$\goth M(P,1,k)=\left(\matrix a_{q-1}t-a_{q-2} &  a_{q-2}t-a_{q-3}  &   \dots & a_{q-k}t-a_{q-k-1} \\ a_{2q-1}t-a_{2q-2} 
 &  a_{2q-2}t-a_{2q-3} &   \dots & a_{2q-k}t-a_{2q-k-1} \\
a_{3q-1}t-a_{3q-2} &  a_{3q-2}t-a_{3q-3}  &   \dots & a_{3q-k}t-a_{3q-k-1} \\
\dots & \dots & \dots  & \dots \\ a_{kq-1}t-a_{kq-2}
 &   a_{kq-2}t-a_{kq-3}  &  \dots & a_{kq-k}t-a_{kq-k-1} \endmatrix \right) \ \ \ (3.2)$
\medskip
{\bf Remark.} $\goth M(P,n,k)$ is (up to a non-essential change of indices) a particular case of the matrix from [FP], (1.5). 
\medskip
{\bf Theorem 3.3.} $L(\goth C^n_P,T)$ is the stable value of $\det(I_k - \goth M(P,n,k)T)$ as $k \to \infty$ (more exactly, for any $k\ge\frac{m+n}{q-1}$, see below). 
\medskip
This follows immediately from the Lefschetz trace formula (see for example (1) of [L], page 2603). The matrix $\goth M(P,n,k)$ is the matrix of $\tau \circ \vf$ on $H^1(P^1, \Cal E)$, see below. Before giving a proof, we need some definitions. The Lefschetz trace formula holds for a slightly different object called $\Cal E$-$\tau$-sheaf ([L], page 2603). Let us recall its definition in the form that we need. Let $P^1$ be the projective line over $\n F_q$ with the function field $\n F_q(\theta)$ and $\fr: P^1\to P^1$ the Frobenius map. The map $(fr, Id): P^1\times \Spec \n F_q[t] \to P^1\times \Spec \n F_q[t]$ is denoted by $fr$ as well. 
\medskip
{\bf Definition.}  A $\Cal E$-$\tau$-sheaf\footnotemark \footnotetext{This is a particular case of the general definition of [L].} is a pair $(\Cal E, \tau)$ where $\Cal E$ is a locally free sheaf on $P^1\times \Spec \n F_q[t]$ and $\tau$ is a $P^1\times \Spec \n F_q[t]$-linear map $\fr^*(\Cal E)\to \Cal E$. 
\medskip
Let $U_0=P^1-\{\infty\}$, $U_1=P^1-\{0\}$ be Zariski open subsets of $P^1$.
\medskip
{\bf Remark.} We do not require that the restrictions of $(\Cal E, \tau)$ to $U_0\times \Spec \n F_q[t]$, $U_1\times \Spec \n F_q[t]$ satisfy (1.2), (1.3), because we do not need this assumption. 
\medskip
The definition of $L$ extends to $\Cal E$-$\tau$-sheaves; clearly the product includes the point $\infty \in P^1$, and --- because $\tau$ is Zariski-locally over $\n F_q[\theta,t]$ (and not over $\n F_q(\theta)[t]$) --- we see that the set of bad points $S$ is empty. 
\medskip
We need also a skew map $\vf: \Cal E \to \fr^*(\Cal E)$ (it is denoted by $\tau$ in [L], p. 2603, 8-th line above the formula (1)). For the affine case its definition is the following. Let $X=\Spec A$ and $L$ be a coherent sheaf on $X$ corresponding to an $A$-module $M$. The sheaf $\fr^*(L)$ corresponds to the module $M\otimes _A A$ respectively the Frobenius map $A\to A$. At the level of modules the map $\vf: M \to M\otimes _A A$ is defined by $m \mapsto m\otimes 1$; we have $\vf(am)=am\otimes 1 = m\otimes a^{(1)} = a^{(1)} \vf(m)$. This definition obviously extends to the case of any scheme, as well as to cohomology. 
\medskip
{\bf Theorem (Lefschetz trace formula)} $$L(\Cal E, \tau,T)= \frac{\det (1- H^1(P^1, \tau \circ \vf ) \cdot T)} {\det (1- H^0(P^1, \tau \circ \vf) \cdot T)}\eqno{(3.4)}$$ 

For a proof of (3.4) see, for example, [B12], Section 9, or the original paper [A00]. 
\medskip
{\bf Proof of Theorem 3.3.} To apply (3.4) to $L(\goth C_P^n,T)$ we should construct first
\medskip
{\bf (3.5)} a $\Cal E$-$\tau$-sheaf whose restriction to $U_0\times \Spec \n F_q[t]$ is $\goth C_P^n$. 
\medskip
Let $\Cal E=\pi^*(O(\goth n))$ where $\pi:P^1\times \Spec \n F_q[t]\to P^1$ is the projection. We have $\fr^*(\Cal E)=\pi^*(O(q\goth n))$. We denote by $e_i$ (resp. $f_i$), $i=0,1$, the only element of a basis of $\Cal E(U_i\times \Spec \n F_q[t])$ (resp. $\fr^*(\Cal E)(U_i\times \Spec \n F_q[t])$ ) over $O(U_i\times \Spec \n F_q[t])$, so $e_1=\theta^\goth n e_0$ in $\Cal E((U_0\cap U_1)\times \Spec \n F_q[t])$, $f_1=\theta^{q\goth n} f_0$ in $\fr^*(\Cal E)((U_0\cap U_1)\times \Spec \n F_q[t])$. Condition (3.5) implies $\tau(f_0)= P(t-\theta)^ne_0$, hence 
$$\tau(f_1)=\theta^{(q-1)\goth n} P(t-\theta)^ne_1\eqno{(3.6)}$$ In order to get a map $\tau: \fr^*(\Cal E)\to\Cal E$, we must have $\theta^{(q-1)\goth n} P(t-\theta)^n\in \n F_q[\theta^{-1},t]$, which is equivalent $\goth n\le-\frac{m+n}{q-1}$. We fix one such $\goth n$ and hence $\Cal E$. 
\medskip
{\bf Remark.} According the terminology of [BP], [B12], the pairs $(\Cal E, \tau)$ for different $\goth n$ are different representatives of the same crystal $j_!(\goth C^n_P)$ where $j_!$ is the extension by zero corresponding to the open immersion $j: U_0=A^1\to P^1$. The formula $\goth n\le-\frac{m+n}{q-1}$ for the case $m=0$ (i.e. the non-twisted Carlitz module) was obtained in [B12], Example 5.12 ($m$ of [B12], Example 5.12 = $\goth n$ of the present paper).\footnotemark \footnotetext{The authors are grateful to a reader who indicated them this information.}
\medskip
It is clear that $\vf: \Cal E \to \fr^*(\Cal E)$ is defined by the formulas $\vf(e_i)=f_i$, $i=0,1$. 
\medskip
We denote $k=-\goth n-1$. We have $H^0(\Cal E)=0$, and elements $\theta^{-1}e_0, \dots, \theta^{-k}e_0$ form a basis of $H^1(\Cal E)$. We have $\vf(\theta^{-i}e_0)=\theta^{-qi}f_0$ and 
$$\tau\circ \vf(\theta^{-i}e_0)=\theta^{-iq}P(t-\theta)^ne_0
= \sum_{j\in \n Z} \left(\sum_{l=0}^n t^{n-l}(-1)^l\binom{n}{l}a_{iq-j-l}\right)\theta^{-j}e_0\eqno{(3.7)}$$ 
hence for $\goth n\le-\frac{m+n}{q-1}$ we have $L(\Cal E, \tau,T)=\det(I_k - \goth M(P,n,k)T)$. 

Finally, $$L(\Cal E, \tau,T)=L(\goth C_P^n,T) \cdot L_\infty(\Cal E, \tau,T)\eqno{(3.8)}$$ 
We have $ L_\infty(\Cal E, \tau,T)=1$ if $\goth n\ne -\frac{m+n}{q-1}$ and $$ L_\infty(\Cal E, \tau,T)=(1-(-1)^n a_mT)^{-1}\hbox{ if }\goth n= -\frac{m+n}{q-1}\eqno{(3.9)}$$ This follows immediately from (3.4), or it can be calculated explicitly as follows. (3.6) is written as (the same calculation as in (3.7))
$$\tau(f_1)=\sum_{j\in \n Z} \left(\sum_{l=0}^n t^{n-l}(-1)^l\binom{n}{l}a_{-(q-1)\goth n-j-l}\right)\theta^{-j}e_1 \eqno{(3.10)} $$ 
The reduction at infinity gives us $\theta^{-1}\mapsto 0$. The cofficient at $\theta^{-j}$ for $j<0$ is 0, hence the coefficient at $(\theta^{-1})^0$ in (3.10) is the only term corresponding to $l=n$, hence $m= -(q-1)\goth n -n$ and $\tilde Q_\infty=(-1)^n a_m$. 
\medskip
 In all cases we get the formula for $L(\goth C^n_P,T)$ (it is clear that $\det(I_k - \goth M(P,n,k)T)$ does not depend on $k$ for $k\ge\frac{m+n}{q-1}$, see also the proof of Proposition 4.4). $\square$
\medskip
{\bf 4. Distinguished coset of rank $\ge 1$ in the group of twists. }
\medskip
{\bf Definition 4.1.} The order of 0 of $L(\goth C_P^n,T)$ at $T=1$ is called the analytic rank at $T=1$ of $P$. It is denoted by $r_1=r_1(P)=r_1(P,n)$. The number $r_\infty:=k+1-\hbox{deg}_T(L(\goth C^n_P,T))$ is called the (deficiency of) the rank of $P$ (or of $a_0,\dots, a_m$) at $T=\infty$. 
\medskip
{\bf Remark 4.2.} $r_\infty$ is not invariant under the natural inclusion of the set of polynomials of degree $m$ to the set of polynomials of degree $m'$, where $m'>m$. Namely, we have $X_\infty(q,n,m-(q-1),i-1)=X_\infty(q,n,m,i)\cap P^{m-(q-1)}$ where $P^{m-(q-1)}\subset P^m$ is the subspace of the last $q-1$ coordinates = 0. This concordance relation will permit us to show that the dimension, degree and other invariants of $X_\infty(q,n,m,i)$ are well-defined, see 0.5. 
\medskip
{\bf Remark 4.3.} Corollary 2.4.4 implies that there is no essential difference between inversigation of the zero at $T=1$ and at $T=c$ for $c\in \n F_q^*$: the order of 0 of $L(\goth C^n_{P_2},T)$ at $T=c$ is equal to the order of 0 of $L(\goth C^n_{P_1},T)$ at $T=1$ if $P_2=cP_1$. Since $\forall M$ we have $L(M,0)=1$, the choices of $T=1$ and $T=\infty$ to find the orders of zero are apparently the simplest ones.
\medskip
{\bf Proposition 4.4.} If $m\equiv -n$ mod $q-1$ and $a_m=(-1)^n$ then $r_1\ge 1$. \footnotemark \footnotetext{See Section 7 for a direct (without using of the Lefschetz trace formula) proof of this proposition.}
\medskip
{\bf Proof.} Follows immediately from (3.8), (3.9). More explicitly, let $i=\frac{m+n}{q-1}$. For $j\ge i$ all elements on the $j$-th line of $I_k - \goth M(P,n,k)T$ to the left from the diagonal are 0, and the diagonal element $(I_k -  \goth M(P,n,k)T)_{jj}$ is 1 for $j> i$, and it is $1-(-1)^na_mT$ for $j=i$. This means that $1-(-1)^na_mT$ is a factor of  $\det(I_k -  \goth M(P,n,k)T)$. $\square$
\medskip
This case corresponds to a coset. Namely, the set of twists of $\goth C$ is isomorphic to $\Hom(\Gal(\n F_q(\theta)), \n Z/(q-1))$. This is a free $\n Z/(q-1)$-module generated by $i_0$, $i_{\goth Q}$ where $i_0$ comes from $\n F_q^*$ and $\goth Q$ runs over the set of places of $\n F_q[\theta]$. Let us consider a homomorpism $\phi$ of this group to $[\n Z/(q-1)]^2=[\n Z/(q-1)]j_1 \oplus [\n Z/(q-1)]j_2$ defined as follows: $i_0 \mapsto j_1$, $i_{\goth Q} \mapsto \deg(\goth Q) j_2$.  
\medskip
{\bf Proposition 4.5.} The set of twists of Proposition 4.4 is $\phi^{-1}(n\frac{q-1}{2}j_1; -n j_2)$ for odd $q$ and $\phi^{-1}(0\cdot j_1; - nj_2)$ for even $q$, i.e. it is a coset of a subgroup of index $(q-1)^2$ of the group of twists. $\square$
\medskip
{\bf 4.6.} We see that if $k:=\frac{m+n}{q-1}-1$ is integer then $L(\goth C^n_P,T)$ is a product of two factors: $$L(\goth C^n_P,T)=\det (I_k-\goth M(P,n,k)T)\cdot (1-(-1)^na_mT)$$ We denote the first factor by $L_{nt}(\goth C^n_P,T)$ --- the non-trivial factor of $L(\goth C^n_P,T)$. Respectively, the order of 0 of $L_{nt}(\goth C_P^n,T)$ at $T=1$ is called the non-trivial part of the analytic rank at $T=1$ of the pair $(P, n)$. It is denoted by $r_{nt\ 1}=r_{nt\ 1}(P,n)$. 
\medskip
{\bf 5. Conjugateness of $\goth M(M)$ and $\goth M(\gamma(M))$ for $\gamma\in GL_2(\n F_q)\times \n Z^+\times (\n F_q[\theta])^{q-1}$.}
\medskip
Corollary 2.4 shows that if $\goth C_{P_2}^{n_2}=\gamma(\goth C_{P_1}^{n_1})$ for $\gamma\in  GL_2(\n F_q)_{(n)} \times \n Z^+\times (\n F_q[\theta])^{q-1}$ then there exists the corresponding relation between their $L$-functions. It is natural to expect that matrices $\goth M(P_1,n_1,k)$ and $\goth M(P_2,n_2,k)$ are conjugate. Let us prove it, separately for 6 types of generators of $GL_2(\n F_q)_{(n)}\times \n Z^+\times (\n F_q[\theta])^{q-1}$. 
\medskip
{\bf Remark.} Equality of characteristic polynomials $|I_k-\goth M_1T|$, $|I_k-\goth M_2T|$ does not imply conjugateness of matrices $\goth M_1$, $\goth M_2$, hence we can consider the contents of the present section as a proof of the theorem that matrices belonging to one orbit of the $GL_2(\n F_q)$-action have the same Jordan type. This is an important invariant --- see, for example, [L], Proposition 2.1, p. 2604, the condition of semi-simplicity of the eigenvalue 1: it means that the lengths of all Jordan blocks having $\lambda=1$ of $\goth M$ are equal to 1. By the way, V. Lafforgue writes (lines 1 -- 2, p. 2604): "Il n'y a aucune raison pour que cette hypoth\`ese de semi-simplicit\'e de la valeur propre 1 soit toujours v\'erifi\'ee"; really, for the case $q=3$, $m=3$, $n=1$ there are 3 polynomials of rank 2 (see table 6.7 below); they form a $\left(\matrix 1&0\\ *&1 \endmatrix \right)$-orbit, and this condition of semi-simplicity does not hold for them. 
\medskip
{\bf 5.1. Type $\mu_d$.} We shall consider infinite matrices with entries in $\n F_q$ whose rows and columns are numbered by 0,1,2,..., all operations over the matrices under consideration will be well-defined. Particularly, the matrices $\goth M(P,n,t):=\underset {\longrightarrow}\to{\hbox{lim }} \goth M(P,n,t,k)$ as $k\to+\infty$ are of this type. We fix $d$ and we define a matrix $W=W_1(d)$ as follows: $$W_{ij}=0 \hbox{ if } i>j,  \ \ W_{ij}=\binom{j}{i}d^{j-i}  \hbox{ if } j\ge i\eqno{(5.1.1)}$$ Obviously $W_1(-d)=W_1(d)^{-1}$.
\medskip
{\bf Proposition 5.1.2.} For $P_2=\mu_{d}(P_1)$ we have $$\goth M(P_2,n,t-d)=W\goth M(P_1,n,t)W^{-1}\eqno{(5.1.3)}$$
\medskip
{\bf Proof.} Let as above $P_1=\sum_{i=0}^\infty a_i \theta^i$ where almost all $a_i$ are 0, we denote by $\goth a$ the infinite-to-bottom vector-column $(a_0, a_1, a_2, ... )^t$ and analogously for $P_2=\sum_{i=0}^\infty b_i \theta^i$, $\goth b=(b_0, b_1, b_2, ... )^t$. Therefore, we have $$\goth b=W\goth a\eqno{(5.1.4)}$$ Further, we denote by $\ve_{ij}$ the $(i,j)$-th elementary matrix (its $(i,j)$-th entry is 1 and all other entries are 0; $\ve_{ij}=0$ if $j<0$) and we denote $$\goth M_l:=\sum_{i=0}^\infty \ve_{i,q(i+1)-1-l}\eqno{(5.1.5)}$$ In this notation (3.1) can be rewritten as follows (warning: in (3.1) the rows and columns are numbered from 1 while here from 0):  $$\goth M(P_1,n,t)=\sum_{k=0}^n(-1)^k\binom{n}{k}(\sum_{i=0}^\infty \goth M_{i+k} a_i )t^{n-k}$$ So, (5.1.3) is equivalent to the following formulas for $k=0,\dots,n$ (coincidence of coefficients at $t^{n-k}$): 

$$\sum_{\gamma=0}^k\binom{n}{\gamma}\binom{n-\gamma}{k-\gamma}d^{k-\gamma}(\sum_{i=0}^\infty \goth M_{i+\gamma}b_i)=\binom{n}{k}W(\sum_{i=0}^\infty \goth M_{i+k}a_i)W^{-1}\eqno{(5.1.6)}$$ 
which is equivalent to
$$\binom{n}{k}\left(\sum_{\gamma=0}^k\binom{k}{\gamma}d^{k-\gamma}(\sum_{i=0}^\infty \goth M_{i+\gamma}b_i)-W(\sum_{i=0}^\infty \goth M_{i+k}a_i)W^{-1}\right)=0\eqno{(5.1.7)}$$ 
 Since (5.1.4) and (5.1.6) are linear by $\goth a$, it is sufficient to prove (5.1.6) for $\goth a=\goth a_j:=(0,0,\dots,0, 1,0,\dots)$ ( 1 at the $j$-th place). For this $\goth a_j$ we have $b_i=\binom{j}{i}d^{j-i}$ if $i\le j$, $b_i=0$ if $i> j$, hence (5.1.7) becomes (omitting a non-essential $\binom{n}{k}$)
$$\sum_{\gamma=0}^k\binom{k}{\gamma}d^{k-\gamma}\sum_{i=0}^j \binom{j}{i}d^{j-i}\goth M_{i+\gamma} =W\goth M_{j+k}W^{-1}\eqno{(5.1.8)}$$
The left hand side of (5.1.8) for $k=0$ is $\sum_{i=0}^j \binom{j}{i}d^{j-i}\goth M_i$, and for general $k$ it is 

$$\sum_{\gamma=0}^k\sum_{i=0}^j \binom{k}{\gamma}\binom{j}{i} d^{k+j-(i+\gamma)}\goth M_{i+\gamma}=\sum_{\delta=0}^{k+j} \binom{k+j}{\delta}d^{k+j-\delta}\goth M_\delta$$ (where $\delta=i+\gamma$), hence (5.1.8) for $j=j_0$, $k=k_0$ coincides with (5.1.8) with $k=0$, $j=j_0+k_0$, and 
hence it is sufficient to prove (5.1.8) for $k=0$. First, we consider the case $j=0$, (5.1.8) becomes $\goth M_0W=W\goth M_0$. By (5.1.1) and (5.1.5), this becomes 
$$\binom{l}{q(i+1)-1}=0\hbox{ if }l\not\equiv -1 \mod q\eqno{(5.1.9)}$$ and 
$$\binom{l}{q(i+1)-1}=\binom{\frac{l+1}q-1}{i}\hbox{ if }l\equiv -1 \mod q\eqno{(5.1.10)}$$ (equalities in $\n F_q$). They are proved as follows. We let $l=\alpha q+c$ where $c\in [0, \dots, q-1]$. Equality $(X+Y)^{\alpha q}=((X+Y)^{\alpha})^ q$ ($X$, $Y$ are abstract letters) implies $$\binom{\alpha q}{\gamma q}=\binom{\alpha}{\gamma}\eqno{(5.1.11)}$$ and$\binom{\alpha q}{\gamma}= 0$ if $\gamma\not\equiv 0$ mod $q$. Further, we have $\binom{l}{q(i+1)-1}=\sum_{\beta=0}^c \binom{c}{\beta} \binom{\alpha q}{q(i+1)-1-\beta}$. If $c\ne q-1$ then all $q(i+1)-1-\beta\not\equiv 0 \mod q$, hence we get immediately (5.1.8). If $c= q-1$ then the only $\beta$ such that $q(i+1)-1-\beta\equiv 0 \mod q$ is $\beta=c$, hence $\binom{l}{q(i+1)-1}= \binom{\alpha q}{q(i+1)-1-(q-1)}$ which is (5.1.10), because of (5.1.11). 
\medskip
The case of any $j$ is similar. We have $(W\goth M_j)_{il}=0$ if $l\not\equiv -j-1 \mod q$ and $(W\goth M_j)_{il}=d^{\alpha-i} \binom{\alpha}{i}$ if $l\equiv -j-1 \mod q$ and $\alpha\ge i$, where $\alpha=\frac{l+1+j}{q}-1$. Further, $(\goth M_\gamma W)_{il}=d^{l-(q(i+1)-1-\gamma)}\binom{l}{q(i+1)-1-\gamma}$ and 
$$(\sum_{\gamma=0}^j d^{j-\gamma}\binom{j}{\gamma}\goth M_\gamma W)_{il}=d^{l+j-(q(i+1)-1)} \sum_{\gamma=0}^j\binom{j}{\gamma}\binom{l}{q(i+1)-1-\gamma}$$ $$=d^{l+j-(q(i+1)-1)}\binom{l+j}{q(i+1)-1}$$ We get immediately the desired taking into consideration that $l+j-(q(i+1)-1)=q(\alpha-i)$ and changing $l$ to $l+j$ in (5.1.9), (5.1.10). $\square$
\medskip
{\bf 5.2. Type $\nu_c$.} Here we let $W_2(c)$  a diagonal matrix whose $i$-th diagonal entry is $c^i$. 
\medskip
{\bf Proposition.} For any $c\in \n F_q^*$ we have: $$\goth M(\nu_c(P), c^{-1}t,n,k)=c^{-n}W_2(c)\ \goth M(P, t,n,k)\ W_2(c)^{-1}$$
\medskip
{\bf Proof.} Obvious. Here rows and columns are numbered from 1; (3.1) gives us $\goth M(\nu_c(P), c^{-1}t,n,k)_{i,j}=c^{qi-n-j}\goth M(P, t,n,k)_{i,j}$. Because of $c^q=c$ we get $\goth M(\nu_c(P), c^{-1}t,n,k)_{i,j}=c^{-n}c^{i-j}\goth M(P, t,n,k)_{i,j}$. $\square$
\medskip
{\bf 5.3. Type $\iota$.}  We choose $m$ such that $(q-1)|(m+n)$, we consider $\iota$ corresponding to this $m$, and we let $k=\frac{m+n}{q-1}-1$. For this case we let: $W_3=\sum_{i=1}^k \ve_{i,k+1-i}$ is the matrix whose elements on the second (non-principal) diagonal are ones and another elements are 0 (the rows and columns are numbered from 1). 
\medskip
{\bf  Proposition.} $W_3\ \goth M(P,t,n,k)\ W_3^{-1}=(-t)^n\ \goth M(\iota(P),t^{-1},n,k)$. 
\medskip
{\bf Proof.} Follows immediately from (3.1) (the conjugation with respect to $W_3$ is the central symmetry with respect to the center of a matrix). $\square$
\medskip
{\bf 5.4. Type $\tau_c$.} We have a trivial equality $\goth M(cP,t,n)= c\goth M(P,t,n)$. 
\medskip
{\bf 5.5. Action of $\n Z^+$.} For this case we define $W=W_5$ by the formula $W_{ij}=t^{j-i}$ if $j\ge i$ and $W_{ij}=0$ if $j< i$, hence $W_5^{-1}$ is defined by the formula $W_{ii}=1$, $W_{i,i+1}=-t$, all other $W_{ij}=0$. 
\medskip
{\bf Proposition.} For any $n$, $P$ we have: $$W^n\ \goth M(P,qn,t)\ W^{-n}=\left(\matrix 0&0\\ \goth A_5& \goth M(P,n,t^q)\endmatrix \right)$$ where sizes of blocks are $n\times n$, $n\times \infty$, $\infty\times n$, $\infty\times \infty$ and $\goth A_5$ is some matrix.
\medskip
{\bf Proof.} Straightforward (induction by $n$, for example). $\square$
\medskip
{\bf 5.6. Multiplication by elements of $\n F_q[\theta]^{q-1}$.} If $P_1 = PQ^{q-1}$ for $Q\in \n F_q[\theta]$ then $\goth C_P$, $\goth C_{P_1}$ are different $F_q(\theta)$-models of a twisted Carlitz module over $\n F_q(\theta)$, hence their $L$functions differ by a factor corresponding to bad points --- irreducible factors of $Q$ which do not enter in $P$. More exactly, if $Q=\prod_i \goth Q_i^{\alpha_i} \cdot \prod_j {\goth Q'_j}^{\alpha'_j}$ is the prime decomposition of $Q$ (where $\goth Q_i$ do not divide $P$ and $\goth Q'_j|P$), then 
$$L(\goth C_{P_1},T)=L(\goth C_P,T)(\prod_i L_{\goth Q_i}(\goth C_P,T))^{-1}\eqno{(5.6.1)}$$

For $Q=\theta$ the matrices $\goth M(P,n,t)$ and $\goth M(P\theta^{q-1},n,t)$ coincide up to a non-essential shift:
$$\goth M(P\theta^{q-1},n,t)=\left(\matrix a_0t^n&0\\ \goth A_6& \goth M(P,n,t)\endmatrix \right)$$ where sizes of blocks are $1\times 1$, $1\times \infty$, $\infty\times 1$, $\infty\times \infty$ and $\goth A_6$ is a matrix column. If deg $Q=1$, i.e. $Q=\theta+b$, $b\in \n F_q$, then Proposition 5.1 gives us immediately the relation between $\goth M(P,n,t)$ and $\goth M(PQ^{q-1},n,t)$. To find this relation for the case of deg $Q>1$ is an exercise for the reader. 

\medskip
{\bf 6. Numerical results and conjectures.}
\medskip
In this section (except 6.18) we shall consider only the case $n=1$, and we shall omit the index $n$. The rank $r_1$ will be denoted by $r$. First, let us mention the following elementary result: 
\medskip
{\bf Proposition 6.1.}  For any $r\le q-1$ there exists $P$ such that the analytic rank of $\goth C_P$ is $r$. 
\medskip
{\bf Proof.} We can take for example $P$ having $a_{i(q-1)-1}=-1$ for $i=1,\dots,r$ and other $a_*=0$. The matrix $\goth M(P,k)$, where $k=r$, is upper-triangular with 1's at the diagonal, hence the proposition. $\square$
\medskip
{\bf Example.} For $q=3$ this polynomial and its $\mu_*$-orbit (see Lemma 2.1.1) are the only polynomials of rank $2$ for $m\le 6$, see table 6.7 below. 
\medskip
{\bf Parameter count.} Let us count the quantity of equations defining 
$X_1(m,r)=X_1(q,1,m,r)$. Let $k$ be the size of $\goth M(P,k)$. Changing variable $U=T^{-1}$ we get $$L(\goth C_P,T)=U^{-k}\sum_{i=0}^k A_iU^i\eqno{(6.2)}$$ where $A_i\in \n F_q[t]$, $\deg A_i=k-i$. Changing variable $V=U-1$ we get $$L(\goth C_P,T)=U^{-k}\sum_{i=0}^k B_iV^i\eqno{(6.2a)}$$ where $B_i\in \n F_q[t]$, $\deg B_i=k-i$, $B_i=\sum_{j=0}^{k-i} H_{ij,1}t^j$, $H_{ij,1}=H_{ij}\in \n F_q[a_0,\dots, a_m]$. \footnotemark \footnotetext{The present $H_{ij}$ do not coincide with $H_{ij}$ of Section 8, although they generate the same ideal.} Condition that the analytic rank is $\ge r_0$ is equivalent to the condition $B_0=\dots=B_{r_0-1}=0$, which gives us $$(k+1)+(k)+\dots+(k+1-(r_0-1))=r_0(k+1)-\frac{r_0(r_0-1)}{2}\eqno{(6.3)}$$ equations in $\n F_q$, where $r_0\le k$.  

Let us find the maximal value of $r$ for which we can find $k$ such that the na\"{\i}ve parameter count predicts existence of $k\times k$ matrix $\goth M$ such that the order of 0 of det$(I_k-\goth MT)$ at $T=1$ is $\ge r$. We take $m=kq-k-1$, $a_{kq-k-1}=-1$, hence the last line of det$(I_k-\goth M(P,k)T)$ gives us a factor $(1-T)$. Formula (6.3) applied to the left-upper $(k-1)\times (k-1)$-minor of $\goth M(P,k)$, $r_0=r-1$, gives us $(r-1)k-\frac{(r-1)(r-2)}{2}$ equations. The quantity of variables $a_i$ is $k(q-1)-1$, hence the question is the following: 
\medskip
{\bf 6.4.} For which $r$ there exists $k\ge r$ such that $$k(q-1)-1\ge (r-1)k-\frac{(r-1)(r-2)}{2}\eqno{(6.5)}$$ 
\medskip
The answer to (6.4) is $r\le 2q-3$, this is the expected maximal value of rank. To formulate a rigorous --- although conditional --- result, we define a projective variety $\bar X_1(q,1,m,r)\subset P^{m+1}(\bar \n F_q)$ --- the projectivization of $X_1(q,1,m,r)(\bar \n F_q)$  --- as the set of zeroes of $\bar H_{ij}$ which are the homogeneization of $H_{ij}$. So, we have got
\medskip
{\bf Proposition 6.6.} For $r\le 2q-3$ there exists $m$ such that dim $\bar X_1(q,1,m,r)\ge 0$. If $\bar X_1(q,1,m,r)\subset P^{m+1}$ is the complete intersection of $\bar H_{ij}$, where $k\ge \frac{m+1}{q-1}$, then for $r> 2q-3$ for any $m$ we have $\bar X_1(q,1,m,r)=\emptyset$. $\square$
\medskip
Analogously, we can ask for which $r$ there exist infinitely many $k$ satisfying (6.5), i.e. for which $r$ we can expect existence of  infinitely many $P$ such that the analytic rank of $\goth C_P$ is $\ge r$. The answer is $r\le q$. But in this case we cannot formulate an analog of the conditional Proposition 6.6, because 5.5 shows that if for some $m$ we have $X_1(q,1,m,r)\ne\emptyset$ then $\forall i>0$ we have $X_1(q,1,m+i(q-1),r)\ne\emptyset$. 
\medskip
{\bf Results of computer calculations.} They are given in the following table 6.7.  We consider the case $q=3$ (case $q=2$ is trivial), and we consider separately the cases of the leading coefficients $a_m=1$ and 2. We consider squarefree $P\in \n F_3[\theta]$. The quantity of these $P$ of the degree $m$ and the leading coefficient $a=a_m$ such that the analytic rank of $\goth C_P$ is $\ge r$ is denoted by $\goth q(m,a,r)$. Table 6.7 covers the case $m\le 15$, $r\ge 2$. The maximal value of $r(P)$ for $m\le 15$ is 3.
\medskip
{\bf Table 6.7.} Numbers $\goth q(m,a,r)$.
\settabs 15 \columns

\+ &&&&&&&$r=2$ \cr
\medskip
\+ $a_m$ & $ m=$  &3&4&5&6&7  &8&9&                 10       &11&12&13&14&15 \cr
\medskip
\+ 1 &                    &   0&0&0&0&0  &3&3                &   0  &      3 & 9 & 12&21&44 \cr
\+ 2 &                    &   3&0&0&0&33&3& 165&   0&717&9&3117&21 &14038\cr
\medskip
\+ &&&&&&&$r=3$ \cr
\medskip
\+ 1 &                 &0&0&0&0&0  &0&  3                  & 0  &     0 & 0 & 0&0&3 \cr
\+ 2 &                   &0&0&0&0&0&0&6                   &0&       0&0& 12 &0&42\cr
\medskip
{\bf Remark.} 1. Corollary 2.4, (2) and (4) implies that $r(\tau_{c^{-1}}\circ \nu_c(P))=r(P)$ hence for even $m$ we have $\goth q(m,1,r)=\goth q(m,2,r)$. 
\medskip
2. Most numbers $\goth q(m,a,r)$ in the above table are multiples of 3, because of Corollary 2.4.1. Exceptions are due to shift-stable polynomials, see below. 
\medskip
3. The ratios $\goth q(m,2,2)/\goth q(m-2,2,2)$ for odd $m=9,11,13,15$ are respectively 5, 4.3455, 4.3472, 4.5037. Does exist its limit as $m \to\infty$, what does it equal? What are the similar limits of $\goth q(m,a,r)/\goth q(m-(q-1),a,r)$ for other $q$?
\medskip
{\bf 6.8. Expected dimensions.} We denote by $X_1(m,a,r)$ the subset of 

\noindent
$X_1(3,1,m,r)(\n F_q)$ consisting of polynomials having $a_m=a$. The above parameter count shows that for the case of complete intersections we have for $m$ odd:
\medskip
dim $X_1(m,2,2)=\frac{m-1}2$, dim $X_1(m,2,3)=0$, $X_1(m,2,r)=\emptyset$ for $r>3$;
\medskip
dim $X_1(m,1,2)=0$, $X_1(m,1,r)=\emptyset$ for $r>2$;
\medskip
and for $m$ even we have dim $X_1(m,2)=0$, $X_1(m,r)=\emptyset$ for $r>2$. 
\medskip
Let us compare these predictions and the entries of Table 6.7 for odd $m$, $a_m=2$ (other cases give us too small numbers).The predicted value of $\#(X_1(15,2,2))$ is $\frac23 \cdot 3^7$ while really it is $14038\sim\frac23 \cdot 3^9$ (the coefficient $\frac23$ appears, because we consider the squarefree polynomials). We see that it is not too likely that they are the complete intersections. See Remark 8.2.4 for the explicit question. 
\medskip
{\bf 6.9. Shift-stable polynomials.}
\medskip
{\bf Definition.} $P\in \n F_q[\theta]$ is called a $\theta$-shift-stable if $\forall d\in\n F_q$ we have $P=\mu_d(P)$. 
\medskip
For the shift-stable case we shall use the same notations as earlier, with the subscript $s$. Obviously $P=\sum_{i=0}^{m_{s}} c_{s,i} (\theta^q-\theta)^i$ where $m_{s}=\frac mq$, $c_{s,i}\in \n F_q$. (2.4.1) implies that for these $P$ \ \ $L(\goth C^n_{P},t,T)$ is $t$-shift-stable, hence $B_i$ of (6.2a) are $t$-shift-stable. We use notations $B_{i}=\sum_{j=0}^{[(k-i)/q]} H_{s,ij}(t^q-t)^j$, $H_{s,ij}\in \n F_q[c_{s,0},\dots,c_{s,m_{s}}]$. 
\medskip
We denote by $\goth q_{s}(m,a,r)$ the quantity of square-free shift-stable polynomials $P$ of a given degree $m$ and the leading coefficient $a$ such that $r(P)\ge r$. 
\medskip
{\bf Table 6.10.} Numbers $\goth q_s(m,a,r)$.
\medskip
\settabs 15 \columns

\+ &&&&&&&&$r=1$ \cr
\medskip
\+ $a_m$ & $ m=$  &3&6&9&12&15  &18&21&                 24       &27&30&33&\ 36&\ \ 39 \cr
\medskip
\+ 1 &                    & 3&0&3&0&36&23&205&97&866&505&3601 &\ 2217 &\ \ 15952\cr
\+2&&3&0&$2\cdot3^2$&0&$2\cdot3^4$&23&$2\cdot3^6$&97&$2\cdot3^8$&505&$2\cdot3^{10}$&\ 2217&\ \ $2\cdot3^{12}$\cr

\medskip
\+ &&&&&&&&$r= 2$ \cr
\medskip
\+ 1 &                    & 0&0&0&0&2&1&15&8&46&24&73&\ 71&\ \ 199\cr
\+ 2 &                    & 0&0&3&0&10&1&93&8&380&24&1747&\ 71&\ \ 7639\cr
\medskip
\newpage
\+ &&&&&&&&$r=3$ \cr
\medskip
\+ 1 &                    & 0&0&0&0&0&1&5&1&7&2&8&\ 3& \ \  12\cr
\+ 2 &                    & 0&0&0&0&0&1&9&1&18&2&43&\ 3& \ \  158\cr
\medskip
\+ &&&&&&&&$r=4$ \cr
\medskip
\+ 1 &                    & 0&0&0&0&0&0&$5^*$&0&3&0&2&\ 0&\ \ 0\cr
\+ 2 &                    & 0&0&0&0&0&0&0&0&4&0&5&\ 0& \ \  $16^*$\cr
\medskip
For odd $m$ and $a_m=2$ all polynomials have $r\ge 1$, hence $\goth q_{s}(m,2,1)=2\cdot3^{m/3-1}$ ($m>3$).
\medskip
{\bf (*) Ranks 5 and 6.} There exists 4 shift-stable squarefree, of degree $\le 39$ polynomials $P_1,\dots,P_4$ of rank 5 and 6 whose coefficients and $L$-functions are given in the following table ($r_{nt}$ is the non-trivial part of the rank, see 4.6). 
\settabs 20 \columns
\medskip
{\bf  Table 6.11.} Numerical data for $P_1,\dots,P_4$ --- polynomials of rank 5, 6.
\medskip
\+ & $m$ & $a_m$&$r$ &$r_{nt}$&$r_\infty$&& $ c_{s,m/3}, \dots, c_{s,0}$&&&&&$k=\frac{m+1}2$ --- the size of $\goth M$ \cr
\medskip
\+ $P_ 1$&21&1&5&5&2&\phantom{2,0,0,0,0,0,}1,0,1,0,1,0,2,0&&&&&& &11&\cr
\medskip
\+ $P_ 2$&39&2&6&5&8&2,0,0,0,0,0,1,0,1,0,0,0,2,0&&&&&&&20& \cr
\medskip
\+ $P_ 3$&39&2&6&5&10&2,0,0,0,2,0,2,0,0,0,2,0,2,0&&&&&&&20& \cr
\medskip
\+ $P_ 4$&39&2&5&4&8&2,0,2,0,1,0,0,0,1,0,2,0,2,0&&&&&&&20& \cr
\medskip
\+&&&$B_k,\dots,B_0$ of (6.2a) &&&&&&&&&&&{\bf Remark.}\cr
\medskip
\+ $P_1$ & \phantom{2,0,0,0,0,0,0,0,0,}1,2,0,2,0,2,2,0,0,0,0,0 &&&&&&&&&& A priory, $B_i\in \n F_3[t]$, and not to $\n F_3$. \cr
\medskip
\+ $P_2$ & 1,2,0,1,2,2,2,1,2,2,1,1,0,0,1,0,0,0,0,0,0&&&&&&&&&&See 6.13.1, 6.16 for a discussion.\cr
\medskip
\+ $P_3$ &1,2,0,0,2,1,0,0,0,1,2,0,0,2,1,0,0,0,0,0,0&\cr
\medskip
\+ $P_4$ & 1,2,0,2,2,2,1,1,2,2,1,1,1,0,1,2,0,0,0,0,0&\cr
\medskip
\+ &&Factorization of $L(\goth C_{P_i},t,T)$&&&&&&&&&&Factorization of $\sum_{i=0}^k B_iV^i$\cr
\medskip
\+ $P_1$ &&$(T-1)^5(T+1)^2(2T^2+2)$&&&&&&&&& $V^5 (1+V)^2(2 + V)^2 (2 + 2 V +
V^2)$ \cr
\medskip
\+ $P_2$ &&$(T-1)^6(T+1)^3(2T^3+2T^2+1)$&&&&&&&&& $V^6 (1 + V)^8 (2 + V)^3 (2 + 2 V + V^3)$ \cr
\medskip
\+ $P_3$ &&$(T-1)^6(T+1)^2(2T^2+T+1)$&&&&&&&&& $V^6 (1 + V)^{10} (2 + V)^2 (1 + V^2)$ \cr
\medskip
\+ $P_4$ &&$-(T-1)^5(T+1)^3(2T^2+T+1)^2$&&&&&&&&& $V^5 (1 + V)^8 (2 + V)^3 (1 + V^2)^2$ \cr
\medskip
See 6.13, 6.16 below for the comments to this table. 
\medskip
{\bf Table 6.12.} Ratios $\goth q_{s}(m,a,r)/\goth q_{s}(m-6,a,r)$.
\settabs 10 \columns
\medskip
\+$(a,r)$& $\ m=$ & 21 & 27&33&39&&24&30&36\cr
\medskip
\+(1,1)&&5.694&4.224&4.158&4.423&&4.217&5.206&4.390\cr
\medskip
\+(1,2)&&&3.067&1.587&2.726&&&3&2.958\cr
\medskip
\+(2,2)&&&4.086&4.597&4.373\cr
\medskip
\+(2,3)&&&2&2.389& 3.674\cr
\medskip
The na\"{\i}ve parameter count (analogous to the one of 6.8) predicts much smaller values of $\goth q_{s}(m,a,r)$. We denote by ${X_1}_s(3,1,m,a,r)$ the affine variety of shift-stable polynomials in variables $c_{s,i}$, $i=0,\dots, m/3-1$ such that $a_m=a$, and the rank at 1 is $\ge r$. For the case $m_{s}$ odd, $a_m=2$ we have: 
\medskip
1. The set of $P$ depends on $m_{s}$ parameters $c_{s,0},\dots,c_{s,m_{s}-1}$;
\medskip
2. $B_0$ has degree $k=\frac{m-1}2$, hence it depends on $[\frac{k}3]=\frac{m_{s}+1}2$ parameters, and the expected dimension of ${X_1}_{s}(3,1,m,2,2)$ for $m_{s}$ odd is $\frac{m-3}6$. 
\medskip
3. $B_1$ has degree $k-1$, hence it depends on the same $[\frac{k-1}3]=\frac{m_{s}+1}2$ parameters, and the expected dimension of ${X_1}_{s}(3,1,m,2,3)$ is negative. 
\medskip
 For $m_{s}$ odd, $a_m=1$ we have the same deg($B_i$), hence the expected dimension of ${X_1}_{s}(m,1,1)$ is $\frac{m_{s}+1}2$ and the one of ${X_1}_{s}(3,1,m,1,r)$ is negative for $r\ge2$. For $m_{s}$ even we have $k=m/2$, hence $B_0$ depends on $\frac{m_{s}}2+1$ parameters, $B_1$ depends on $\frac{m_{s}}2$ parameters, and the expected dimensions of ${X_1}_{s}(3,1,m,a,r)$ are $\frac{m}6-1$ for $r=1$ and negative for $r\ge 2$. 
\medskip
Table 6.10 gives evidence that ${X_1}_{s}(3,1,m,a,r)$ is far to be the complete intersection of the polynomials $H_{s,ij}$.
\medskip
{\bf 6.13. $\tau_c\circ \nu_{c^{-1}}$-stability, and comments to 6.11.}
\medskip
2.4.2, 2.4.4 imply that $P$ and $\tau_c\circ \nu_{c^{-1}}(P)$ have the same rank. We see that shift-stable ( = $\mu_c$-stable) polynomials give us exemples of jump of rank; the same occurs for $\tau_c\circ \nu_{c^{-1}}$-stable polynomials $P=\sum_i a_{i}\theta^{(q-1)i}$. Really, all $P_1,\dots,P_4$ of 6.11 are $\tau_c\circ \nu_{c^{-1}}$-stable, as well as one of 3 elements of ${X_1}_{s}(3,1,27,1,4)$ and four of 12 elements of ${X_1}_{s}(3,1,39,1,3)$. Other observations related to high-rank polynomials are the following: 
\medskip
{\bf 1. Property $B_i\in \n F_q$.} A priory, $B_i\in \n F_q[t]$, namely $\deg_t B_i=k-i$. We see that for all $P_1,\dots,P_4$ of 6.11 we have: all $B_i\in \n F_3$ (or, equivalently, $H_{s,ij}=0$ unless $j=0$). Is it typical for other $P$ of high rank? The same phenomenon occurs for all elements of ${X_1}_{s}(3,1,21,1,4)$ and ${X_1}_{s}(3,1,27,1,4)$, but not for ${X_1}_{s}(3,1,33,1,4)$, see {\bf 3} below. Subsection 6.16 gives some evidence in favour of this phenomenon.
\medskip
{\bf 2. Factorization of  $L(\goth C_P,T)$.} For $P_1,\dots,P_4$ we have a "nice" factorization of their $\sum_{i=0}^k B_iV^i$, which is practically the same as the factorization of $L(\goth C_{P_i},T)$ (the factorization of a "random" polynomial in $\n F_q[V]$ is much "worse"). Why? Does this factorization come from a natural partition of the set of local $\goth P$-factors of $L$, i.e. is it possible to find a natural partition of the set of local $\goth P$-factors such that the product of the factors in any set of this partition will give us a factor of $L(\goth C_P,T)$? 
\medskip
{\bf 3.} There exists $P_{33,1,4}\in {X_1}_{s}(3,1,33,1,4)$ (a representative of the only $\tau_c\circ \nu_{c^{-1}}$-orbit in ${X_1}_{s}(3,1,33,1,4)$ ) --- a polynomial of sufficiently high rank 4 --- whose $B_i\not\in \n F_3$ and whose $\sum_{i=0}^k B_iV^i$ does not have a "nice" factorization: 
\medskip
\noindent
$\sum_{i=0}^k B_iV^i=V^4(1+V)^6(2+V)\cdot [V(1+V)^2(2+V^2+V^3)+2(t^3-t)(2+2V+V^3)]$
\medskip
What is the relation between the properties of 
\medskip
(a) High rank of $P$;
\medskip
(b) $B_i\in \n F_q$;
\medskip
(c) $L(\goth C_P,T)$ has a factorization to factors of small degree? 
\medskip
{\bf 4. Coincidence of $L$-functions.} The set ${X_1}_{s}(3,1,21,1,4)$ contains one $\tau_c\circ \nu_{c^{-1}}$-stable polynomial $P_1$ of rank 5 and two $\tau_c\circ \nu_{c^{-1}}$-orbits. $L$-functions of polynomials of these two orbits coincide; we do not see a reason of this coincidence. Analogously, the set ${X_1}_{s}(3,1,27,1,4)$ contains one $\tau_c\circ \nu_{c^{-1}}$-stable polynomial and one $\tau_c\circ \nu_{c^{-1}}$-orbit; their $L$-functions coincide. 
\medskip
{\bf 6.14. Rank at infinity.} To get evidence that Conjecture 8.7 is true, we present results of calculations of orders of $X_\infty(3,1,m,i)(\n F_3)$ for $k=6$ ($m=13$),  $k=7$ ($m=15$),  $k=19$ ($m=39$).  More exactly, for $m=13, \ 15$ the below table gives $\goth q_{sf\infty}(m,i)$ --- the quantity of squarefree monic $P$ of degree $m=13$ and $15$ belonging to $X_\infty(3,1,m,i)(\n F_3)$ (for $a_m=2$ the quantities are obviously the same as for $a_m=1$). For $m=39$ we use the Monte Carlo method: $r_\infty$ was calculated for about 100 000 random points in $\n F_3^{39}$ --- the set of monic polynomials of degree 39 (there was no selection of squarefree polynomials), and the table gives approximate values of $\goth q_\infty(m,i)$ --- the quantity of all monic $P$ of degree $m=39$ whose $r_\infty \ge i$.

 For $m=13$ the numbers  $\goth q_{sf\infty}(m,i)$ have remarcable factorization. 
\medskip
{\bf Table 6.15.} Numbers $\goth q_{sf\infty}(m,i)$, $\goth q_\infty(m,i)$
\medskip
\+ $m=$&13 (squarefree)&&&15 (squarefree)&&&39 (all)\cr
\medskip
\+$i=0$&$2\cdot 3^{12}=1062882$ &&&$2\cdot 3^{14}=9565938$&&&$3^{39}$ \cr
\medskip
\+$i=1$&$2\cdot 3^{5}\cdot 5 \cdot 7^2$=119070&&&$2^2\cdot3\cdot88771$=1065252 &&&$(1.90\pm0.03)\cdot 3^{37}$\cr
\medskip
\+$i=2$ &$2\cdot 3^{5}\cdot 5 \cdot 7$=17010&&&$2\cdot 3^{2}\cdot 5^3 \cdot 71$=159750 &&& $(3.64\pm0.06)\cdot 3^{35}$ \cr
\medskip
\+$i=3$&$2\cdot 3^{2}\cdot 5 \cdot 7$=630&&&$2\cdot3^2\cdot431$=7758 &&&$(6.1\pm0.2)\cdot 3^{33}$\cr
\medskip
\+$i=4$ &0&&&$2\cdot3^2$=18&&&$(9\pm0.8)\cdot 3^{31}$\cr
\medskip
\+$i=5$ &0&&&0&&&$(14\pm3)\cdot 3^{29}$\cr
\medskip
\medskip
{\bf 6.16. Correlation between $r_1(P)$ and $r_\infty(P)$.} Theorem 8.6 means that polynomials $H_{ij\infty}$ ($i=0,\dots,k-1$, $j=0,\dots, k-i$, see (8.3.1) ) are highly dependent. Condition $r_1(P)\ge c$ can be written in the form that some linear combinations of $H_{ij\infty}$ are 0 (for example, $r_1(P)\ge 1 \iff \forall j=0,\dots,k$ we have $ \sum_{i=0}^{k-j}H_{ij\infty}=0$), and the condition $B_i\in \n F_q$ from 6.13.1 is equivalent to (the reduction of) the condition $A_i=$const $\iff H_{ij\infty}=0$ for $j>0$, where $A_i$ are from (6.2). So, it is not too surprising that the condition $B_i\in \n F_q$ is often fulfilled, and there exists a correlation between $r_1(P)$ and $r_\infty(P)$.

A numerical example of this correlation is the following. We considered the set $S$ of squarefree shift-stable $P=\sum_{i=0}^{13} c_{s,i} (\theta^q-\theta)^i$ of degree 39 having $c_{s,13}=2$, $c_{s,12}=1$ and having $r_1\ge 3$ ( $\iff r_{1,nt}\ge 2$). We have $\#(S)=46$. The quantity of polynomials in $S$ having the given value of $r_\infty$ is given in the following table: 
\settabs 20 \columns
\medskip
\+{\bf Table 6.17.} Value of $r_\infty$ &&&&&&&0&1&2&3&4&5&6&7&8&9&10&$\ge11$\cr
\medskip
\+ $\#$ of polynomials $\in S$\cr
\+ having this $r_\infty$ &&&&&&&0&0&15&0&0&8&11&9&1&1&1&0\cr
\medskip
Table 6.15 shows that in the absence of correlation the lower line of Table 6.17 should have the form $\sim$ 40, 5, 1, 0, ..., 0. 
\medskip
Another argument in the favour of correlation is the fact that polynomials $P_2$, $P_3$, $P_4$ of Table 6.11 have very high values of $r_\infty$.
\medskip
Exact statements (although as conjectures) on dimensions and degrees of the varieties of $(a_0,\dots,a_m)$ such that  $A_i$ from (6.2) are constants, on correlation between $r$ and $r_\infty$, etc., are not known yet. 
\medskip
{\bf 6.18. Case of any $n$.} If $n>1$ then the same dimension considerations give us much smaller values of the expected dimension of $X_1(q,n,m,r)$. For example, for $n=2$, $q=3$ we cannot expect to get $r\ge 3$, and we can expect $r=2$ only for the case $m$ even, $a_m=1$ (the distinguished coset). A computation is concordant with this prediction: for $m=8,10,12$ we have respectively $\goth q(m,1,2)=9,21, 81$, all other $\goth q(m,a,2)$ for $m\le 12$ are 0, there is no $P$ of degree $\le12$ such that $r(\goth C_P^2)\ge 3$. 
\medskip
{\bf 7. Results obtained without application of the Lefschetz trace formula.}
\medskip
We shall prove in this section that it is possible to prove without application of the Lefschetz trace formula the two following propositions: 
\medskip
1. $L(\goth C^n_P,t,T)\in\n F_q[t][T]$ (recall that a priory $L(\goth C^n_P,t,T)\in\n F_q[t][[T]]$ );
\medskip
2. Proposition 4.4 on a coset of $r\ge 1$. 
\medskip
Let $\goth P=\sum_{i=0}^d c_i\theta^i\in \n F_q[\theta]$ be an irreducible monic polynomial of degree $d$. We let $\goth P[t]:=\sum_{i=0}^d c_it^i\in \n F_q[t]$. 
The $\goth P$-local factor $L_\goth P(\goth C^n_P,t,T)$ of $L(\goth C^n_P,t,T)$ is $(1-[\tilde P(t-\tilde \theta)^n] ^{[d]}T^d)^{-1}$. Since the map $*\mapsto * ^{[d]}$ is multiplicative, we have $[\tilde P(t-\tilde \theta)^n]^{[d]}=\tilde P^{[d]}((t-\tilde \theta)^{[d]})^n$. Obviously  $(t-\tilde \theta)^{[d]}=\goth P[t]$. Further, we denote the roots of $P=\sum_{i=0}^m a_i\theta^i$ in $\bar \n F_q$ by $r_1, \dots, r_m$ (with multiplicities). We have $\tilde P^{[d]}=((-1)^ma_m)^d \cdot \goth P(r_1) \cdot\dots\cdot\goth P(r_m)$, hence $$L_\goth P(\goth C^n_P,t,T)=(1-((-1)^ma_m)^d \cdot\goth P(r_1) \cdot\dots\cdot\goth P(r_m)\cdot\goth P[t]^n\cdot T^d)^{-1}$$ This expression is multiplicative, therefore we can apply the converse of the Euler product formula: 
$$L(\goth C^n_P,t,T)=\sum_{\goth N} ((-1)^ma_m)^d \cdot\goth N(r_1) \cdot\dots\cdot\goth N(r_m)\cdot\goth N^n\cdot T^d$$ $$ =\sum_ {d=0}^\infty ((-1)^ma_m)^d \ \  [\ \sum_{\goth N} \goth N(r_1) \cdot\dots\cdot\goth N(r_m)\cdot\goth N^n] \ \  T^d $$ where the sum runs over all monic $\goth N\in F_q[t]$ and $d\ge 0$ is the degree of $\goth N$. 

We see that for $P=1$ we have: $$L(\goth C^n,t,T)=\zeta(-n,T)\eqno{(7.1)}$$ where $\zeta$ is the Goss' zeta function of the ring $\n F_q[t]$ (see, for example, [T], p. 233, middle of the page; $T$ of the present paper is $X$ of [T]). 
\medskip
{\bf 7.2. Second proof of $L(\goth C^n_P,t,T)\in\n F_q[t][T]$.} We must prove that for $d>>0$ $\sum_{\goth N} \goth N(r_1) \cdot\dots\cdot\goth N(r_m)\cdot\goth N^n=0$, where $\goth N$ is monic of degree $d$. Let $s$ be a number such that all $r_1, \dots, r_m\in \n F_{q^s}$. For $b_1,\dots,b_m\in \n F_{q^s}$ we denote by $W_{d,\beta}$ (where $\beta = (b_1,\dots,b_m)\in \n (F_{q^s})^m$) the set of all monic $\goth N\in\n F_q[t]$ of degree $d$ such that $\goth N(r_i)=b_i$, $i=1,\dots,m$. We have $$\sum_{\goth N} \goth N(r_1) \cdot\dots\cdot\goth N(r_m)\cdot\goth N^n=\sum_{\beta\in \n (F_{q^s})^m} b_1 \cdot\dots\cdot b_m\sum_{\goth N\in W_{d,\beta}} \goth N^n$$ Clearly $W_{d,\beta}$ is a $\n F_q$-affine space in $\n F_q[t]$, possibly empty. Now we can apply a Goss' lemma ([T], p. 234, Theorem 1):
\medskip
Let $W$ be a $\n F_q$-affine space in a ring over $\n F_q$. If dim $W> n$ then $\sum _{\goth N \in W}\goth N^n=0$. 
\medskip
{\bf Remark.} The statement of [T], p. 234, Theorem 1 gives a more low bound for dim $W$ with respect to $n$, we do not need it. Further, the same statement requires $0\not\in W$, but really this condition is excessive. 
\medskip
Hence, in order to prove our result, we must prove that for all $\beta\in \n (F_{q^s})^m$ either $W_{d,\beta}=\emptyset$ for all $d$, or dim $W_{d,\beta}\to\infty$ as $d\to\infty$. This is obvious: if $\goth N_0\in W_{d,\beta}$ for some $d$, $\beta$, then --- since $P(r_i)=0$ --- for any $Y=\sum y_it^i\in \n F_q[t]$ of degree $m'>>0$ with the leading coefficient $y_{m'}=a_m^{-1}$ we have $\goth N=\goth N_0+ P[t]Y\in W_{m+m',\beta}$. $\square$
\medskip
{\bf Remark.} Using the exact bound of dim $W$ such that $\sum _{\goth N \in W}\goth N^n=0$ in Goss' lemma, the reader can try to find the upper bound of the degree of $L(\goth C^n_P,t,T)$ as a polynomial in $T$. Will be got the same value that is given by Proposition 3.3? 
\medskip
{\bf 7.3. Second proof of Proposition 4.4.} Condition $m\equiv -n$ mod $q-1$ implies $(-1)^m=(-1)^n$ in $\n F_q$, hence $a_m=(-1)^n$ implies $((-1)^ma_m)^d=1$. Further, condition $m\equiv -n$ mod $q-1$ implies that for any $c\in \n F_q^*$ we have $$\goth N(r_1) \cdot\dots\cdot\goth N(r_m)\cdot\goth N^n=c\goth N(r_1) \cdot\dots\cdot c\goth N(r_m)\cdot(c\goth N)^n$$
hence $$L(\goth C^n_P,t,1)=-\sum_{\goth N} \goth N(r_1) \cdot\dots\cdot\goth N(r_m)\cdot\goth N^n$$ where the sum runs over the set of all $\goth N\in F_q[t]$. More exactly, $L(\goth C^n_P,t,1)$ is the stable value as $d\to\infty$ of $$L(\goth C^n_P,t,1)_d:=-\sum_{\goth N\in Pol_{\le d}} \goth N(r_1) \cdot\dots\cdot\goth N(r_m)\cdot\goth N^n$$ where $Pol_{\le d}$ is the set of all elements of $\n F_q[t]$ of degree $\le d$. 

The same arguments as the ones of the proof of (7.2) (sets $W_{d,\beta}$ now are the sets of all $\goth N$ of degree $\le d$ such that $\goth N(r_i)=b_i$) show that for $d>>0$ we have $L(\goth C^n_P,t,1)_d=0$. $\square$
\medskip
{\bf Part II. Resultantal varieties (characteristic 0 case). } 
\medskip
{\bf Introduction.} In order to make the present part independent on the Part I, we repeat some definitions. The object of research are varieties $X_c(q,n,m,i)$ defined over $\n Z$ where $q\ge2$, $n\ge0$, $m\ge1$, $i\ge1$ are integer parameters and $c\in(\n Z[t]\cup\infty)$ is a polynomial parameter. $X_c(q,n,m,i)$ are sets of zeroes of some explicitly defined polynomials $H_*$ coming from the determinant of a matrix. It turns out that for $c=\infty$ the polynomials $H_*$ are highly dependent, and for most valies of $q,n,m,i$ finding the dimension of $X_\infty(q,n,m,i)$ is a non-trivial problem. For other $c\in\n Z[t]$ we can only expect that $H_*$ are dependent. 
\medskip
Importance of study of $X_c(q,n,m,i)$ comes from the results of Part I. Namely, let $q=p^\ga$ be a power of a prime $p$. We denote by tilde the reduction at $p$. We consider the $n$-th tensor power of the Carlitz module over $\n F_q$ and its twists by polynomials of degree $\le m$. It turns out that $\widetilde{X}_c(q,n,m,i)(\n F_q)$ describes the set of twists such that the order of zero of their $L$-functions at the point $\tilde c$ is $\ge i$. Particularly, we have the following reformulaton of the characteristic $p$ analog of a famous conjecture of non-boundedness of ranks of twists of any elliptic curve over $\n Q$:
\medskip
{\bf II.1. Reformulation.} Let $q,n$ be fixed. If $\forall \ i \ \exists \ m $ such that $\widetilde{X}_1(q,n,m,i)(\n F_q)$ $\ne\emptyset$ then the order of zero at 1 of the $L$-functions of twists of the $n$-th tensor power of the Carlitz module over $\n F_q$ is not bounded. 
\medskip
As an approach to solve this problem, we should answer
\medskip
{\bf II.2. Question.} Let  $q,n,i$ be fixed. Whether $\exists \ m $ such that $X_1(q,n,m,i)\ne\emptyset$, or not? Moreover, whether dim $X_1(q,n,m,i)\to\infty$ as $m\to\infty$, or not? 
\medskip
If $H_*$ were independent then the answer to this question would be negative, but we expect that they are dependent. See Remark 8.2.4 --- a discussion of a particular case of $X_1(3,1,15,1)$. 
\medskip
Part II is organized as follows. Section 8 contains results that are valid for any $q$, and in Section 9 we get much more detailed results for $q=2$. The definition of $X_c(q,n,m,i)$ is given in (8.1.1). Later we shall consider exclusively the case $c=\infty$, and we omit the index $\infty$. We repeat the definition of $X(q,n,m,i)$ for this case.
\medskip
Theorem 8.6 shows that for $c=\infty$ the polynomials $H_*$ are highly dependent. In 8.7 - 8.9 we state problems on dimension of $X(q,n,m,i)$ and of its irreducible components; this is a subject of further research. Since 8.6a states that $X(q,n,m,i)$ has a subvariety $X_r(q,m,i)$ of maximal (?) dimension and $X_r(q,m,i)$ does not depend on $n$, and since for $q=2$ the conjecture 9.3 (it has a strong numerical evidence) claims that $\forall \ n\ge0$ \ $X(2,n,m,i)$ do not depend on $n$, we state an open question 8.8: whether for all $q$, for all sufficiently large $n$ the varieties $X(q,n,m,i)$ do not depend on $n$, or not? 
\medskip
We see that for $q=2$ the subject of further research is: 
\medskip
{\bf (II.3)} To prove Conjecture 9.3; 
\medskip
{\bf (II.4)} To study varieties $X(2,0,m,i)$. 
\medskip
Let us consider both these problems.
\medskip
An important difference between the cases $q=2$ and $q>2$ is the following: for $q=2$ we have $X(2,n,m,i)$ are (conjecturally) equal for all $n$, including $n=0$, while for $q>2$ they are (maybe) equal only for sufficiently large $n$. Varieties $X(q,0,m,i)$ are much simpler objects than $X(q,n,m,i)$, because the set of $H_*$ defining them, depends on 1 parameter and $X(q,0,m,i)$ are (conjecturally) complete intersections, while for a general $X(q,n,m,i)$ the set of $H_*=H_{\al\be}$ defining them, depends on 2 parameters $\al,\be$, see below. 
\medskip
The set of $H_*$ defining $X(2,0,m,i)$ is denoted by $D(m,0), D(m,1),\dots,D(m,i-1)$, and the set of $H_*$ defining $X(2,n,m,i)$ is denoted by $H_{\al\be n}$ where $\al,\be$ run over a $\n Z$-triangle $\Delta$ (see (8.3.2), (9.12) --- (9.15) for details). To prove (9.3) we must prove that $\forall \  n,m,i$ $$\forall \ \al,\be\ \exists \ \ga\suchthat H_{\al\be n}^\ga\in \langle D(m,0),\dots, D(m,i-1)\rangle \eqno{(\hbox{II}.5)}$$ where $\langle D(m,0),\dots, D(m,i-1)\rangle $ is the ideal generated by $D(m,*)$. 
\medskip
{\bf II.6.} At the moment (II.5) is proved / conjectured / verified for the cases: 
\medskip
(1) $i=1$ ( $\iff \al=0$), any $m,n$ --- "base of $\Delta$". Proved in Proposition 9.12 or Part III. 
\medskip
(2) $\be=0$, $\be=$ max, any $m,n$ --- "the lateral sides of $\Delta$". Proved in Proposition 9.14. 
\medskip
(3) Some "interior points of $\Delta$ near vertices". Conjecture based on calculation, without proof. Remark 9.15. 
\medskip
In all these cases $\ga$ of (II.5) are 1. 
\medskip
(4) Case $n=1$, $m=4$, $i=2$. Explicit calculation 9.17b. In this case $\ga=2$. 
\medskip
Problem (II.4). Conjectures on this problem are given in 9.7. They describe the quanity of irreducible components of $X(2,0,m,i)$, their degrees, multiplicities etc. At the moment we have no even conjectural values of most of these numerical characteristics, and the proofs look a much more difficult problem. Moreover, although $X(2,0,m,i)$ and $X(2,n,m,i)$ coincide as the sets of points, they are different as schemes, particularly, the multiplicities of their irreducible components are different. See Conjectures 9.7.8, 9.7.10. 
\medskip
{\bf 8. Case of any $q$.}
\medskip
Let $q\ge2$, $n\ge0$, $m\ge1$ be integers such that $k:=\frac{m+n}{q-1}-1$ is integer $\ge1$. Let $P=\sum_{i=0}^m a_i\theta^i\in\n Z[a_0,\dots,a_m][\theta]$ be a polynomial.
\medskip
The $k\times k$-matrix $\goth M(P,n,k)$ is defined by the formula (3.1). For the reader's convenience, its explicit form is given for several particular cases: case $n=1$ in (3.1), case $q=2$, $n=0$, $m=6,7$ in (9.1), case $q=2$, $n=2$, any $m$ here:  
$$\goth M(P,2,k)^t=$$ $$\left(\matrix a_{1}t^2-2a_{0}t &  a_{3}t^2-2a_{2}t+a_{1}&  a_{5}t^2-2a_{4}t+a_{3}  &   \dots & 0 \\ a_{0}t^2 &  a_{2}t^2-2a_{1}t+a_{0}  &  a_{4}t^2-2a_{3}t+a_{2}&   \dots & 0 \\  0 &  a_{1}t^2-2a_{0}t &  a_{3}t^2-2a_{2}t+a_{1}&   \dots & 0 \\ \dots&
\dots & \dots & \dots  & \dots \\ 0 &  0  &   \dots &\dots & -2a_{m}t+a_{m-1} \endmatrix \right) $$ 
In Part I the elements $a_0,\dots,a_m$ were considered as elements of $\n F_q$, where $q$ was a power of a prime; now we consider them as abstract elements or elements of $\bar \n Q$. Entries of $\goth M(P,n,k)$ belong to $\n Z[a_0,\dots,a_m][t]$. Let $$CH(\goth M(P,n,k),T):=\det(I_k - \goth M(P,n,k)T)\in\n Z[a_0,\dots, a_m][t][T]\eqno{(8.1)}$$ 
be (a version of)  its characteristic polynomial. Let $c\in \n Z[t]$ be fixed (for the most interesting case $c=\infty$ see (8.3)). 
\medskip
{\bf Definition 8.1.1.} $X_c(q,n,m,i)\subset A^{m+1}$ is an affine algebraic variety defined by the condition that the order of 0 of $CH(\goth M(P,n,k),T)$ at $T=c$ is $\ge i$. 
\medskip
{\bf Remark 8.1.2.} If $a_0,\dots,a_m\in \n F_q$ then $CH(\goth M(P,n,k),T)$ is the non-trivial factor of $L(\goth C^n_P,T)$ (see 4.6). 
\medskip
We have:  
$$CH(\goth M(P,n,k),T)=\sum_{i=0}^k A_iT^{k-i}\eqno{(8.2.1)}$$
where $A_i\in \n Z[a_0,\dots, a_m][t]$, $\deg_t A_i\le n(k-i)$ (because entries of $\goth M$ are of degree $\le n$). After the substitution $T=c+W$ (8.2.1) becomes 
$$CH(\goth M(P,n,k),T)=\sum_{i=0}^k C_i W^{i}\eqno{(8.2.2)}$$
and 
$$C_i=\sum_{j=0}^{*}H_{ij,c}t^j\eqno{(8.2.3)}$$ (we are not interested to know $\deg_t(C_i)$, because we shall not work with them) and $H_{ij,c}\in \n Z[a_0,\dots,a_m]$. 

This means that $X_c(q,n,m,i)$ is really an affine variety, it is the set of zeroes of $\{H_{\al j,c}\}$, $0\le\al<i$.
\medskip
{\bf Remark.} In Section 6 for a particular case $c=1$ we use a version of $H_{\al j,1}$ which differ from $H_*$ defined above (but clearly generate the same ideal), because $H_{\al j,1}$ of Section 6 are of lesser degree. Since we shall not more work with the case $c\ne\infty$, there will be no confusion. 
\medskip
{\bf Remark 8.2.4.} It would be interesting to find dim $X_c(q,n,m,i)$ in some simplest cases. For example, $X_1(3,1,m,1)$ is the set of zeroes of $k+1=(m+1)/2$ polynomials in $A^{m+1}$. For example, for $m=15$ we have 8 polynomials and 16 variables. From another side, $\#(\tilde X_1(3,1,15,1)(\n F_3))\approx 3^{10}$ (see Table 6.7) which gives evidence that dim $X_1(3,1,15,1)=10$. 
\medskip
{\bf 8.3. Case $c=\infty$.} We get from (8.2.1): 
$$A_i=\sum_{j=0}^{n(k-i)}H_{ij,n}t^j\eqno{(8.3.1)}$$
(here and until the end of the paper we omit the index $c=\infty$, but we write the index $n$ if necessary). $H_{ij,n}$ are homogeneous of degree $k-i$, hence 
$$CH(\goth M(P,n,k), T)=(H_{00}+H_{01}t+\dots +H_{0,nk}t^{nk})T^k+$$ $$+(H_{10}+H_{11}t+\dots +H_{1,n(k-1)}t^{n(k-1)})T^{k-1}+ ... +(H_{k-1,0}+\dots +H_{k-1,n}t^n)T+H_{k0}\eqno{(8.3.2)}$$

{\bf 8.3.3.} There is a symmetry of order 2: $a_i\mapsto a_{m-i}$ sends $H_{ij}$ to $\pm H_{i,n(k-i)-j}$.
\medskip
{\bf Remark 8.3.4.} If $n=0$ then the formulas (8.3.1) are much simpler, because there is no variable $t$. 
\medskip
Let us repeat the definition of $X(q,n,m,i)$ (recall that for $c=\infty$ the number $i$ is the complement of the order of zero of $CH(\goth M(P,n,k), T)$ at $T=\infty$): 
\medskip
{\bf 8.4. Definition.} $X(q,n,m,i)$ is a projective subvariety of $P^m(\bar \n Q)=\{(a_0:...: a_m)\}$ defined by the condition 
$$\deg_T(CH(\goth M(a_0,\dots, a_m),T))\le k-i$$
\noindent
Particularly, $X(q,n,m,0)=P^m(\bar \n Q)$, 
\medskip
\noindent
$X(q,n,m,1)$ is the set of zeroes of $H_{00}, H_{01}, \dots, H_{0,nk}$ in $P^m(\bar \n Q)$; 
\medskip
\noindent
$X(q,n,m,2)$ is the set of zeroes of $H_{00}, \dots, H_{0,nk}$, $H_{10},\dots, H_{1,n(k-1)}$  in $P^m(\bar \n Q)$;
\medskip
\noindent
$X(q,n,m,i)$ is the set of zeroes of in $P^m(\bar \n Q)$ of $$H_{00}, \dots, H_{0,nk}, \ \  H_{10}, \dots, H_{1,n(k-1)}, \ \ \dots, \ \ H_{i-1,0}, \dots, H_{i-1,n(k-(i-1))} \eqno{(8.5)}$$ 
\medskip
\noindent
and $P^m(\bar \n Q)=X(q,n,m,0)\supset X(q,n,m,1)\supset X(q,n,m,2)\supset \dots \supset X(q,n,m,k)\supset X(q,n,m,k+1)=\emptyset$ (because $H_{k0}=1$). 
\medskip
Sometimes we shall use the same notation $X(q,n,m,i)$ for the scheme 
$$\Proj \bar \n Q[a_0,\dots, a_m]/\{H_{00}, \dots, H_{0,nk}, \ H_{10}, \dots, H_{1,n(k-1)}, \ \dots, $$ $$\ H_{i-1,0}, \dots, H_{i-1,n(k-(i-1))}\}$$
Equivalently, $r_\infty$ (see 4.1) is the maximal $i$ such that $a_0,\dots, a_m\in X(q,n,m,i)$. 
\medskip
It turns out that polynomials $H_{ij}$ are highly dependent. We have:
\medskip
{\bf Theorem 8.6.}  (a) $X(q,n,m,i)$ has a subvariety $X_r(q,m,i)$ (not depending on $n$) of codimension $ i(q-1)$ in $P^m$.
\medskip
(b) $X(q,n,m,i)$ has a subvariety $X_l(q,n,m,i)$ of codimension $\le i(q-1)+(q-2)(n-1)$ in $P^m$ (for $q=2$ we have $X_r(q,m,i)= X_l(q,n,m,i)$, see 9.9). 
\medskip
(c) $X(q,n,m,i)$ is of codimension $\le i(i+n)$ in $P^m$. 
\medskip
{\bf Conjecture 8.7.} Codim $X(q,n,m,i)=$ min ( $ i(q-1), \  i(i+n)$ ), i.e. it is equal to $ i(q-1)$ for $i\ge q-1-n$, \ $ i(i+n)$ for $i\le q-1-n$.
\medskip
{\bf Open question 8.8.} Is it true that $\exists \ n_0=n_0(q)$ such that for $n\ge n_0$ we have: $X(q,n,m,i)$ does not depend on $n$? 
\medskip
For $q=2$ we have evidence that this is true, see 9.3, and $n_0(2)=0$. For $q>2$ the truth of (8.8) is supported by (8.6a, b). 
\medskip
{\bf Conjecture 8.9.} If $n=1$, $q>2$, $i\ge q-2$ then $X_r(q,m,i)\ne X_l(q,n,m,i)$, they are different components of $X(q,n,m,i)$ of maximal dimension.
\medskip
Let us prove Theorem 8.6 and other results. We denote by
$\goth N(P,n,k)$ the following $(k+n)\times k$-matrix (we apologize that it is transposed with respect to $\goth M$) $$\goth N(P,n,k):=\left(\matrix a_{q-1} & a_{2q-1}& a_{3q-1}&...&a_{kq-1}\\ a_{q-2} & a_{2q-2}& a_{3q-2}&...&a_{kq-2}\\ a_{q-3} & a_{2q-3}& a_{3q-3}&...&a_{kq-3}\\ \dots & \dots & \dots & \dots & \dots \\  a_{q-(k+n)} & a_{2q-(k+n)}& a_{3q-(k+n)}&...&a_{kq-(k+n)}\endmatrix \right)\eqno{(8.11)}$$
(here like in (3.1), $a_*=0$ if $*\not\in [0,\dots, m]$). 
\medskip
{\bf Lemma 8.12.} (a) For any $q$ and $n=1$ we have: $(a_0: ... : a_m)\in X(q,n,m,1)\iff $ the rank of $\goth N(P,n,k)$  is $<k$; 
\medskip
(b) For any $q,\ n$, if the rank of $\goth N(P,n,k)$  is $\le k-i$ then $(a_0: ... : a_m)\in X(q,n,m,i)$. 
\medskip
{\bf Proof.} (a) $(a_0: ... : a_m)\in X(q,n,m,1)\iff |\goth M(P,n,k)|=0 \iff \forall \ i \ \ H_{0i}=0$. For $n=1$ \ $H_{0i}= \pm$ (minor of $\goth N(P,n,k)$ obtained by elimination of its $(i+1)$-th row). 

(b) $H_{jl}$ is a linear combination of minors of $\goth N(P,n,k)$ of size $k-j$. $\square$
\medskip
{\bf Remark 8.13.} (a) Most likely, we can omit the restriction $n=1$ at 8.12(a), i.e. for any $n$ we have: if $(a_0: ... : a_m)\in X(q,n,m,1) $ then the rank of $\goth N(P,n,k)$  is $<k$. From another side, even for $n=1$ the condition $(a_0: ... : a_m)\in X(q,n,m,2) $ does not imply that the rank of $\goth N(P,n,k)$  is $\le k-2$. 
\medskip
(b) Property of small codimension of $X(q,n,m,i)$ does not hold for a more general matrix $\goth N$. See Appendix, 3. 
\medskip
{\bf Definition 8.14.} We denote by $X_\goth N(q,n,m,i)$ the set of $P$ such that the rank of $\goth N(P,n,k)$  is $\le k-i$. 
\medskip
Lemma 8.12 means $X_\goth N(q,n,m,i)\subset X(q,n,m,i)$.
\medskip
{\bf 8.15.} In order to formulate further results, let us consider the general $\de$-th determinantal variety $D(\de,\alpha,\ga)$ of the $\alpha\times\ga$-matrices, i.e. the subset of the space $P^{\alpha\ga-1}=\{(c_{11}:...:c_{\alpha\ga})\}$ of $\alpha\times\ga$-matrices $C_{\al\ga}=\left(\matrix c_{11}&\dots&c_{1\ga}\\\dots&\dots&\dots\\c_{\alpha1}&\dots& c_{\alpha\ga} \endmatrix \right)$ whose rank is $\le \de$. By the general theory of determinantal varieties, $$\hbox{Codim $D(\de,\alpha,\ga)$ in $P^{\alpha\ga-1}$ is $(\alpha-\de)(\ga-\de)$}\eqno{(8.15.1)}$$ and for $\de=\alpha-1<\ga$ we have $$\deg D(\alpha-1,\alpha,\ga)=\binom{\ga}{\alpha-1}\eqno{(8.15.2)}$$ (see, for example, [FP], (1.2), p. 4). 
\medskip
For the matrix $\goth N(P,n,k)$ we have $\al=k+n$, $\ga=k$, and it defines a linear inclusion $P^m\to P^{(k+n)k-1}$. We have $X_\goth N(q,n,m,i)=D(k-i,k+n,k)\cap P^m$, hence 
$$\hbox {codim $X_\goth N(q,n,m,i)$ in $P^m$ is $\le i(i+n)$} \eqno{(8.16)}$$
and if the equality holds and $n=1$ then deg $X_\goth N(q,n,m,i)=\binom{k+n}{n+1}$. The case of the equality in (8.16) will be called the trivial case. We shall see that in many cases the codimension is much less. 
\medskip
To show that $\goth N(P,n,k)$ is not of the maximal rank, it is sufficient to find a matrix such that the product of $\goth N(P,n,k)$ and this matrix is 0. Depending on the side of the product, we get two different subvarieties of $X_\goth N(q,n,m,i)$. 
\medskip
For simplicity we consider first the case when $\goth N(P,n,k)$ has a simple form (all below results are true for the general case, see 8.32). This is the case $$m\equiv -1 \mod q\eqno{(8.17)}$$ We let $\be:=\frac{m+1}q$. In this case $l:=\frac{k+n}q$ is integer. The reader should consider $q$, $m$ and hence $\be$ as constants, and $l$, $k$, $n$ as variables.
\medskip
Let $B$ be a $(q\times \be)$-block $\left(\matrix a_{q-1} & a_{2q-1}& a_{3q-1}&...&a_{m}\\ a_{q-2} & a_{2q-2}& a_{3q-2}&...&a_{m-1}\\ a_{q-3} & a_{2q-3}& a_{3q-3}&...&a_{m-3}\\ \dots & \dots & \dots & \dots & \dots \\  a_{0} & a_{q}& a_{2q}&...&a_{m-(q-1)}\endmatrix \right)$ (the upper
\medskip
\noindent
left $(q\times \be)$-submatrix of $\goth N(P,1,k)$ ).
\medskip
{\bf Definition 8.18.} Let $\goth B(P,\mu)$ be a matrix having the following block form: $$\goth B(P,\mu)=\left(\matrix B&0&\dots\\0&B&0\\ \dots & \dots & \dots\\\dots&0&B \endmatrix \right)\eqno{(8.19)}$$ of $\mu$ block rows, each block row consists of $q$ ordinary rows, the left 0-block of the $j$-th block row contains $j-1$ ordinary columns and the right  0-block of the $j$-th block row contains $\mu-j$ ordinary columns. 
\medskip
Particularly, $\goth N(P,n,k)=\goth B(P,l)$. 
\medskip
{\bf Lemma 8.20.} (a) $\forall l,j$ Corank $\goth B(P,l-j)>0\Rightarrow$ rank $\goth B(P,l)\le ql-j-1$. 
\medskip
(b) If $q=2$ and the quantity of rows of  $\goth B(P,l-j)$ is $\le$ the quantity of its columns, then $\forall l,j$ Corank $\goth B(P,l-j)>0\iff$ rank $\goth B(P,l)\le ql-j-1$. 
\medskip
{\bf Proof.} (a) If Corank $\goth B(P,l-j)$ $> 0$ then its rows are linearly dependent. Let 
$$\Lambda:=\left(\matrix \la_1&\dots&\la_{q(l-j)}\endmatrix \right)$$ be a row matrix of dependence, i.e. $\Lambda\goth B(P,l-j)=0$. We define a $(j+1)\times (ql)$-matrix $\Lambda_j$ as follows (this is a block structure): 
$$\Lambda_j:=\left(\matrix \Lambda&0&\dots\\0&\Lambda&0\\ \dots & \dots & \dots\\\dots&0&\Lambda \endmatrix \right)$$ where block rows are ordinary rows, the $\al$-th row consists (from the left to the right) of $1\times q(\al-1)$ zero block, $1\times (q(l-j))$ block equal to $\Lambda$, and a  $1\times q(j+1-\al)$ zero block. We have $\Lambda_j \goth B(P,l)=0$, hence --- because $\Lambda_j$ is of the maximal rank --- the rank of $\goth B(P,l)$ is $\le ql-j-1$. 
\medskip
(b) Let $\goth B(P,l-j)$ be of the maximal rank. It is equal to $2(l-j)$ --- the quantity of its rows. If $a_m\ne0$ we consider a submatrix of $\goth B(P,l)$ generated by a maximal non-zero minor of $\goth B(P,l-j)$ and elements $a_m$ in the last $j$ columns of $\goth B(P,l)$. It is clear that its determinant is $\ne0$ and its size is $2l-j$, hence we get the desired. If $a_m=0$ and $a_{m-1}\ne0$ then we consider elements $a_{m-1}$ instead of $a_m$. If both $a_{m-1}=a_m=0$ and $a_{m-2}\ne0$, then we consider elements $a_{m-2}$ in the right $(j+1)$-th, ... , 2-nd columns of $\goth B(P,l)$ etc. $\square$
\medskip
We let $j=i+n-1$. 
\medskip
{\bf Definition 8.21.} $X_l(q,n,m,i)$ the set of $(a_0: ... : a_m)\in P^m$ such that Corank $\goth B(P,l-j)>0$ (subscript $l$ because of the left multiplication $\Lambda_j \goth N(P,n,k)$ of $\goth N(P,n,k)$ ). 
\medskip
{\bf 8.21a.} Lemma 8.20 implies $X_l(q,n,m,i)\subset X_\goth N(q,n,m,i)\subset  X(q,n,m,i)$, and if $q=2$, $i\ge1$ then $X_l(q,n,m,i)=X_\goth N(q,n,m,i)$. 
\medskip
We apply (8.15) to the matrix $\goth B(P,l-j)$. In notations of 8.15, we have for $i>0$: $$\al=k-(q-1)n-q(i-1)\le \ga=k-n-(i-1), \ \ \ \de=\al-1$$ hence $$(\alpha-\de)(\ga-\de)= i(q-1)+(q-2)(n-1)$$ $$\binom{\ga}{\alpha-1}=\binom{k-(i+n-1)}{ i(q-1)+(q-2)(n-1)}$$ The space $P^m=\{(a_0: ... : a_m)\}$ is a linear subspace of $P^{\alpha\ga-1}$, hence we have 
\medskip
{\bf Proposition 8.22.}  Codim $X_l(q,n,m,i)$ in $P^m$ is $\le  i(q-1)+(q-2)(n-1)$.
\medskip
{\bf Conjecture 8.23.} In the above notations, $P^m=\{(a_0: ... : a_m)\}$ and $D(\de,\alpha,\ga)$ are of general intersection in $P^{\alpha\ga-1}$. 
\medskip
Hence, if 8.23 holds then codim $X_l(q,n,m,i)=  i(q-1)+(q-2)(n-1)$ and $$\deg X(q,n,m,i)= \binom{k-(i+n-1)}{ i(q-1)+(q-2)(n-1)}\eqno{(8.23.1)}$$
Let us compare 8.22 and 8.16. Since 
\medskip
$i(i+n)\underset{>}\to{\overset{<}\to{=}} i(q-1)+(q-2)(n-1)$ for $i\underset{>}\to{\overset{<}\to{=}}q-2$, we get a 
\medskip
{\bf Proposition 8.24.} (a) If $i\ge q-2$ then $X(q,n,m,i)$ has a subvariety $X_l(q,n,m,i)$ of codimension $\le i(q-1)+(q-2)(n-1)$ in $P^m$. 
\medskip
(b) If $i\le q-2$ then $X(q,n,m,i)$ is of codimension $\le i(i+n)$ in $P^m$. 
\medskip
(c) Conjecturally, for the case (a) we have equality. 
\medskip
{\bf Remark 8.25.} There exists a question: Let $q$, $\be$ (and hence $m$) and $\mu$ be fixed. For which $k$, $n$ the condition Corank $\goth B(P,\mu)>0$ defines a non-trivial $X_l(q,n,k,i)$? 
\medskip
Answer: An elementary calculation shows that we have 
$$i=-n(1-\frac1{q-1})+\frac{\be-1}{q-1}-\mu+1$$
hence if $q\ne2$ then we have only finitely many such $k$, $n$. 
\medskip
{\bf Remark 8.26.} To prove Conjecture 8.23 for $X_l(q,1,m,1)$ it is sufficient to find a linear subspace $L'\subset P^m$ such that Codim $X_l(q,1,m,1)\cap L'$ in $L'$ is $q-1$. There are many ways to construct these $L'$, we give one of them. This is elementary but tedious (see Appendix, 2). 
\medskip
{\bf Right multiplication part of $X(q,n,m,i)$.} Let $\goth x:=(x_1,\dots,x_i)$ be a set of different elements. As earlier we consider the case 8.17. Let $\goth x(\be)$ be a $\be\times i$ matrix whose $(\ga,\de)$-th entry $\goth x(\be)_{\ga\de}$ is $x_\de^{\ga-1}$. We denote by $H_\goth x$ a linear subspace of $P^m=\{(a_0: ... : a_m)\}$ defined by the condition $B\cdot\goth x(\be)=0$. The matrix of coefficients of linear equations defining $H_\goth x$ is block diagonal. Each block is a Vandermonde matrix, hence this matrix is of the maximal rank, hence the codimension of $H_\goth x$ is $qi$. For any $\goth x$ we have: $H_\goth x \subset X_\goth N(q,n,m,i)$. $H_\goth x$ does not depend on $n$. 
\medskip
{\bf Definition 8.27.} $X_r(q,m,i)$ is the closure of the union $\cup_\goth x H_\goth x $ as a subvariety of $X_\goth N(q,n,m,i)$.
\medskip
We must consider the closure, because the condition that all $x_*$ are different implies that the union $\cup_\goth x H_\goth x $ is not a closed subvariety of $X(q,n,m,i)$. The subscript $r$ because of the right multiplication $\goth N(P,n,k)\cdot \goth x(k)=0$ of $\goth N(P,n,k)$. 
\medskip
{\bf 8.28.} Let for $\ga\in \{0,\dots,q-1\}$ we denote $P_{[\ga]}:=\sum_{i=0}^{\be-1} a_{qi+\ga}x^i$ (a polynomial whose coefficients are from the $\ga+1$-th row of $B$, counting from the bottom). Equivalently, $X_r(q,m,i)$ is the closure of the set of $(a_0: ... : a_m)$ such that all $P_{[\ga]}$ have $\ge i$ common roots counting with multiplicities. 

The set of $H_\goth x$ is $i$-dimensional, hence we can expect that codim $X_r(q,m,i)$ in $P^m$ is $(q-1)i$. This is really so:
\medskip
{\bf Lemma 8.30.} codim $X_r(q,m,i)=(q-1)i$.
\medskip
{\bf Proof.} Let all $P_{[\ga]}$ have roots $x_1,\dots,x_i$ and let $P$ be the monic polynomial of degree $i$ having roots $x_1,\dots,x_i$. We denote $P_\ga:=P_{[\ga]}/P$. Let us consider a map of affine spaces $\vf: A^{m+1-(q-1)i}\to A^{m+1}$ defined as follows: the first $i$ coordinates of an element $t\in A^{m+1-(q-1)i}$ form coefficients of $P$ (all except the leading one who is 1), other coordinates of $t$ form (subsequently) coefficients of $P_0,\dots, P_{q-1}$, including their leading coefficients, and the matrix row $\vf(t)$ is formed by coefficients of $PP_0,\dots, PP_{q-1}$ (the same order as in the matrix $B$). If the first $i$ coordinates of $t$ do not belong to the discriminant variety defined by the condition that $P$ has no multiple roots, then Proj$(\vf(t))\in X_r(q,m,i)$. Obviously fibers of $\vf$ are finite, hence the lemma. $\square$
\medskip
The above results give us Theorem 8.6.
\medskip
Clearly $X_r(q,m,i)$ are irreducible. For $i=1$ we have (see Appendix, 1)
\medskip
{\bf Proposition 8.31.} deg $X_r(q,m,1)=m+1-q$. 
\medskip
{\bf 8.32. Case $m\not\equiv -1 \mod q$.} To define $X_r(q,m,i)$ we consider $B$ as a submatrix of $\goth N(P,n,k)$ formed by the first $q$ rows and the first $\lceil \frac{m+1}q \rceil$ columns, where $\lceil x \rceil:=\min \{n\in \n Z|n\ge x\}$ is the ceiling function. Matrices $\goth B(P,\mu)$ are defined for fractional $\mu\equiv l \mod 1$, where $l:=\frac {k+n}q$. We let $\goth B(P,l):=\goth N(P,n,k)$ and $\goth B(P,l-j)$ is obtained from $\goth N(P,n,k)$ by elimination of the last $qj$ rows and the right $j$ columns. It is easy to check that all the above definitions and results hold for this case. 
\medskip
We see that 3 formulas for the codimension of (sub)varieties of $X(q,n,m,i)$ in $P^m$, namely $ i(i+n)$, $ i(q-1)$, $i(q-1)+(q-2)(n-1)$, coincide for $n=1$, $q=i+2$. As an example, we consider 
\medskip
{\bf 8.33. Case $q=3$, $n=1$.} 
\medskip
\noindent
We have $m=2k+1$ is odd, $l=\frac{m+1}6$. $X_r(3,m,1)$ is of codimension 2, conjecturally $X_l (3,1,m,1)$ and $X(3,1,m,1)$ are also of codimension 2. In this case we have deg $X(3,1,m,1)=\binom{k+1}{2}$, deg $X_r(3,m,1)=2k-1$, deg $X_l(3,1,m,1)=\binom{k-1}{2}$ (see (8.23.1)), hence if $X_l (3,1,m,1)$ is irreducible (we are sure that it is, but we have no proof) then $$X(3,1,m,1)=X_r(3,m,1)\cup X_l(3,1,m,1)$$

{\bf Conjecture 8.33.1.} $X_r(3,m,1)\cap X_l(3,1,m,1)$ is of codimension 3 in $P^m$. 
\medskip
{\bf Idea of the proof.} There exists a $1\times (k-2)$-matrix $\Lambda=(\la_1, \la_2,\dots, \la_{k-2})$ such that $\Lambda \goth B(P,l-1)=0$. This is equivalent to $\La_2P_{[2]}+\La_1P_{[1]}+\La_0P_{[0]}=0$ where $\La_2:=\la_1+\la_4x+\la_7x^2+\dots$, $\La_1:=\la_2+\la_5x+\la_8x^2+\dots$, $\La_0:=\la_3+\la_6x+\la_9x^2+\dots$. Let us consider now the resultantal variety $R\subset P^m$ of elements $(a_0,\dots, a_m)$ such that $P_{[2]},\ P_{[1]}$ have a common root. It is a hypersurface in $P^m$. We have $R\cap X_l(3,1,m,1)$ is of codimension 3 in $P^m$. If  $(a_0,\dots, a_m)\in R\cap X_l(3,1,m,1)$ and $r$ is a common root of $P_{[2]},\ P_{[1]}$, then either $\La_0(r)=0$ or $P_{[0]}(r)=0$. These two cases correspond to two irreducible components of $R\cap X_l(3,1,m,1)$. It is possible to show that both these components have codimension 3 in $P^m$. The second of them is contained in $X_r(3,1,m,1)$. $\square$
\medskip
{\bf For $k=2$} we have $X_l(3,1,5,1)=\emptyset$, $X_r(3,5,1)=X(3,1,5,1)=P^2\times P^1\subset P^5$ is the Segre inclusion, which is smooth of degree 3 and dimension 3. 
\medskip
{\bf For $k=3$} we have $X_l(3,1,7,1)$ is of degree 1, i.e. $P^5\subset P^7$, it is the set of zeroes of $a_2=a_5=0$, and  $X_r(3,7,1)$ is of degree 5. It is easy to see $X_l(3,1,7,1)\cap X_r(3,7,1)$ is of codimension 3 in $P^7$ and degree 4: it is a resultantal variety given by the equation Res$(a_1+a_4x+a_7x^2, a_0+a_3x+a_6x^2)=0$ in $P^5=\{(a_0:a_1:a_3:a_4:a_6:a_7)\}$. 
\medskip
Let us describe $X(3,1,7,2)$. It is the intersection of $X(3,1,7,1)$ with the set $Y_2$ of zeroes of polynomials $H_{1i}$ (see 8.3.1), $i=0,1,2$. It is easy to see that we have: $X_l(3,1,7,1)\cap Y_2=P^2\times P^1\subset P^5$ is the Segre inclusion. Computer calculations show that $X_r(3,7,1)\cap Y_2$ is irreducible of dimension 3 and degree 9. 
\medskip
{\bf For $k=4$} we have $X_l(3,1,9,1)$ is of degree 3, it is a 7-dimensional cone over the image of the Segre inclusion, because $\goth B(P,l-1)=\left(\matrix a_2 & a_5& a_8\\ a_1&a_4 &a_7\endmatrix \right)$, and the set of $(a_2 : a_5: a_8: a_1:a_4 :a_7)$ such that rank $\goth B(P,l-1)=1$  is $P^2\times P^1$. We have $X_r(3,9,1)$ is of degree 7. As earlier we have $X_l(3,1,9,1)\cap X_r(3,9,1)$ is of codimension 3 in $P^9$. 
\medskip
Let us describe $X(3,1,9,2)$. It is the intersection of $X(3,1,9,1)$ with the set $Y_2$ of zeroes of polynomials $H_{1i}$, $i=0,1,2,3$. Computer calculations show that $X_r(3,9,1)\cap Y_2$ is the union of 3 components $C_7$, $C_{10}$, $C_{15}$ of codimension 4 and degrees 7, 10, 15 respectively, and $X_l(3,1,9,1)\cap Y_2$ is $C_7\cup C_{10}$. The same phenomenon of coincidence of components of different intersections occurs for $q=2$, see 9.7.9. 
\medskip
{\bf Remark.} Since deg $X_l(3,1,m,1)=$ deg $X(3,1,m-4,1)$, it is possible to conjecture that there is a relation between $X_l(3,1,m,1)$ and $X(3,1,m-4,1)$. The above examples show that this is true for $k=3,4$. Nevertheless, we do not know how to interpret this equality for $k=5$. For example, while $X(3,1,7,1)=P^5\cup X_r(3,1,7,1)$, it is known that the variety $X_l(3,1,11,1)$ does not contain $P^9$. 
\medskip
{\bf 9. Case $q=2$.} 
\medskip
We define in this section a variety $X(m,i)_{pr}$ defining a condition that two polynomials $P_{[0]}$, $P_{[1]}$ in one variable of degrees $\approx \frac m2$ have $i$ common roots. $X(m,i)_{pr}$ is an irreducible component of a variety $X(m,i)$ which is a complete intersection defined by coefficients of the characteristic polynomial of a modified Sylvester matrix. The meaning of other irreducible components is unclear. We do not know how to define $X(m,i)$ for the case $\deg (P_{[0]})-\deg(P_{[1]})\ne 0,\pm1$. 

In order to make the present section (almost) independent on the previous ones, we repeat some definitions. Let $(a_0: ...:a_m)\in P^m(\n C)$, $a_i=0$ for $i\not\in \{0,\dots,m\}$, 
\medskip
$P:=\sum_{i=0}^m a_i\theta^i$, $P_{[0]}:=a_0+a_2x+a_4x^2+a_6x^3+\dots$, 
\medskip
$P_{[1]}:=a_1+a_3x+a_5x^2+a_7x^3+\dots$ polynomials.  We consider a $(m-1)\times (m-1)$-matrix $\goth M(P,m)$ whose entries $\goth M(P,m)_{ij}=a_{2j-i}$. For example, for $m=6,7$ they are the following: 
\medskip
\centerline{$\goth M(P,m)$ for $m=6$ \ \ \ \ \ \ \ \ \ \ \ \ \ \ \ $\goth M(P,m)$ for $m=7$}
\medskip
$$\left(\matrix a_1&a_3&a_5&0&0\\ a_0&a_2&a_4&a_6&0 \\ 0&a_1&a_3&a_5&0 \\ 0&a_0&a_2&a_4&a_6\\ 0&0&a_1&a_3&a_5\endmatrix \right) \ \ \ \ \ \   \left(\matrix a_1&a_3&a_5&a_7&0&0\\ a_0&a_2&a_4&a_6&0& 0 \\ 0&a_1&a_3&a_5&a_7&0 \\ 0&a_0&a_2&a_4&a_6&0\\ 0&0&a_1&a_3&a_5&a_7\\ 0&0&a_0&a_2&a_4&a_6\endmatrix \right) \eqno{(9.1)}$$
\medskip
$\goth M(P,m)$ is a permutation of rows and columns of the Sylvester matrix of $P_{[0]}, P_{[1]}$. It coincides with $\goth M(P,0,k)^t$ of (3.1) and with $\goth N(P,0,m-1)$ of 8.11, case $n=0$, $q=2$, $k=m-1$. 
\medskip
Let $Ch(\goth M(P,m))$ be the $(-1)^{m-1}\cdot$ characteristic polynomial of $\goth M(P,m)$, i.e. 
$$Ch(\goth M(P,m))=|\goth M(P,m)-U\cdot I_{m-1}|=D(m,0)+D(m,1)\ U+$$ $$+D(m,2)\ U^2+\dots + D(m,m-2)\ U^{m-2}+(-U)^{m-1}$$
where $D(m,i)\in \n Z[a_0,\dots, a_m]$ are homogeneous polynomials of degree $m-1-i$.
\medskip
{\bf Definition 9.2.} $X(m,i)$ is a projective scheme 

Proj $\n C[a_0,\dots,a_m]/\{D(m,0),\dots,D(m,i-1)\}$. 
\medskip
{\bf Conjecture 9.3.} $\forall n$ \  Supp $X(m,i)$ = Supp $X(2,n,m,i)$ where $X(2,n,m,i)$ is defined in 8.4 (we consider $X(2,n,m,i)$ as a scheme as well). 
\medskip
{\bf Remark 9.4.} $X(m,i) = X(2,n,m,i)$ as sets, but not as schemes. This means that the multiplicity of an irreducible component of Supp $X(m,i)$ = Supp $X(2,n,m,i)$ depends on $n$, hence not all $H_{\al\be n}$ defining $X(2,n,m,i)$ (see (8.3.1)) belong to $[D(m,0),\dots, D(m,i-1)]$ --- the ideal generated by $D(m,*)$,  but only their powers $H_{\al\be n}^*$. The below tables 9.7.7 etc. give this multiplicity for $n=0$ and 1. See Remark 9.17.1 for an explicit example for $m=4$. 
\medskip
We see that the sets $X(2,n,m,i)$ depend on 2 parameters but not of 3 parameters. See 9.12 --- 9.17 for a justification of 9.3. They cover a simple case when some $H_{\al\be n}$ belong to $[D(m,0),\dots, D(m,i-1)]$. A non-trivial case $H_{\al\be n}\not\in [D(m,0),\dots, D(m,i-1)]$ is much more complicated, the authors have no proof that $\forall \ n,i,j$ \ $\exists \ga$ such that $H_{\al\be n}^\ga\in [D(m,0),\dots, D(m,i-1)]$.
\medskip
{\bf Remark 9.5.} Let us consider the $\n F_2$-case (see Sections 3, 5). Let $Q:=\sum_{i=0}^m b_i\theta^i$ where $b_i\in\n F_2$, and $Q=\prod_i\goth Q_i^{\al_i}$ its prime factorization in $\n F_2[\theta]$. We let $j:=m-\sum_i \deg \goth Q_i$. We have 
\medskip
{\bf Corollary 9.6.} Conjecture 9.3 implies: let $\be$ be the minimal number such that $D(m,\be)(b_0,\dots, b_m)\ne0$. Then $\be=j$. 
\medskip
{\bf Proof.} We apply (5.6.1) to our case ($P=1$, $P_1=Q$). We have $L(\goth C,T)=1+T$, $L_{\goth Q_i}(\goth C_P,T)^{-1}=1+\goth Q_i(t)T^{\deg \goth Q_i}$, hence $\deg_T L(\goth C_P,T)=1+\sum_i \deg \goth Q_i$ and $r_\infty=j$ (the summand 1 comes from the trivial part of $L(\goth C_P,T)$). $\square$
\medskip
{\bf 9.7. Properties of $X(m,i)$.} Most of the below properties are conjectural. Below we give some proofs. Source of evidence: computer calculations of resultants, see Appendix, 4. 
\medskip
{\bf Conjecture 9.7.1.} $X(m,i)$ is the complete intersection of the hypersurfaces $\{D(m,0)=0\}$,  $\{D(m,1)=0\}$, ... ,  $\{D(m,i-1)=0\}$. Particularly, codim $X(m,i)=i$, deg $X(m,i)=(m-1)(m-2)...(m-i)$ (not all multiplicities of irreducible components of $X(m,i)$ are 1, see below). 
\medskip
{\bf Conjecture 9.7.2.} All irreducible components of $X(m,i)$ have the same codimension $i$.
\medskip
Let $(a_0:...:a_m)\in \Supp X(m,i)$ be a generic point. The Jordan form of $\goth M(P,m)$ has $i$ zeroes on the diagonal, hence its block structure defines a partition of $i$. We denote the set of partitions of $i$ by $P(i)$. This means that we have a map $\pi: \Supp X(m,i)_g\to P(i)$ (subscript $g$ means generic). We can consider $\pi$ as a map from the set of irreducible components of $X(m,i)$ to $P(i)$. 
\medskip
{\bf Remark.} It is possible to consider analogs of $X(m,i)$ for general determinantal varieties $D(\de,\al,\al)$, see (8.15). Namely, in notations of (8.15), we can consider $D(\al+1,i)$ for $C_{\al\al}$ and $X(m,i)\subset P^{\al^2-1}$ --- the varieties of their zeroes. We have the same map $\pi$. But for this case there is no equidimensionality: the codimension in $ P^{\al^2-1}$ of matrices whose Jordan form has a $\left(\matrix 0&1\\ 0&0\endmatrix\right)$-block (resp. a $\left(\matrix 0&0\\ 0&0\endmatrix\right)$-block) is 2 (resp. 4). 
\medskip
{\bf Definition 9.7.3.} $X(m,i)_{pr}\subset P^m$ is the closure of $\pi^{-1}$ of the partition $i=1+...+1$, i.e. the Zariski closure of the set of $a_0, \dots, a_m$ such that the 0-Jordan block of $\goth M(P,m)$ is the 0 matrix. 
\medskip
$X(m,i)_{pr}$ is called the principal component of $X(m,i)$. Obviously $(a_0,\dots, a_m)\in X(m,i)_{pr} \iff$ corank of $\goth M(P,m)\ge i$. 
\medskip
{\bf Proposition 9.7.4.} $X(m,i)_{pr}=X_r(2,m,i)=X_l(2,n,m,i)=X_\goth N(2,n,m,i)$ for all $n$ where $X_r(2,m,i)$, $X_l(2,n,m,i)$, $X_\goth N(2,n,m,i)$ are defined in Section 8. We repeat here their definition. Let $\goth S(P,m,j)$ be a submatrix of $\goth M(P,m)$ obtained by elimination of its last $2j$ rows and last $j$ columns (in notations of 8.18, 8.32,  $\goth S(P,m,j)=\goth B(P,\frac{m-1}2-j)$ ). Formally, we can consider $\goth S(P,m,j)$ for negative $j$ as well --- the above formula continues to hold. 
\medskip
$X_r(2,m,i)$ is the set of $(a_0:...:a_m)$ such that $P_{[0]}$, $P_{[1]}$ have $\ge i$ common roots, counting with multiplicities;
\medskip
$X_l(2,n,m,i)$ (it does not depend on $n$) is the set of $(a_0:...:a_m)$ such that $\goth S(P,m,i-1)$ is not of the maximal rank. 
\medskip
$X_\goth N(2,n,m,i)$ is the set of $(a_0:...:a_m)$ such that the rank of $\goth S(P,m,-n)$ is $\le m+n-1-i$. 
\medskip
{\bf Conjecture 9.7.5.} $X(m,i)_{pr}$ is one of the irreducible components of 

\noindent
$\Sing_{i-1}(X(m,1))$ 
(as usual, for a variety $Y$ we let $\Sing_0(Y)=Y$, $\Sing_{k+1}(Y)=\Sing(\Sing_k(Y))$ ). Particularly, for $i=2$ we have $X(m,2)_{pr}=\Sing (X(m,1))$, while for $i=3$ \ Sing$(X(m,2)_{pr})$ has two irreducible components: $X(m,3)_{pr}$ and a component corresponding to the case of double common root of $P_{[0]}$, $P_{[1]}$ (the authors are grateful to A.I. Esterov who indicated them that $\Sing_{i-1}(X(m,1))$ has components distinct from $X(m,i)_{pr}$). 
\medskip
{\bf Proposition 9.7.6.} $\deg X(m,i)_{pr}=\binom{m-i}{i}$.
\medskip
{\bf Conjecture 9.7.7.} For $i=2$ we have $X(m,2)$ is the union of 2 irreducible components $C_{21}$ and $X(m,2)_{pr}$. We have the following table: 
\settabs 5 \columns
\medskip
\+ Components &&$C_{21}$&  $X(m,2)_{pr}$\cr
\medskip
\+ Degrees&&$2(m-2)$& $\binom{m-2}{2}$\cr
\medskip
\+ $\pi$-images &&$\{2=2\}$&  $\{2=1+1\}$\cr
\medskip
\+ Multiplicities  in $X(m,2)$&&1&2\cr
\medskip
\+ Multiplicities  in $X(2,1,m,2)$&&1&1\cr
\medskip
{\bf Conjecture 9.7.8.} For $i=3$ we have $X(m,3)$ is the union of 4 irreducible components $C_{31}$,  $C_{32}$,  $C_{33}$,  $X(m,3)_{pr}$. We have the following table: 
\medskip
\settabs 6 \columns
\+ Components &&$C_{31}$&  $C_{32}$&  $C_{33}$&  $X(m,3)_{pr}$\cr
\medskip
\+ Degrees&&$4(m-3)$& $m-3$& $4\binom{m-3}{2}$& $\binom{m-3}{3}$\cr
\medskip
\+ $\pi$-images &&$\{3=3\}$&$\{3=2+1\}$ & $\{3=2+1\}$& $\{3=1+1+1\}$\cr
\medskip
\+ Multiplicities in $X(m,3)$&&1&2&3&6\cr
\medskip
\+ Multiplicities &&1&1&1&2\cr
\+ in $X(2,1,m,3)$\cr
\medskip
 Further, $C_{21}\cap \{D(m,2)=0\}=C_{31}\cup C_{33}$, 
\medskip
$X(m,2)_{pr}\cap \{D(m,2)=0\}=C_{32}\cup C_{33}\cup 3\ X(m,3)_{pr}$. 
\medskip
{\bf 9.7.9.} This means that it can happen that different irreducible components of $X(m,i)$ crossed with a hypersurface $\{D(m,i)=0\}$ have coinciding irreducible components of their intersection. 
\medskip
{\bf Conjecture 9.7.10.} For any fixed $i$ and varying $m$ the quantity of irreducible components of $X(m,i)$, their multiplicities and $\pi$-images do not depend on $m$ (exception: if the degree $=0$ then the corresponding component is empty). The degree of any irreducible component is $c\binom{m-i}{j}$ where $c$ and $j\le i$ do not depend on $m$. Moreover $j<i$ unless of the principal component. This information is presented in the below table:
\settabs 12 \columns
\medskip
\+ Components &&&$C_{i1}$& \dots&  &  $C_{i*}$& &\dots&& $X(m,i)_{pr}$ \cr
\medskip
\+ Degrees&&&$2^{i-1}(m-i)$&&$c_*\binom{m-i}{j_*}$, where $1\le j_*<i$&&&&&  $\binom{m-i}{i}$\cr
\medskip
\+ $\pi$-images &&&$\{i=i\}$&&&$\ne\{i=i\}$,&&& $\{i=1+1+...+1\}$\cr
\medskip
\+&&&&&$\ne\{i=1+1+...+1\}$\cr
\medskip
\+ Multiplicities &&&1&&&?&&&&$i!$\cr
\+ in $X(m,i)$\cr
\medskip
\+ Multiplicities &&&1&&&?&&&&?\cr
\+ in $X(2,n,m,i)$\cr
\medskip
{\bf Conjecture 9.7.11.} Description of components having $j=1$ ($j$ from the above table).
\medskip
Let $i$ be fixed. We consider components of $X(m,i)$ having $j=1, \ 2$. We denote by $\al_i$, resp. $\be_i$ the quantity of irreducible components of $X(m,i)$ having $j=1$, resp. $j=2$. We denote by $$c_{i11}\binom{m-i}{1},\  \dots ,\ c_{i1\al_i}\binom{m-i}{1}, \ \ c_{i21}\binom{m-i}{2},\dots,\ c_{i2\be_i}\binom{m-i}{2}$$ the degrees of these irreducible components of $X(m,i)$. Then for $i+1$ we have $$\al_{i+1}=\al_i+\be_i, \hbox{ the numbers $c_{i+1,1,*}$ are }2c_{i11}, \dots, \ 2c_{i1\al_i}, \ c_{i21}, \dots, c_{i2\be_i}$$ For example, there are 3 irreducible components of $X(m,4)$ whose degrees are $8(m-4)$, $2(m-4)$, $4(m-4)$. Particularly, $\forall i$ there exists the only irreducible component of $X(m,i)$ whose $\pi$-image is $\{i=i\}$. Its degree is $2^{i-1}(m-i)$ and its multiplicity in both $X(m,i)$ and $X(2,n,m,i)$ is 1 (the first column of the above tables). 
\medskip
Let us denote by $OP(i)$ the set of ordered partitions of $i$ (for example, 3=2+1 and 3=1+2 are two different ordered partitions of 3). There is a map $f: OP(i)\to P(i)$ forgetting ordering. Let $IR(X(m,i))$ be the set of irreducible components of $X(m,i)$. 
\medskip
{\bf Supposition 9.7.12.} Let $m\ge 2i$. There is an isomorphism $\al: IR(X(m,i))\to OP(i)$ such that $f\circ\al=\pi$. Particularly, $\# \ IR(X(m,i))=2^{i-1}$ (if $m\ge 2i$). 
\medskip
{\bf Supposition 9.7.13.} All irreducible components of $X(m,i)$ are rational varieties. 
\medskip
{\bf 9.7.14. Problems.} To find quantity of irreducible components of $X(m,i)$ (is 9.7.12 true?), their $\pi$-images, degrees, singularities, multiplicities, nilpotent part of rings, intersections etc. What are  multiplicities of components of $X(2,n,m,i)$?
\medskip
{\bf Remark 9.8.} We can consider a more general situation. Let $A\in GL_{m-1}(\n C)$ be any fixed matrix, $S=S(P_{[0]}, P_{[1]})$ the Sylvester matrix of $P_{[0]}, P_{[1]}$. Instead of the matrix $\goth M(P,m)$ we can consider the matrix $AS$, its characteristic polynomial and varieties of zeroes of its coefficients. Shall we get some interesting results? 

Obviously $X(m,i)_{pr}$ does not depend on $A$, other $C_{ij}$ are different (for different $A$) as sets of points. $X(3,2)$ is a non-singular plane conic for all $A$, while for $A=I_3$ we have: $X(4,3)$ is a normcubic $\cup$ a triple $P^1$, i.e. we have a type distinct from the one described above in 9.7.8. 
\medskip
{\bf Lemma 9.9.} For any $m$, $n$ we have: $X_l(2,n,m,i)=X_r(2,m,i)$. 
\medskip
{\bf Proof.} Let $(a_0,\dots, a_m)\in X_r(2,m,i)$, i.e. $P_{[0]}$, $P_{[1]}$ have (at least) $i$ common roots (counting with multiplicities). This means that there exists a polynomial $P$ of degree $i$, polynomials $P_0$, $P_1$ such that $P_{[0]}=PP_0$, $P_{[1]}=PP_1$. This means that $P_{[0]}P_1=P_{[1]}P_0$. This gives us a non-trivial linear dependence of rows of $\goth S(P,m,i-1)$, i.e. it is not of the maximal rank. All these arguments are convertible, i.e. if $(a_0,\dots, a_m)\in X_l(2,n,m,i)$ then $(a_0,\dots, a_m)\in X_r(2,m,i)$. $\square$
\medskip
Lemma 9.9 and 8.21a give us a proof of 9.7.4. 
\medskip
Let us justify some other conjectures. Results of Section 8 imply codim $X(m,i)=i$, hence we can apply (8.15.2) to $\goth S(P,m,i)$. This gives us 9.7.6. Further, let $C_{ik}$ be the $k$-th irreducible component of $X(m,i)$. We denote its degree by $d(C_{ik})$ and its multiplicity in $X(m,i)$ by $\mu(C_{ik})$. We have 
$$\sum_k d(C_{ik})\mu(C_{ik})=\deg X(m,i)=(m-1)(m-2)\dots(m-i)\eqno{(9.11)}$$ According 9.7.10, for all $k$ except $k_{max}$ --- the one that corresponds to the principal component, $d(C_{ik})$ is a polynomial in $m$ of degree $<i$, and $\mu(C_{ik})$ does not depend on $m$. Comparing the leading coefficients of the both sides of 9.11 we get that for $C_{ik_{max}}=X(m,i)_{pr}$ we have $\mu(C_{ik_{max}})=i!$ . 
\medskip
Let us justify 9.7.11. We fix $m$ and we consider the cases $i=m-2$, $i=m-1$. $X(m,m-2)$ is a surface whose components have $j=1,\ 2$ from 9.7.11, and $X(m,m-1)$ is the intersection of these components with a hyperplane $\{D(m,m-2)=0\}$. We can expect that all these intersections are distinct, hence there is 1 -- 1 correspondence between components of $X(\goth m,i)$ having $j=1, \ 2$, and components of $X(\goth m,i+1)$ having $j=1$ (recall that the set of components depends only on $i$ and not of $\goth m$). For $\goth m=m$ their degrees coincide, i.e. for the $k$-th component of $X(m,m-2)$ having $j=1$ its degree is $c_{i1k}\binom{m-(m-2)}{1}=c_{i+1,1k}\binom{m-(m-1)}{1}$, i.e. $c_{i+1,1k}=2c_{i1k}$. Analogically for components of $X(m,m-2)$ having $j=2$. 
\medskip
{\bf Proposition 9.12.} Conjecture 9.3 is true for $i=1$. 
\medskip
{\bf Proof.} This follows immediately from Theorem III. Really, Theorem III is a much stronger result. To prove the present proposition, it is sufficient to read the first 4 lines of the proof of Lemma 7.1 of Part III: 
\medskip
"This means that if $D(m,0)=0$ then $\forall L$ we have $|A_L|=0$." 
\medskip
$A_L$ is defined in (5), Part III, and the above lines. According (5), Part III, $|\goth M(P,n,k)|$ is a linear combination of $|A_L|$ for all $L$, hence if $D(m,0)=0$ then $|\goth M(P,n,k)|=0$. This means that $X(m,1)\subset X(2,n,m,1)$. The converse inclusion is "obvious" (see Theorem III for justification). $\square$
\medskip
{\bf Proposition 9.14.} $\forall \  m, \ n, \ i$  we have $$H_{i0n}=\pm D(m,i-n) \pm a_0  D(m,i-n+1)$$  $$H_{i,n(k-i),n}=\pm  D(m,i-n) \pm a_m  D(m,i-n+1)$$ (we let $D(m,j)=0$ if $j\not\in \{0,\dots,m-1\}$ ). 
\medskip
{\bf Proof.} To find $H_{i0n}$ we let $t=0$ in $\goth M(P,n,k)$. We get $\{I_k-\goth M(P,n,k)T\}_{t=0}=\left(\matrix *_{11} &  *\\0 & I_{m-1}-\goth M(P,m)T\endmatrix\right)$ (the $(n, \ m-1)$-block form) where $*_{11}$ is an upper-triangular $n\times n$-matrix with $(1,1,\dots,1,1\pm a_0T)$ at the diagonal. This gives us immediately the formula for $H_{i0n}$. The formula for $H_{i,n(k-i),n}$ follows from (8.3.3). $\square$
\medskip
{\bf Remark 9.15.} (a) For $n=1$ we have $H_{m-2,1}=\pm D(m,m-3)\pm D(m,m-2)^2$. 
\medskip
\noindent
(b) For $n=2$, any $m$ we have $H_{112}=2(\pm a_0^2D(m,1)\pm(a_0+a_1)D(m,0))$ 
\medskip
\noindent
(c) and, by (8.3.3), $H_{1,2m-1,2}=2(\pm a_m^2D(m,1)\pm(a_m+a_{m-1})D(m,0))$. 
\medskip
{\bf 9.16. Example: $m=3$. (a)} $n=0$. We have $\goth M(P,3)=\left(\matrix a_1 &  a_3\\ a_0&a_2\endmatrix\right)$, 

\noindent
$D(3,0)=\det \goth M(P,3)=-\Res(a_3t+a_1, a_2t+a_0)= -a_0a_3+a_1a_2$; $X(3,1)$ is the quadric surface $\{D(3,0)=0\}$. It is non-singular. $D(3,1)=a_1+a_2$, $X(3,2)_{pr}=\emptyset$, $X(3,2)=C_{21}$ is a non-singular conic line in $P^2$. 
\medskip
{\bf (b)} $n=1$. (9.12.1), (9.14), (9.15a) imply that Conjecture 9.3 is true for $n=1$, $m=3$, all $i$. 
\medskip
{\bf (c)} $n=2$. For this case $H_{0i2}$ can be found using (9.13.1). (9.15b,c) show that $H_{1i2}$ belong to the ideal generated by $D(3,0)$ and $D(3,1)$ for $i=1,7$. Explicit calculations show that the same is true for all $H_{1i2}$ (compare with (9.17c)), hence $X(3,1)\subset X(2,2,3,1)$. Further calculations show that $X(3,1)= X(2,2,3,1)$, $X(3,2)= X(2,2,3,2)=\emptyset$. This means that Conjecture 9.3 is true for $m=3$, $n=1,2$. 
\medskip
{\bf 9.17. Case $m=4$.  (a)} $n=0$. We have $\goth M(P,4)=\left(\matrix a_1 &  a_3&0\\a_0&a_2&a_4\\ 0&a_1 &  a_3\endmatrix\right)$, 
\medskip
\noindent
$D(4,0)=\det\goth M(P,4)=\pm \Res(a_3t+a_1, a_4t^2+a_2t+a_0)=a_1a_2a_3-a_0a_3^2 -a_1^2a_4$;
\medskip
\noindent
 $X(4,1)$ is the cubic threefold $\{D(4,0)=0\}$.  
\medskip
\noindent
 We have $\goth S(P,4,1)=\left(\matrix a_1 &  a_3\endmatrix\right)$, i.e. $\Sing(X(4,1))$ is the plane $a_1=a_3=0$.
\medskip
\noindent
We have $D(4,1)=-a_0a_3+a_1a_2+a_1a_3-a_1a_4+a_2a_3$. Variety $\{D(4,1)=0\}$ is a cone whose vertex is a point $s:=(1:0:1:0:1)\in \Sing(X(4,1))$. We have: $X(4,2)=X(4,1)\cap \{D(4,1)=0\}$ (complete intersection), 
\medskip
$X(4,2)=C_{21}\cup X(4,2)_{pr}$, where $X(4,2)_{pr}=\Sing(X(4,1))$ and $C_{21}$ is isomorphic to $P^1\times P^1$ with two glued points. $C_{21}$ has degree 4 in $P^4$. It is given by parametric equations $t: P^1\times P^1 \to P^4$:
\medskip
$t(\la_0:\la_1, \ c_1:c_3):=\la_0(0:-c_3^2:-c_1c_3:c_1c_3:c_1^2)+\la_1(-c_3^2:-c_1c_3:c_1c_3:c_1^2:0)$ 
\medskip
\noindent
where $(\la_0:\la_1) \in P^1$,  $(c_1:c_3)\in P^1$ are parameters. Two points $t(\la_0:\la_1, \ c_1:c_3)$ are glued in $s$: $t(1:-\zeta_3, \ 1:\zeta_3)=t(1:-\zeta_3^2, \ 1:\zeta_3^2)$ where $\zeta_3$ is a primitive cubic root of 1. There is no more glueing. 
\medskip
We have $C_{21}\cap X(4,2)_{pr}$ is a singular cubic on $X(4,2)_{pr}$, whose singular point is $s$. Its parametric equation is $(c_3^3:0:-c_1^2c_3-c_1c_3^2:0:c_1^3)$. 
\medskip
{\bf Description of $X(4,3)$.} We have $D(4,2)=a_1+a_2+a_3$, hence $\{D(4,2)=0\}$ is a $P^3$, 
$X(4,3)=C_{31}\cup C_{32}$, where $C_{31}=C_{21}\cap \{D(4,2)=0\}$, $C_{32}=X(4,2)_{pr}\cap \{D(4,2)=0\}$, $C_{33}=X(4,3)_{pr}=\emptyset$. 
\medskip
\noindent
$C_{31}$ is a non-singular rational curve of degree 4 given by parametric equations $(-c_3^4:-c_1^2c_3^2-c_1c_3^3:-c_1^3c_3+c_1c_3^3:c_1^3c_3+c_1^2c_3^2:c_1^4)$, where $(c_1:c_3)$ are as above. 
\medskip
\noindent
$C_{32}$ is a straight line, its equations are $a_1=a_2=a_3=0$. $C_{31}\cap C_{32}$ consists of 2 points 
$(0:0:0:0:1)$ and $(1:0:0:0:0)$. 
\medskip
{\bf (b)} $n=1$. (9.12.1), (9.14) imply $H_{0i}=\pm a_i D(4,0)$ ($i=0,\dots,4$); 
$H_{10}=\pm D(4,0)\pm a_0 D(4,1)$, $H_{13}=\pm D(4,0)\pm a_4 D(4,1)$. We have
\medskip
\noindent
$H_{11}=-a_0a_1a_3+a_0a_1a_4+a_0a_3^2-a_0a_3a_4+a_1^2a_2+a_1^2a_3-a_1^2a_4+a_1a_2a_3-a_1a_2a_4+a_2^2a_3$
\medskip
\noindent
$H_{12}=a_0a_1a_4+a_0a_2a_3+a_0a_3^2-a_0a_3a_4-a_1^2a_4-a_1a_2^2-a_1a_2a_3-a_1a_3^2+a_1a_3a_4- a_2a_3^2$
\medskip
We have $H_{11}, H_{12}\not\in [D(4,0), D(4,1)]$, $H_{11}^2, H_{12}^2\in [D(4,0), D(4,1)]$, hence $X(4,2)=X(2,1,4,2)$ as the sets of points, i.e. Conjecture 9.3 is true for this case. 
\medskip
{\bf Remark 9.17.1.} $\{D(4,1)=0\}\cap \{H_{11}=0\}$ is $C_{21}\cup X(4,2)_{pr}\ \cup $ \{the $P^2$ having equations $a_2=a_0$, $a_4=a_0+a_3$\}. This means that $$X(4,2)=\hbox{Proj }\n C[a_0,\dots, a_4]/\{D(4,0), D(4,1)\}\ne X(2,1,4,2)\hbox{ as schemes}$$ 

Equations of $H_{2i}$ are the following (see 9.14, 9.15a):
\medskip
$H_{20}=D(4,1)+a_0D(4,2)=a_0a_1+a_0a_2+a_1a_2+a_1a_3-a_1a_4+a_2a_3$, it is a cone whose vertice is 1 point $(0:1:0:-1:0)$; 
\medskip
$H_{21}=D(4,1)-D(4,2)^2=-a_0a_3-a_1^2-a_1a_2-a_1a_3-a_1a_4-a_2^2-a_2a_3-a_3^2$, it is non-singular; 
\medskip
$H_{22}=D(4,1)+a_4D(4,2)=-a_0a_3+a_1a_2+a_1a_3+a_2a_3+a_2a_4+a_3a_4$, it is a cone whose vertice is 1 point $(0:1:0:0:-1)$;
\medskip 
The set of singular linear combinations $\la_0 H_{20} + \la_1 H_{21} + \la_2 H_{21}$ is a curve of degree 5 in $P^2$. It is the union of a non-singular conic and 
a triple $P^1$ given by the equation $\la_0 + \la_1 + \la_2=0$. Quadrics which correspond to this $P^1$ are of rank 2, i.e. they 
are the union of two $P^3$. One of these $P^3$ is $\{D(4,2)=0\}$ and another $P^3$  contains the plane 
$\{a_0=a_4, \ a_1+a_2+a_3+a_4=0\}$. 
\medskip
{\bf (c)} $n=2$. According (9.13), all $H_{0i2}\in [D(4,0)]$. According (9.14), (9.15b,c) $H_{1i2}\in [D(4,0), D(4,1)]$ for $i=0,1,7,8$. It is possible to check that $H_{122}\ne C_1 D(4,0)+C_2 D(4,1)$ where $C_1$, $C_2$ are polynomials of degrees 1, 2 respectively.

\medskip
{\bf Part III. Calculation of a determinant. }
\medskip
Theorem III --- the result of the present part --- grew from a proof of Conjecture 9.3 for $i=1$, see (II.6) and (9.12) for details. In order to make this part independent on the rest of the paper, we repeat and slightly modify definitions. Let $q\ge2$, $n\ge0$, $m\ge1$ be integers such that $k=k_n=k(m,n,q):=\frac{m+n}{q-1}-1$ is integer $\ge1$. Let $a_0,a_1,\dots,a_m$ and $t$ be variables. 
\medskip
The $k\times k$-matrix $\widehat {\goth M}(a_*,n,k)$ whose entries depend on $a_*$, $t$ is defined by the formula (we shall need only the case $q=2$)
$$\widehat {\goth M}(a_*,n,k)_{i,j}=\sum_{l=0}^n \binom{n}{l}a_{qj-i-l}\ t^{n-l}$$ (throughout all this part, $a_*=0$ if $*\not\in [0,\dots, m]$). For the reader's convenience, we give the explicit form of $\widehat {\goth M}(a_*,2,k)$ for $q=2$, $n=2$: 
$$\left(\matrix a_1t^2+2a_0t&a_3t^2+2a_2t+a_1&a_5t^2+2a_4t+a_3&\dots&0\\
a_0t^2&a_2t^2+2a_1t+a_0&a_4t^2+2a_3t+a_2&\dots&0\\0&a_1t^2+2a_0t&a_3t^2+2a_2t+a_1&\dots&0\\
\dots&\dots&\dots&\dots&\dots\\
0&0&\dots&\dots&2a_mt+a_{m-1}\endmatrix\right)$$
For $q=2$, $n=0$ the matrix $\widehat {\goth M}(a_*,0,k)$ is equal to $\goth M(P,m)$, see (II.9.1) for $m=6,7$.

\medskip
{\bf Theorem III.} For $q=2$, any $m,n$ we have $$|\widehat {\goth M}(a_*,n,k_n)|=(2t)^{n\choose 2}\cdot (a_0+a_1t+a_2t^2+\dots +a_mt^m)^n\cdot|\widehat {\goth M}(a_*,0,k_0)|\eqno{(1)}$$

{\bf Remark 2.} $\widehat {\goth M}(a_*,0,k_0)$ does not contain $t$. This means that $|\widehat {\goth M}(a_*,0,k_0)|$ can be considered as a common factor of coefficients of $|\widehat {\goth M}(a_*,n,k_n)|$ at $t$.
\medskip
{\bf Remark 3.} The matrix $\widehat {\goth M}(a_*,n,k)$ is a version of a matrix $\goth M(P,n,k)$ defined in (I.3.1) ($P$ of (I.3.1) is $(a_0,a_1,\dots,a_m)$ of the present part) obtained by changing all minus signs in $\goth M(a_*,n,k)$ by the plus signs, and transposition: $\widehat {\goth M}^t(a_*,n,k)(t)=(-1)^n\goth M(a_*,n,k)(-t)$. Further, for $q=2$ the matrix $\goth M(P,k_0)$ defined in (II.9.1) is $\widehat {\goth M}(a_*,0,k_0)$.
\medskip
We give two proofs of the theorem. Proof B in (17) is much shorter, but intermediate results of the Proof A are of independent interest and can be used for generalization of the theorem. 
\medskip
{\bf Proof A.} First, we consider the formula for the determinant of a $k\times k$-matrix $\widetilde {\goth M}(*,n,k)$ depending on $a_{ij}$, $i=1,\dots,k$, $j=1,\dots,k+n$, defined by the formula 
$$\widetilde {\goth M}(*,n,k)_{i,j}=\sum_{l=0}^n \binom{n}{l}a_{i,j+l}t^{n-l}$$ which is more general than the matrix $\widehat{\goth M}^t(*,n,k)$ (its version $\widetilde {\goth M}^-(*,1,k)$ is given in (A3.1) ).  For the reader's convenience, here we give the explicit form of $\widetilde {\goth M}(*,n,k)$ for $n=2$: 
\medskip
\newpage
$$\widetilde {\goth M}(*,2,k)=$$ $$\left(\matrix a_{11}t^2+2a_{12}t+a_{13} &  a_{12}t^2+2a_{13}t+a_{14}  &   \dots & a_{1k}t^2+2a_{1,k+1}t+a_{1,k+2} \\ a_{21}t^2+2a_{22}t+a_{23} &  a_{22}t^2+2a_{23}t+a_{24}  &   \dots & a_{2k}t^2+2a_{2,k+1}t+a_{2,k+2} \\  a_{31}t^2+2a_{32}t+a_{33} &  a_{32}t^2+2a_{33}t+a_{34}  &   \dots & a_{3k}t^2+2a_{3,k+1}t+a_{3,k+2} \\
\dots & \dots & \dots  & \dots \\ a_{k1}t^2+2a_{k2}t+a_{k3} &  a_{k2}t^2+2a_{k3}t+a_{k4}  &   \dots & a_{kk}t^2+2a_{k,k+1}t+a_{k,k+2} \endmatrix \right) \eqno{(4)}$$ (for the case of a general $n$ the numerical coefficients of any entry of the matrix are ${n\choose 0}$, ${n\choose 1}$, ${n\choose 2}, \dots, {n\choose n-1}, {n\choose n}$   ). 

We denote by $A_l$, $l=1,\dots, k+n$, the $l$-th column of $\{a_{ij}\}$, i.e. $A_l=\left(\matrix a_{1l}&a_{2l}&\dots&a_{kl}\endmatrix \right)^t$, and for an ordered sequence $L=(l_1,\dots, l_k)$ we denote by $|A_L|:=|\matrix A_{l_1}&A_{l_2}&\dots&A_{l_k}\endmatrix |$ the determinant of the $k\times k$ matrix formed by columns $A_{l_1}, A_{l_2},\dots,A_{l_k}$. Hence, we consider $L$ satisfying the condition $1\le l_1<l_2<\dots <l_k\le n+k$. We denote the set of these $L$ by $\goth L$. Further, for these $L$ we denote by $\bar L:=(\mu_1,\dots,\mu_n)$ the complement to $L$ in $[1,\dots,n+k]$, and we denote $d(L):=\mu_1+\mu_2+...+\mu_n-{n+1\choose 2}$. 
We have $$|\widetilde {\goth M}(*,n,k)|=\sum_{L\in \goth L}c(L) \  |A_L |\ t^{d(L)}\eqno{(5)}$$
where $c(L)\in \n Z$ are some coefficients. 
\medskip
{\bf Proposition 6.} $c(L)= \frac{\prod_{1\le i<j\le n}(\mu_j-\mu_i)}{(n-1)!!}$ where $n!!:=1!\cdot2!\cdot...\cdot n!$. 
\medskip
{\bf Proof.} Let us consider a $k\times (n+k)$-matrix $B_1$ whose $i$-th row is the set of coefficients of the $i$-th column of $\goth M(*,n,k)$ as a linear combination of $A_1,\dots, A_{k+n}$: 
$$B_1=\left(\matrix {n\choose 0} & {n\choose 1} & {n\choose 2}&\dots&{n\choose n}&0&0&\dots&0\\ 0&{n\choose 0} & {n\choose 1} & {n\choose 2}&\dots&{n\choose n}&0&\dots&0\\ \dots&\dots&\dots&\dots&\dots&\dots&\dots& \dots&\dots \\ 0&0&\dots&0&{n\choose 0} & {n\choose 1} & {n\choose 2}&\dots&{n\choose n}
\endmatrix \right)$$ We have $c(L)=$ the $(l_1,\dots,l_k)$-th minor of $B_1$. 
Now we use [K], page 30, Theorem 26, formula 3.13 (the authors are grateful to Suvrit Sra, Gjergji Zaimi, Christian Stump, Christian Krattenthaler who indicated them this information): 

Let $\nu,A,B,L_1,\dots,L_\nu$ be numbers. Then
$$\left| \matrix {BL_1+A\choose L_1+1}&{BL_1+A\choose L_1+2}&\dots&{BL_1+A\choose L_1+\nu}\\ \\ {BL_2+A\choose L_2+1}&{BL_2+A\choose L_2+2}&\dots&{BL_2+A\choose L_2+\nu}\\ \\ \dots&\dots&\dots&\dots \\ \\ {BL_\nu+A\choose L_\nu+1}&{BL_\nu+A\choose L_\nu+2}&\dots&{BL_\nu+A\choose L_\nu+\nu}\endmatrix \right|=\prod_{1\le i<j\le\nu}(L_i-L_j)\prod_{i=1}^\nu(BL_i+A)! \ \cdot$$ $$\cdot \ \prod_{j=1}^{\nu-1}(A-B(j+1)+1)_{j} \left(\prod_{i=1}^\nu(L_i+\nu)!\right)^{-1}\left(\prod_{i=1}^\nu((B-1)L_i+A-1)!\right)^{-1}\eqno{(6.1)}$$ where $(x)_j:=x(x+1)(x+2)\dots(x+j-1)$. 
We transpose $B_1$ and write its columns in the inverse order. As a result, $B_1$ will have a form of the matrix of (6.1) with $\nu=k$, $A=n$, $B=0$, $L_i=l_i-k-1$.  Substituting these values to (6.1) we get: 
$$\prod_{1\le i<j\le\nu}(L_i-L_j)=\frac{(-1)^{n\choose2}\cdot(k+n-1)!!\cdot\prod_{1\le i<j\le n}(-(\mu_i-\mu_j))}{((\mu_1-1)!\cdot... \cdot(\mu_n-1)!)\cdot((k+n-\mu_1)!\cdot...\cdot(k+n-\mu_n)!)}\eqno{(6.2)}$$
$$\prod_{i=1}^\nu(BL_i+A)!=(n!)^{k}\eqno{(6.3)}$$
$$\prod_{j=1}^{\nu-1}(A-B(j+1)+1)_{j}=(n+1)\cdot(n+1)(n+2)\cdot\dots\cdot(n+1)(n+2)\dots(n+k-1)\eqno{(6.4)}$$
hence the product of (6.3), (6.4) is 
$$n!\cdot (n+1)!\cdot\dots\cdot(n+k-1)!=\frac{(n+k-1)!!}{(n-1)!!}\eqno{(6.5)}$$
Further, $$\prod_{i=1}^\nu(L_i+\nu)!=\frac{(n+k-1)!!}{(\mu_1-1)!\cdot... \cdot(\mu_n-1)!}\eqno{(6.6)}$$
$$\prod_{i=1}^\nu((B-1)L_i+A-1)!=\frac{(n+k-1)!!}{(k+n-\mu_1)!\cdot...\cdot(k+n-\mu_n)!}\eqno{(6.7)}$$ Substituting (6.2), (6.5), (6.6), (6.7) to (6.1), we get the desired. $\square$
\medskip
Now let us return to the case of $\widehat{\goth M}(a_*,n,k)$. For $q=2$ we have $k=m+n-1$, and $\widehat{\goth M}(*,n,k)=\widetilde {\goth M}^t(*,2,k)$, where $a_{ij}$ of $\widetilde {\goth M}(*,2,k)$ are equal to $a_{2i-j}$ of the statement of the theorem, and $a_*=0$ if $*\not\in [0,\dots,m]$. The columns $A_1,\dots, A_{k+n}$ of $\widetilde {\goth M}(*,2,k)$ after transposition become lines, and we consider a $(k+n)\times k$-matrix $\goth N(a_*,n,k)$ whose $i$-th line is $A_i^t$: $\goth N(a_*,n,k)_{ij}=a_{2j-i}$. For the reader's convenience, we give the explicit form of $\goth N(a_*,n,k)$: 
$$\left( \matrix a_1&a_3&a_5&\dots&\dots&\dots&\dots&0\\
a_0&a_2&a_4&\dots&\dots&\dots&\dots&0\\
0&a_1&a_3&a_5&\dots&\dots&\dots&0\\
0&a_0&a_2&a_4&\dots&\dots&\dots&0\\
\dots&\dots&\dots&\dots&\dots&\dots&\dots&\dots\\
0&\dots&\dots&\dots&\dots&a_{m-3}&a_{m-1}&0\\
0&\dots&\dots&\dots&\dots&a_{m-4}&a_{m-2}&a_{m}\\
0&\dots&\dots&\dots&\dots&a_{m-5}&a_{m-3}&a_{m-1}\endmatrix \right)$$
For $n=0$ we have $\goth N(a_*,0,k_0)=\widehat{\goth M}(a_*,0,k_0)$ is a square matrix which is a permutation of rows and columns of the Sylvester matrix of two polynomials $$P_{[0]}:=a_0+a_2x+a_4x^2+a_6x^3\dots $$ $$ P_{[1]}:=a_1+a_3x+a_5x^2+a_7x^3\dots $$ (compare Section 9). We denote $D(m,0):=|\widehat{\goth M}(a_*,0,k_0)|$. Let us change notations: from here $L=\{\mu_1,\mu_2,\dots, \mu_n\}$ will mean the object that earlier was denoted by $\bar L$, and $A_L$ will mean the transposed to the matrix that was denoted by $A_{\bar L}$ earlier.\footnotemark \footnotetext{The authors apologise for this inconvenience.} Hence, $A_L$ is a maximal square submatrix of $\goth N(a_*,n,k)$ obtained by elimination of its $\mu_1,\mu_2,\dots, \mu_n$-th rows. 
\medskip
{\bf Proposition 7.} $\forall \ L $ \  $|A_L|=|W_L| D(m,0)$, where $$W_L=\left( \matrix a_{\mu_1-1}&a_{\mu_2-1}&\dots&a_{\mu_n-1}\\  a_{\mu_1-3}&a_{\mu_2-3}&\dots&a_{\mu_n-3}\\  \dots&\dots&\dots&\dots \\ a_{\mu_1-(2n-1)}&a_{\mu_2-(2n-1)}&\dots&a_{\mu_n-(2n-1)}\endmatrix \right)$$

{\bf Proof.} We need a chain of lemmas. 
\medskip
{\bf Lemma 7.1.} $|A_L|$ is a multiple of $D(m,0)$.
\medskip
{\bf Proof.} Let $C_i$ be the $i$-th column of $\goth N(a_*,n,k)$. If $D(m,0)=0$ then $P_{[1]}, \ P_{[0]}$ have a common root $r$. In this case we have $\sum_{i=1}^k r^{i-1}C_i=0$, i.e. the columns of $\goth N(a_*,n,k)$ are linearly dependent and hence all its maximal minors are 0. This means that if $D(m,0)=0$ then $\forall L$ we have $|A_L|=0$. Since both $D(m,0)$, $|A_L|$ are homogeneous polynomials in $a_0,\dots,a_m$, this means that $\exists \be$ such that $|A_L|^\be\in\langle D(m,0)\rangle$ (here $\langle D(m,0)\rangle$ is the ideal generated by $D(m,0)$). To prove that $D(m,0)$ is a factor of $|A_L|$ it is sufficient to prove that $D(m,0)$ is squarefree in $\n C[a_0,\dots,a_m]$, which is equivalent to: dim Sing$(X(m,1))<m-1$ where $X(m,1)\subset P^m$ is the variety of zeroes of $D(m,0)$. For a proof it is sufficient to find a straight line $P^1\subset P^m$ that crosses $X(m,1)$ in deg $X(m,1)=m-1$ distinct points. 

We define this line as the line joining two points $t=(P_{[0]}, P_{[1]})$ and $t'=(P'_{[0]}, P'_{[1]})$ where $P_{[0]}$, $ P_{[1]}$ are from above and $P'_{[0]}$, $ P'_{[1]}$ come from $a'_0,\dots,a'_m$. Let $(u:u')$ be projective coordinates of a point on $P^1$. This point belongs to $X(m,1)$ iff $P_{[0]}u+P'_{[0]}u'$, $P_{[1]}u+P'_{[1]}u'$ have a common root. If $x$ is this common root then $\goth D(P)=0$ where $\goth D(P):=\left|\matrix P_{[0]}(x) &P'_{[0]}(x) \\ P_{[1]}(x) &P'_{[1]}(x)\endmatrix \right|$. Hence, we have to find $ P_{[*]},  P'_{[*]}$ such that the equation $\goth D(P)=0$ has distinct roots $r_1,\dots,r_{m-1}$ and moreover the numbers $(-u'(r_i):u(r_i))= (P_{[0]}(r_i) :P'_{[0]}(r_i))$ are also distinct. 

For even $m=2\nu+2$ we can choose $$P_{[0]}=x^{\nu+1}-1,\ \ P'_{[0]}=1,\ \ P_{[1]}=-x^\nu,\ \ P'_{[1]}=x^\nu+1$$ $\goth D(P)$ is $x^{m-1}+x^{\nu+1}-1$ and $(u':u)=-x^{\nu+1}+1$. It is obvious that all roots of $\goth D(P)$ are distinct (the roots of its derivative are not the roots of $\goth D(P)$) and the ratios of roots are not $\zeta^*_{\nu+1}$. For odd $m=2\nu+1$ we can choose $$P_{[0]}=x^\nu-1,\ \ P'_{[0]}=1,\ \ P_{[1]}=x,\ \ P'_{[1]}=x^\nu+1$$ $\goth D(P)$ is $x^{m-1}-x-1$ and $(u':u)=-x^\nu+1$. So, $D(m,0)$ is a factor of $|A_L|$. $\square$
\medskip
By technical reasons, we impose temporary a condition on $m$: 
$$m>4n^2+6n\eqno{(7.2)}$$

{\bf Lemma 7.3.} If $m$ satisfies $(7.2)$ then $|A_L|$ is a multiple of $|W_L|D(m,0)$.
\medskip
{\bf Proof. } First, let us prove that $|W_L|$ also is a factor of $|A_L|$. We denote $\goth N=\goth N(a_*,n,k)$, and let $B$ be the following $(n\times (k+n))$-matrix: $$\left(\matrix a_{0}&-a_{1}&a_2&\dots&(-1)^m a_{m}&0&0&0&\dots&0\\  0&0&a_{0}&-a_{1}&a_2&\dots&(-1)^m a_{m}&0&\dots&0\\  \dots&\dots&\dots&\dots&\dots&\dots&\dots&\dots&\dots&\dots \\ 0&0&0&\dots&0&a_{0}&-a_{1}&a_2&\dots&(-1)^m a_{m} \endmatrix\right)$$ We have $B\goth N=0$. Let $\hat B$, resp. $\hat \goth N$ be matrices obtained from $B$, resp. $\goth N$ by a permutation of columns (resp. rows) sending columns (resp. rows) $\{\mu_1,\mu_2,\dots, \mu_n\}$ to $\{1,\dots,n\}$. We have again $\hat B\hat \goth N=0$. Let us denote $\hat B=(B_1|B_2)$, $\hat \goth N=\left(\matrix \goth N_1\\ \goth N_2 \endmatrix\right)$ the block partitions: $B_1,B_2,\goth N_1, \goth N_2$ are respectively $n\times n$, $n\times k$, $n\times k$, $k\times k$-matrices. We have $\hat B\hat \goth N=B_1\goth N_1+B_2\goth N_2$, hence $B_1\goth N_1=-B_2\goth N_2$. We have $|W_L|=\pm |B_1|$ and $|A_L|=\pm|\goth N_2|$. 
\medskip
{\bf Sublemma 7.3.1.} If $|B_1|=0$ then $|\goth N_2|=0$. 
\medskip
If $|B_1|=0$ then the rank of $B_1\goth N_1$ is $<n$. If $B_2$ has the maximal rank then $|\goth N_2|=0$. Really, if $m$ is small with respect to $n$ it can happen that $B_2$ contains a row of zeroes; it can happen also that $a_i$ take such values that the rank of $B_2$ is not maximal. Let us give a rigorous proof of the sublemma. 

We denote by $H\subset [1,\dots,k+n]$ the set of $h$ such that $a_h$ is an entry of $W_L$, and let $\bar H=[\bar h_1,\dots,\bar h_\al]$ be its complement in $[1,\dots,k+n]$. Let $a_i$, $i\in H$, take values $b_i\in \n C$ such that $|W_L|=0$ at these values. Formally, we consider a ring homomorphism $\vf: \n C[a_0,\dots,a_m]\to \n C[a_{\bar h_1},\dots,a_{\bar h_\al}]$ defined by $\vf(a_i)= b_i$ for $i\in H$, $\vf(a_i)=a_i$ for $i\in \bar H$. 

We have $\vf(B_1)\vf(\goth N_1)=-\vf(B_2)\vf(\goth N_2)$. Condition $|\vf(B_1)|=0$ implies that there exists a non-zero $(1\times n)$-matrix $C=(c_1\ \dots \ c_n)$ with entries $c_i\in \n C$ such that $C \cdot \vf(B_1)=0$, hence $C \cdot \vf(B_2)\vf(\goth N_2)=0$. Condition $(7.2)$ implies that $\exists \bar h\in \bar H$ such that $$[\bar h-1, \bar h-3,\dots,\bar h-(2n-1)]\cap H=\emptyset, \ \ [\bar h-1, \bar h-3,\dots,\bar h-(2n-1)]\subset [0,\dots,m]$$ This means that the corresponding element of the row matrix $C \cdot \vf(B_2)$ is equal to $c_1a_{\bar h-1}+c_2a_{\bar h-1}+\dots+ c_na_{\bar h-(2n-1)}\ne 0$, i.e. $C \cdot \vf(B_2)\ne0$. Let $\vf(\goth N_2)^{adj}$ be the adjoint matrix. We have $C \cdot \vf(B_2)\vf(\goth N_2)\vf(\goth N_2)^{adj}=0=|\vf(\goth N_2)|\cdot C \cdot \vf(B_2)I_k$, hence $|\vf(\goth N_2)|=0$, i.e. Sublemma 7.3.1 is proved.
\medskip
As above, to get that $|W_L|$ is a factor of $|A_L|$ we must prove that $|W_L|$ is squarefree. Really, it is irreducible, we prove it by induction for $n$. $|W_L|$ is linear as a polynomial in $a_{\mu_n-1}$: $|W_L|=\goth C_1a_{\mu_n-1}+\goth C_0$, hence its possible factor is free from $a_{\mu_n-1}$ and divides both $\goth C_1$, $\goth C_0$. We have: $\goth C_1$ is a $((n-1)\times(n-1))$-determinant of the same type as $|W_L|$, and hence it is irreducible by the induction hypothesis. Let us prove that $\goth C_1$ does not divide $\goth C_0$. We consider the lexicographic order on $\n Z[a_0,\dots,a_m]$ defined by the condition $a_0<\dots<a_m$. The highest term of $\goth C_0$ is $\pm a_{\mu_{n-1}-1}\cdot a_{\mu_{n}-3}\cdot a_{\mu_{n-2}-5}\cdot\dots \cdot a_{\mu_{1}-(2n-1)}$, and the highest term of $\goth C_1$ is $\pm a_{\mu_{n-1}-3}\cdot a_{\mu_{n-2}-5}\cdot\dots \cdot a_{\mu_{1}-(2n-1)}$ corresponding to its antidiagonal elements. It is not a factor of the highest term of $\goth C_0$, hence the desired. 

This means that to prove that $|A_L|$ is a multiple of $|W_L|D(m,0)$ we must prove that $|W_L|$ and $D(m,0)$ are coprime. Since $|W_L|$ is irreducible it is sufficient to prove that it is not a factor of $D(m,0)$. Again  $\pm a_m^{(m-1)/2}\cdot a_0^{(m-1)/2}$ for odd $m$ and $\pm a_m^{(m-2)/2}\cdot a_1^{m/2}$ for even $m$ --- the  highest term of $D(m,0)$ --- is not a multiple of the highest term of $|W_L|$. Lemma 7.3 is proved. $\square$
\medskip
This means (because of equality of degrees) that $|A_L|=c|W_L| D(m,0)$ where $c$ is a constant. We must prove that $c=1$ (it is important and not obvious that $c$ cannot be $-1$). Let us consider first the case 
\medskip
{\bf 7.4.} All entries of $W_L$ are not 0 and for odd $m$ they are not $a_0, a_m$ (i.e. $\mu_1\ge 2n$, $\mu_n\le m$), for even $m$ they are not $a_0, a_1, a_m$ (i.e. $\mu_1\ge 2n+1$, $\mu_n\le m$). 
\medskip
We shall use a terminology (following N.N. Luzin): a set of entries of a square matrix such that every row and column contains exactly one element of this set is called a lightning, and the product of these elements is called the value of the lightning.
\medskip
{\bf Lemma 7.5.} For odd $m$ satisfying (7.2) and (7.4) we have $c=1$.
\medskip
{\bf Proof. Step 7.5.1: Construction of the highest lightning of $|A_L|$.} The highest term of $|W_L| D(m,0)$ is the product of the highest terms of factors, it is equal to $$(-1)^{{n\choose 2}+{(m+1)/2\choose 2}} a_m^{(m-1)/2}\cdot a_0^{(m-1)/2} \cdot a_{\mu_n-1}\cdot  a_{\mu_{n-1}-3} \cdot\dots \cdot a_{\mu_{1}-(2n-1)}$$ Let $\la$ be a lightning of $A_L$ of value (without sign) $$a_m^{(m-1)/2}\cdot a_0^{(m-1)/2} \cdot a_{\mu_n-1}\cdot  a_{\mu_{n-1}-3} \cdot\dots \cdot a_{\mu_{1}-(2n-1)}$$ Let us prove that there exists only one such $\la$, and that its sign is $(-1)^{{n\choose 2}+{(m+1)/2\choose 2}}$. 

Let $[i_1,\dots,i_\ga]\subset [1,\dots,n]$ be numbers such that $\mu_{i_1},\dots,\mu_{i_\ga}$ are even and its complement $[j_1,\dots,j_{n-\ga}]\subset [1,\dots,n]$ be numbers such that $\mu_{j_1},\dots,\mu_{j_{n-\ga}}$ are odd. We subdivide the matrix $\goth N$ to 3 submatrices: 
\medskip
$\goth N_l$ (left) formed by the 1-st --- $\frac{m-1}2$-th columns of $\goth N$; 
\medskip
$\goth N_{md}$ (middle) formed by $\frac{m+1}2$-th --- $(\frac{m-1}2+n)$-th columns of $\goth N$; 
\medskip
$\goth N_r$ (right) formed by $(\frac{m+1}2+n)$-th --- $k$-th columns of $\goth N$. 
\medskip
A column of $\goth N_l$, resp. of $\goth N_{md}$, $\goth N_r$, contains the element $a_0$ and does not contain the element $a_m$, resp. contains both the elements $a_0$ and $a_m$, resp. contains the element $a_m$ and does not contain the element $a_0$. We consider the same subdivision of the matrix $A_L$, it is the union of $A_{L,l}$, $A_{L,md}$, $A_{L,r}$. 

Any row of $\goth N$ contains exactly one of the elements $a_0$, $a_m$. Conditions $\mu_1\ge 2n$, $\mu_n\le m$ imply that $\forall c, \ 1\le c\le\ga$, the $\mu_{i_c}$-th row of $\goth N$ contains $a_0$ in $\goth N_l$, and $\forall c, \ 1\le c\le n-\ga$, the $\mu_{j_c}$-th row of $\goth N$ contains $a_m$ in $\goth N_r$. This means that $A_{L,l}$ contains $\frac{m-1}2-\ga$ elements $a_0$ and no $a_m$, $A_{L,md}$ contains $n$ elements $a_0$ and $n$ elements $a_m$, and $A_{L,r}$ contains $\frac{m-1}2-(n-\ga)$ elements $a_m$ and no $a_0$. $\la$ contains $m-1$ elements $a_0$ and $a_m$. It cannot contain more than $n$ of these elements from $A_{L,md}$, hence $\la$ must contain all $\frac{m-1}2-\ga$ elements $a_0$ from $A_{L,l}$, all $\frac{m-1}2-(n-\ga)$ elements $a_m$ from $A_{L,r}$, $\ga$ elements $a_0$ and $n-\ga$ elements $a_m$ from $A_{L,md}$. 

The only columns of $A_L$ such that the element of $\la$ of these columns are neither $a_0$ nor $a_m$ are columns $\mu_{i_1}/2,\dots,\mu_{i_\ga}/2$ in $A_{L,l}$ and $(\mu_{j_1}+m)/2,\dots,(\mu_{j_{n-\ga}}+m)/2$ in $A_{L,r}$. An element of $\la$ in a column $\mu_{i_c}/2$, where $1\le c\le\ga$, must be $a_\de$ for $\de$ odd, because elements $a_\de$ for $\de$ even are in the rows containing $a_0$ in $A_{L,l}$, and all $a_0$ in $A_{L,l}$ belong to $\la$ --- a contradiction. Analogically, an element of $\la$ in a column $(\mu_{j_c}+m)/2$, where $1\le c\le n-\ga$, must be $a_\de$ for $\de$ even. 

Further, an element of $\la$ in a column $\mu_{i_c}/2$, where $1\le c\le\ga$, must be in the $\tau$-th row where $\tau\le 2n$. Really, if $\tau>2n$ is odd (numbering of $\goth N$) then the $\tau$-th row contains the element $a_m$ in $A_{L,r}$, it belongs to $\la$ --- a contradiction. If $\tau>2n$ is even (numbering of $\goth N$) then the $\tau$-th row contains the element $a_0$ in $A_{L,r}$ or $A_{L,md}$. If $a_0$ is in $A_{L,r}$ then it belongs to $\la$ --- a contradiction. If $a_0$ is in $A_{L,md}$ then the $(\mu_{i_c}/2,\tau)$-th entry of $\goth N$, and hence $A_L$, is 0. Finally, $\tau$ must be odd, because for even $\tau$ the $\tau$-th line contains $a_0$ in $A_{L,l}$ belonging to $\la$. 

Therefore, we consider a $(n\times \ga)$-submatrix $U_l$ ($l$ means left) of $A_L$ formed by $1,3,\dots,2n-1$-th rows and by $\mu_{i_1}/2,\dots,\mu_{i_\ga}/2$-th columns. By the symmetry with respect to the center of $\goth N$, we get that all elements of $\la$ in $(\mu_{j_1}+m)/2,\dots,(\mu_{j_{n-\ga}}+m)/2$-th columns have the number of row $\in [k,k-2,\dots, k-(2n-2)]$ (numbering of $A_L$). We denote by $U_r$ ($r$ means right) the $(n\times (n-\ga))$-submatrix of $A_L$ formed by these rows and columns. Let $U$ be their union --- a $n\times n$-matrix. Let us show that the elements of $\la$ in $U_l$, $U_r$, treated as elements of $U$ (we denote this set by $w$), form a lightning in $U$. Really, each column of $U$ contains only one element of $w$. Let a $z$-th row of $U$ contains an element of $w$. This means that the element $a_m$ on $(2z-1,(m-1)/2+z)$-th position in $A_L$ does not belong to $\la$. Since the $((m-1)/2+z)$-th column of $A_L$ belongs to $A_{L,md}$, we get that $a_0$ in this column belongs to $\la$. It is easy to see that this $a_0$ is in a row of $A_L$ which corresponds to the $z$-th row of $U_r$, hence the $z$-th row of $U_r$ does not contain elements of $w$. By the symmetry of properties of $U_l$ and $U_r$ (or because $\#(w)=n$) we conclude that $w$ is a lightning in $U$. 

It is easy to see that $U$ is $W_L$ with a permutation of columns (first are columns having even $\mu_*$, second are columns having odd $\mu_*$). Since the antidiagonal of $W_L$ is the only lightning with its value, we get unicity of $\la\subset A_L$, hence $c=\pm1$. We shall need a fact (which follows immediately) that the above $\tau$ is $2n-2i_c+1$.
\medskip
Let us calculate the sign of $\la$. It is easier first to show that the neighbor $L$, $L'$ have the same sign, and then to find this sign for one fixed $L$. 
\medskip
{\bf Step 7.5.2. Equality of signs of neighbor $L$.} Two sets $L=(\mu_1,\dots,\mu_n)$ and $L'=(\mu'_1,\dots,\mu'_n)$, where $\mu_1<\dots<\mu_n$, $\mu'_1<\dots<\mu'_n$, are called neighbor if $\exists i$ such that $\forall j\ne i$ we have $\mu_j=\mu'_j$ and $\mu_i+1=\mu'_i$. We can assume that $\mu_i$ is odd. The matrices $A_L$, $A_{L'}$ differ in only one row --- the $\mu_i-(i-1)$-th row (here and below numbering of $A_L$). We denote by $\la$, $\la'$ the above lightnings of $A_L$, $A_{L'}$ respectively. According the above construction of $\la$, we get that $\la'$ contains the $(\mu_i-(i-1), \psi)$-th entry of $A_{L'}$ where $\psi$ is some number, this entry is $a_0$. $\la$ contains a $(\mu_i-(i-1), \psi+(m-1)/2)$-th entry of $A_{L}$, this entry is $a_m$. We denote --- as above --- by $\tau=2n-2i+1$ the number of row such that $\la$ contains the $(\tau,\psi)$-th entry of $A_L$. This means that $\la$ does not contain the element $a_m$ at the $\tau$-th row of $A_L$. This element is at the $(\tau,(\tau+m)/2)$-th position in $A_L$. Therefore, the element of $\la$ in the $(\tau+m)/2)$-th column of $A_L$ is the element $a_0$ of this column, it is at the $(\tau+m-n,(\tau+m)/2)$-th position of $A_L$. 

Analogically, $\la'$ contains (in addition to the $(\mu_i-(i-1), \psi)$-th entry mentioned above) also the $(\tau,(\tau+m)/2)$-th entry of $A_{L'}$ which is $a_m$, and the $(\tau+m-n, \psi+(m-1)/2)$-th entry of $A_{L'}$, and all other entries --- except these 3 entries --- of $\la$ and $\la'$ coincide. This can be shown by explicit calculation of $\psi$ (it is a function of the quantities of odd and even $\mu_j$ for $j<i$), or we can use the fact that if $\mu_i$ is odd then the $(2n-2i+1,\frac{m+1}2+n-i)$-th entry of $A_L$ (it is $a_m$) belongs to $\la$, and if $\mu_i$ is even then the $(m+n-2i+1,\frac{m+1}2+n-i)$-th entry of $A_L$ (it is $a_0$) belongs to $\la$ (this follows easily from the above fact that $U_l\cup U_r$ is $W_L$ with a permutation of columns, and because the (image under the permutation of columns of the) lightning $w$ is the antidiagonal of $W_L$). 

We get that there are two triples 
\medskip
$(t_1,t_2,t_3):=(\tau,\mu_i-(i-1),\tau+m-n)$ and 
\medskip
$(s_1,s_2,s_3):=(\psi,(\tau+m)/2,\psi+(m-1)/2)$ 
\medskip
such that $\la$ contains the $(t_1,s_1)$-th, $(t_2,s_3)$-th, $(t_3,s_2)$-th entries of $A_L$, 
\medskip
$\la'$ contains the $(t_1,s_2)$-th, $(t_2,s_1)$-th, $(t_3,s_3)$-th entries of $A_{L'}$, 
\medskip
and other $k-3$ entries of $\la$, $\la'$ coincide. This implies that $\la$ and $\la'$ have the same parity. Since any two $L$ can be joined by a chain of neighbors, we get that for fixed $m$, $n$ the coefficent $c$ does not depend on $L$. 
\medskip
{\bf Step 7.5.3. Calculation of sign of some given $L$.} We consider $L=(2n,2n+2,\dots,4n-2)$. The description of $\la$ given in Step 1 shows that for this case $\la$ is the disjoint union of 5 sets (we indicate positions and values of their entries):
\medskip
1. $(2\al,\al)$, value $a_0$, $\al=1,\dots,n-1$;
\medskip
2. $(2\al-1,2n-\al)$, value $a_{4n+1-4\al}$, $\al=1,\dots,n$;
\medskip
3. $(2n-1+\al,(m-1)/2+n+\al)$, value $a_m$, $\al=1,\dots,n$;
\medskip
4. $(3n+2\al,2n+\al)$, value $a_0$, $\al=0,\dots,(m-1)/2-n$;
\medskip
5. $(3n-1+2\al,(m-1)/2+2n+\al)$, value $a_m$, $\al=1,\dots,(m-1)/2-n$.
\medskip
We interchange the elements of $\la$ such that the numbers of columns form the sequence $(1,\dots,k)$. The numbers of rows are the following (we indicate to which of the above sets (1) - (5) they belong):
\medskip
\noindent
$2,4,\dots,2n-2$ (set 1), $2n-1,2n-3,\dots,1$ (set 2), $3n,3n+2,\dots,m+n-1$ (set 4), $2n,2n+1,\dots,3n-1$ (set 3), $3n+1,3n+3,\dots, m+n-2$ (set 5). 
\medskip
The quantity of inversions of this sequence is ${n\choose2}+{(m+1)/2\choose2}$. Lemma 7.5 is proved. $\square$
\medskip
{\bf End of the proof of Proposition 7.} For even $m$ the proof is similar. Throughout the proof $a_0$ will be changed to $a_1$, the highest term of $D(m,0)$ is $(-1)^{m/2\choose2}a_m^{(m-2)/2}\cdot a_1^{m/2}$, etc. The proof is omitted (we can use also the below reduction from $m+1$ to $m$). 

To proceed from $m+1$ to $m$, we consider a map $\vf:\n Z[a_0,\dots,a_{m+1}]\to \n Z[a_0,\dots,a_{m}]$ sending $a_{m+1}$ to 0. As earlier let $L=(\mu_1,\dots,\mu_n)$ where $\mu_n\le m+2n-1$ be the same for $m$ and for $m+1$. Dependence of $W_L$, $A_L$ on $m$ will be indicated explicitly. We have 
\medskip
$\vf(|W_L(m+1)|)=|W_L(m)|$, $\frac{\vf(D(m+1,0))}{a_m}=D(m,0)$, $\frac{\vf(|A_L(m+1)|)}{a_m}=|A_L(m)|$.
\medskip
This implies immediately that if the lemma is true for $m$ then it is true for $m-1$. Now we use symmetry: the lemma is stable with repect to the $m$-symmetry denoted by $\sigma_m$: $\sigma_m(a_i)=a_{m+1-i}$, for $L=(\mu_1,\dots,\mu_n)$ we have $\sigma_m(L)= (k+n+1-\mu_n,\dots,k+n+1-\mu_1)$ ($\sigma_m$ does not change the sign of all involved terms). 
\medskip
Therefore, let $m,n,L$ be arbitary. For a sufficiently large $m_1$ the $\sigma_{m_1}(L)$ "satisfies (7.4) from the left", i.e. its $\mu_1$ is $\ge 2n+1$, and for a sufficiently large odd $m_2$ the $\sigma_{m_1}(L)$ "satisfies (7.4) from the right", i.e. its $\mu_n$ is $\ge m_2$, and the condition $m_2>4n^2+6n$ also holds. So, Lemma 7.5 holds for $m_2,n,\sigma_{m_1}(L)$. The above decreasing from $m_2$ to $m_1$ shows that the lemma holds for $m_1,n,\sigma_{m_1}(L)$. The symmetry shows that the lemma holds for $m_1,n,L$. The decreasing from $m_1$ to $m$ shows that Proposition 7 holds for all $m,n,L$.
$\square$

{
\medskip
Substituting the values of $c(L)$ and $|A_L|$ to (5) we get that (1) is equivalent to the formula 
$$\sum_{L\in\goth L} \frac{\prod_{1\le i<j\le n}(\mu_j-\mu_i)}{(n-1)!!} \left| \matrix a_{\mu_1-1}&a_{\mu_2-1}&\dots&a_{\mu_n-1}\\  a_{\mu_1-3}&a_{\mu_2-3}&\dots&a_{\mu_n-3}\\  \dots&\dots&\dots&\dots \\ a_{\mu_1-(2n-1)}&a_{\mu_2-(2n-1)}&\dots&a_{\mu_n-(2n-1)}\endmatrix \right|t^{(\mu_1+...+\mu_n)}=$$ $$=2^{n\choose 2}t^{n^2}(a_0+a_1t+a_2t^2+\dots +a_mt^m)^n\eqno{(8)}$$
It is sufficient to prove that $\forall \ r_1,r_2,\dots,r_n$, $0\le r_i\le m$, the numerical coefficient at $a_{r_1}a_{r_2}\dots a_{r_n}$ in both left and right hand sides of (8) are equal (it is clear that the degrees of $t$ entering to the coefficient at $a_{r_1}a_{r_2}\dots a_{r_n}$ in both left and right hand sides of (8) are equal). We denote $R=(r_1,r_2,\dots,r_n)$. Changing the order of $r_i$ if necessary we can assume that $$r_1=r_2=\dots=r_{\al_1}$$ $$r_{\al_1+1}=r_{\al_1+2}=\dots=r_{\al_1+\al_2}$$ $$\dots$$ $$r_{n-\al_c+1}=r_{n-\al_c+2}=\dots=r_{n}$$ and there is no more equalities between $r_i$. The segments of consecutive length $\al_1,\al_2,\dots,\al_c$ in the segment $[1,\dots,n]$ will be called the standard segments. The coefficient at $a_{r_1}a_{r_2}\dots a_{r_n}$ in the right hand side of (8) is $2^{n\choose 2}\cdot {n\choose{\al_1,\dots,\al_c}}$.

We denote by $S_\al$ the subgroup of $S_n$ consisting of permutations that stabilize the standard segments.  We have $S_\al=S_{\al_1}\times S_{\al_2}\times\dots\times S_{\al_c}$. For all $\s\in S_n$ ($S_n$ acts on $1,2,\dots, n$) we define an ordered sequence
$$R_\s=\{r_1+(2\cdot\s(1)-1), r_2+(2\cdot\s(2)-1), \dots, r_n+(2\cdot\s(n)-1)\}$$

We consider a subset $\goth L_R$ of $\goth L$ consisting of $L=(\mu_1,\mu_2,\dots,\mu_n)$, where $\mu_1<\mu_2<\dots<\mu_n$, such that $\exists \s\in S_n$ such that $L$ is a permutation (denoted by $\de=\de(\s)$) of $R_\s$ (since all $\mu_*$ are different, $\de$ is defined uniquely by $\s$). It is clear that the left hand side of (8) does not depend on the order of $\mu_*$. 

For $L\in \goth L_R$ we denote by $N(L)=N_R(L) \subset S_n$ the set of $\s$ such that $L$ is a permutation of $R_\s$. Moreover for $\chi=+$ or $-$ we denote by $N_{\chi}(L)\subset N(L)$ the set of the above $\s$ such that the parity of $\s\cdot\de(\s)$ is $\chi$. We denote $\nu_\chi(L)=\#(N_{\chi}(L))$.
\medskip
{\bf Proposition 9.} The matrix $W_L$ contains $\frac{\nu_{+}(L)}{\al_1!\al_2!\cdot\dots\cdot\al_c!}$, resp. $\frac{\nu_{-}(L)}{\al_1!\al_2!\cdot\dots\cdot\al_c!}$ even, resp. odd lightnings of value $a_{r_1}a_{r_2}\cdot\dots\cdot a_{r_n}$. 
\medskip
{\bf Proof.} Let $\s\in N(L)$. We fix the direction of the action of $\de$ by the formula $\mu_{\de(i)}=r_i+(2\cdot\s(i)-1)$, $i=1,\dots,n$. We define a map $\la$ form $N(L)$ to the set of lightnings of $W_L$ of order $a_{r_1}a_{r_2}\cdot\dots\cdot a_{r_n}$ as follows: for $\s\in N(L)$ we have: $\la(\s)$ is the set of $(1,\s\de^{-1}(1))$-, 
$(2,\s\de^{-1}(2))$-, $\dots$, $(n,\s\de^{-1}(n))$-th entries of $W_L$. Surjectivity of $\la$ is obvious, as well as the fact that if $W(L)$ contains a lightning of value $a_{r_1}a_{r_2}\cdot\dots\cdot a_{r_n}$ then $L\in \goth L_R$. Let us find the order of a fiber of $\la$. Let $\la(\s)=\la(\s')$. We have $\s\de^{-1}=\s'{\de'}^{-1}$ (here $\de'=\de(\s')$ ). Further, $\forall i$ the $(i,\s\de^{-1}(i))$-th entry of $W_L$ is $a_{r_{\de^{-1}(i)}}$, hence --- because for $\s$, $\s'$ not only positions of the elements of the lightning coincide, but the elements themselves, we get $r_{\de^{-1}(i)}=r_{{\de'}^{-1}(i)}$,  $r_{{\de}^{-1}\de'(i)}=r_i$. This means that ${\de}^{-1}\de'\in S_\al$. 

Conversely, for all $\de'\in S_\al\de$ and $\s':=\s\de^{-1}\de'$ we have: $R_{\s'}$ is a permutation of the same $L$, and $\la(\s)=\la(\s')$. This means that the order of all fibers of $\la$ is $\#S_\al=\al_1!\al_2!\cdot\dots\cdot\al_c!$  Finally, the parity of the lightning $\la(\s)$ is the parity of $\s\de$, hence the proposition . $\square$
\medskip
This proposition  implies that to prove (8) we have to prove that $\forall \ R$ 
$$\sum_{L=(\mu_1,\dots,\mu_n)\in\goth L_R} \left(\prod_{1\le i<j\le n}(\mu_j-\mu_i)\right) (\nu_{+}(L)-\nu_{-}(L))=2^{n\choose 2}n!!\eqno{(10)}$$
\medskip
{\bf Proposition 11.} $\forall \ R=(r_1,\dots,r_n)$ the left hand side of (25) is 
$$\sum_{\s\in S_n} \sg(\s)\prod_{1\le i<j\le n} ((r_j+(2\s(j)-1))-(r_i+(2\s(i)-1)))\eqno{(12)}$$

{\bf Proof.} If $\s\not\in\cup_{L\in \goth L_R}N_R(L)$, i.e. if $R_\s$ contains two equal numbers, then the corresponding term of the sum (12) is 0. The union $\cup_{L\in \goth L_R}N_R(L)$ is disjoint, i.e. $L$ is defined by $\s$ uniquely, and it is sufficient to prove that $\forall \ L=(\mu_1,\dots,\mu_n)\in \goth L_R$ we have 
$$\left(\prod_{1\le i<j\le n}(\mu_j-\mu_i)\right) (\nu_{+}(L)-\nu_{-}(L))=$$ $$=\sum_{\s\in N_R(L)} \sg(\s)\prod_{1\le i<j\le n} ((r_j+(2\s(j)-1))-(r_i+(2\s(i)-1)))\eqno{(13)}$$
This is clear, because the set of $\mu_i$ is a $\de$-permutation of the set $r_i+(2\s(i)-1)$, and $\nu_{+}(L)$, resp. $\nu_{-}(L)$ are the quantities of $\s\in N_R(L)$ such that $\s\de$ is even, resp. odd. $\square$
\medskip
To prove that (12) is $2^{n\choose 2}n!!$ we need two lemmas. 
Let $x_1,\dots,x_n$ be abstract variables. Let us consider a $n!\times{n\choose 2}$-matrix $A$ whose lines are numbered by $\sigma\in S_n$ and whose columns are numbered by pairs $(i,j)$ such that $1\le i<j\le n$, and defined as follows: $$A_{\sigma,(ij)}=x_{\s(j)}-x_{\s(i)}$$ Let, further, $S$ be a subset of $B:=$ the set of the columns of $A$. 
\medskip
{\bf Lemma 14.} If $S\ne B$ then $\sum_{\s\in S_n} \sg(\s)\prod_{\ga\in S}A_{\s\ga}=0$. 
\medskip
{\bf Proof.} $\sum_{\s\in S_n} \sg(\s)\prod_{\ga\in S}A_{\s\ga}$ is an alternating polynomial in $x_1,\dots,x_n$ of degree $\#S$. There is no alternating polynomials in $n$ variables of degrees $<{n\choose 2}$. $\square$
\medskip
{\bf Lemma 15.} Let $r_1,\dots,r_n$ be abstract variables, $1,3,\dots, 2n-1$ the set of odd numbers, the group $S_n$ acts on this set, and $i,j$ as above. Then $$\sum_{\s\in S_n} \sg(\s)\prod_{(i,j)\in B} ((r_j+\s(2j-1))-(r_i+\s(2i-1)))=2^{n\choose 2}n!!\eqno{(16)}$$ (particularly, the left hand side of (16) does not depend on $r_1,\dots,r_n$). 
\medskip
{\bf Proof.} We shall prove a more general equality. Instead of $1,3,\dots, 2n-1$ we consider any set of numbers $x_1,\dots,x_n$, and instead of $r_j-r_i$ we consider independent variables $\la_\al$, where $\al=(i,j)\in B$. The left hand side of (16) becomes $$\sum_{\s\in S_n} \sg(\s)\prod_{\al=(i,j)\in B}(\la_\al+(x_{\s(j)}-x_{\s(i)})$$
The above lemma shows that $\forall S\ne\emptyset$ coefficients at $\prod_{\al\in S} \la_\al$ are 0, and the $\la_*$-free term (for the case $x_i=2i-1$) is $n!\cdot \prod_{(i,j)\in B} \ ((2j-1)-(2i-1))=2^{n\choose 2}n!!$ $\square$
\medskip
The formula (16) implies the theorem III. $\square$
\medskip
{\bf 17. Proof B.} (5) and (7.1) show that $|\widehat {\goth M}(a_*,n,k)|$ is a multiple of $D(m,0)$. Further, $|\widehat {\goth M}(a_*,n,k)|$ is a polynomial in $t$; its $\la$-th derivative is a linear combination of $|\widehat {\goth M}(a_*,n,k) ^{(\la_1)(\la_2)\dots(\la_n)}|$ where $\la_1+\dots+\la_n=\la$ and for a matrix $M$ the $M^{(\la_1)(\la_2)\dots(\la_n)}$ is a matrix whose $i$-th row is the $\la_i$-th derivative of the $i$-th row of $M$, $i=1,\dots,n$. Hence, to prove that $|\widehat {\goth M}(a_*,n,k)|$ is a multiple of $(\sum_{i=0}^m \ a_it^i)^n$ it is sufficient to prove
\medskip
{\bf Proposition 18.} $\forall \ \la<n$, $\forall \ \la_1,\dots,\la_n\ge0$ such that $\la_1+\dots+\la_n=\la$ we have $|\widehat {\goth M}(a_*,n,k) ^{(\la_1)(\la_2)\dots(\la_n)}|$ is a multiple of $\sum_{i=0}^m \ a_it^i$. 
\medskip
{\bf Proof.} \ \ $\widehat {\goth M}(a_*,n,k) ^{(\la_1)(\la_2)\dots(\la_n)}\cdot \left(\matrix 1\\ t^2\\t^4\\ \dots \\t^{2k-2}\endmatrix \right)$ is a multiple of $\sum_{i=0}^m \ a_it^i$. $\square$
\medskip
\noindent
$\sum_{i=0}^m \ a_it^i$ is not a divisor of $D(m,0)$, because $D(m,0)$ is $t$-free, hence $|\widehat {\goth M}(a_*,n,k)|$ is a multiple of $(\sum_{i=0}^m \ a_it^i)^n\cdot D(m,0)$. Factor $t^{n\choose2}$ appears in any of $k!$ terms of $|\widehat {\goth M}(a_*,n,k)|$. Finding of the numerical coefficient $2^{n\choose2}$ is like in the Proof A. $\square$
\medskip
{\bf Appendix. Some auxiliary results and remarks.}
\medskip
{\bf 1. Proof of Proposition 8.31.} We consider the affine part $a_0\ne 0$ of $P^m$. We can take $a_0=1$. We consider an affine $q-1$-dimensional linear space $H_{q-1}\subset A^m$ defined parametrically: $$a_i=c_{i0}+c_{i1}t_1+\dots+c_{i,q-1}t_{q-1}\eqno{(A1.1)}$$ where $i=1,\dots,m$, $t_1,\dots, t_{q-1}$ are parameters of this space and $c_{ij}$ are arbitrary constants. We find $\#(H_{q-1}\cap X_r(q,m,1)$. We substitute (A1.1) to polynomials $P_{[j]}$ from 8.28, $j=0,\dots,q-1$. The $j$-th equation $P_{[j]}(x)=0$ becomes 
$$P_{j0}(x)+P_{j1}(x)t_1+\dots +P_{j,q-1}(x)t_{q-1}=0\eqno{(A1.2)}$$ where $P_{ji}(x)$ are polynomials in $x$ of degree $\gamma_j$, where $\gamma_j$ --- the degree of $P_{[j]}$ --- is the maximal number such that $q\gamma_j+j\le m$. The system (A1.2) has a solution iff $|P_{ji}(x)|=0$. Since $|P_{ji}(x)|$ is a polynomial of degree $d:=\sum_{j=0}^{q-1}\gamma_j$, we get that $\#(H_{q-1}\cap X_r(q,m,1))=d$. Obviously $\sum_{j=0}^{q-1}\gamma_j=m+1-q$.
\medskip
{\bf 2. Proof of Conjecture 8.23 for $X(3,1,m,1)$. }
\medskip
We need a definition: an $l$-quasidiagonal of a $k\times k$-matrix is the set of its $(i,i+l)$-entries, where $i$ runs over $1,\dots, k$, and $i+l$ mod $k$. Analogous definition holds for a $(k+1)\times k$-matrix ($i$ runs over $1,\dots, k+1$). A $(k,l_1,l_2)$-biquasidiagonal matrix $A=(a_{ij})$ is a $k\times k$-matrix having non-zero entries only on its $l_1$- and $l_2$-quasidiagonals: $a_{i,i+l_1}=c_i$, $a_{i,i+l_2}=d_i$. Obviously if $(k,l_1-l_2)=1$ then $|A|=\pm\prod_i c_i\pm\prod_i d_i$.

Let us consider for simplicity the case $q=3$. After a permutation of rows of $\goth N(P,1,k)$ we get the following Sylvester-type matrix $\goth N'(P,1,k)$ whose block structure is $\left(\matrix B_2\\B_1\\B_0\endmatrix \right)$ and the blocks have the form 

$B_i=\left(\matrix a_i&a_{3+i}&a_{6+i}&\dots &a_{\mu}&0&0 &\dots &0\\ 0&a_i&a_{3+i}&a_{6+i}&\dots &a_{\mu}&0 &\dots &0 \\ 0&0&a_i&a_{3+i}&a_{6+i}&\dots &a_{\mu} &\dots &0\\ \dots &  \dots &  \dots &  \dots &  \dots &  \dots &  \dots &  \dots &  \dots \\ 0&\dots&0&0&a_i&a_{3+i}&a_{6+i}&\dots &a_{\mu}\endmatrix \right)$ 
\medskip
where $\mu$ is the maximal number satisfying $\mu\le m$, $\mu\equiv i$ mod 3. If $k\not\equiv 0$ mod 3 then $\exists l_1, l_2$ such that $(k,l_1-l_2)=1$ and all elements on $l_1$- and $l_2$-quasidiagonals of $\goth N'(P,1,k)$ are not 0. Moreover, there is a subset $\{i_1,\dots,i_6\}\subset \{1,\dots,m\}$ such that $a_i\in l_1$-quasidiagonal $\cup\ l_2$-quasidiagonal $\iff i\in \{i_1,\dots,i_6\}$. 

We consider $L'$ given by the equations $a_i=0$ if $i\not\in\{i_1,\dots,i_6\}$, and we consider the corresponding $(k,l_1,l_2)$-biquasidiagonal $(k+1)\times k$-matrix $\goth N''(P,1,k)$. It is sufficient to prove that it contains two $k\times k$-minors whose determinants are coprime. We can take the extreme 0-th and $k$-th minors $M_0$, $M_k$: they are $(k,l_1,l_2)$-biquasidiagonal $k\times k$-matrices, their determinants are $\pm a_{i_1}^*a_{i_2}^*a_{i_3}^*\pm a_{i_4}^*a_{i_5}^*a_{i_6}^*$, where for brevity we do not give here the exact values of * --- they are $\sim k/3$, they depend on $k$ mod 3 and on the minor (0-th and $k$-th). It is faster to the reader to check himself that in all cases  $|M_0|$, $|M_k|$ are coprime than to understand a written proof. 

All other cases ($k$ is a multiple of 3; $q>3$) are treated by the similar manner. We omit the details. 
\medskip
{\bf 3. Cyclicity of rows of $\goth M$ is essential.} 
\medskip
Analogs of the determinantal varieties $X(q,n,m,i)=X(i)$ have meaning for the following matrices $\widetilde {\goth M}^-(*,n,k)$ which are more general than $\goth M(P,n,k)$ (in the below example we take $n=1$; its version without minus signs for the case $n=2$ is given in (III.4), the general form is clear):

$$\widetilde {\goth M}^-(*,1,k)=\left(\matrix a_{11}t-a_{12} &  a_{12}t-a_{13}  &   \dots & a_{1k}t-a_{1,k+1} \\ a_{21}t-a_{22} 
 &  a_{22}t-a_{23} &   \dots & a_{2k}t-a_{2,k+1} \\
a_{31}t-a_{32} &  a_{32}t-a_{33}  &   \dots & a_{3k}t-a_{3,k+1} \\
\dots & \dots & \dots  & \dots \\ a_{k1}t-a_{k2}
 &   a_{k2}t-a_{k3}  &  \dots & a_{kk}t-a_{k,k+1} \endmatrix \right) \eqno{(A3.1)}$$
where $a_{ij}$ are arbitrary. Lemma 8.12 shows that for this case we have Codim $X(1)=n+1$ as well. It is natural to ask:
\medskip
whether Codim $X(i+1)$ in $X(i)$ is $n+1$, or not? 
\medskip
Answer: not, this can be shown for $n=1$, $k=3$, $i=1$ even by hand calculation:
\medskip
{\bf Proposition.} Let us consider the space $P^{11}$ --- the variety of matrices (A3.1) for $k=3$. The subvarieties $X(i)$ are defined for it, and we have: $X(2)$ is a complete intersection in $X(1)$, i.e. the codimension of $X(2)$ in $X(1)$ is 3 unlike 2 for the case of $X(i)$ of Definition 8.4. 
\medskip
{\bf Proof.} According Lemma 8.12, $(a_{11}:...:a_{34})\in X(1)\iff \goth N(P,1,3)$ has rank 2, where for the present case $\goth N(P,1,3)^t=\left(\matrix a_{11} &\dots & a_{14}\\ \dots &\dots &\dots \\ a_{31}&\dots &a_{34}\endmatrix \right)$. On an open part of $X(1)$ this implies that $$\left(\matrix  a_{31}&\dots &a_{34}\endmatrix \right)=\lambda_1\left(\matrix  a_{11}&\dots &a_{14}\endmatrix \right)+\lambda_2\left(\matrix  a_{21}&\dots &a_{24}\endmatrix \right)\eqno{(A3.2)}$$ We denote by $d_{ij}$, $1\le i < j\le 4$, the determinant of the $(i,j)$-th minor of $\left(\matrix a_{11} &\dots & a_{14}\\  a_{21}&\dots &a_{24}\endmatrix \right)=(i,j)$-th Pl\"ucker coordinate of $( a_{11} \dots a_{14})\wedge (a_{21}\dots  a_{24})$. 
\medskip
For $(a_{11}:...:a_{34})\in X(1)$ we have $(a_{11}:...:a_{34})\in X(2)\iff $ the coefficient at $U$ of det $(U\cdot I_k-\tilde\goth M(P,1,3))$ is 0. Taking into consideration A3.2, this is equivalent to the condition $$(1 \ \lambda_1 \ \lambda_2)D=0 \hbox{ where } D=\left(\matrix d_{12}&-d_{13}&d_{23}\\ -d_{23}&d_{24}&-d_{34}\\ d_{13}&-d_{14}-d_{23} & d_{24}\endmatrix \right)\eqno{(A3.3)}$$
\medskip
If the codimension of $X(2)$ in $X(1)$ is 2 then $\forall a_{11}, \dots ,a_{24}\ \exists \lambda_1,\lambda_2$ satisfying A3.3. This means that $\forall a_{11}, \dots ,a_{24}$ we have det $D=0$. This can happen only if det $D$ is a multiple of $d_{12}d_{34} -d_{13}d_{24}+d_{14}d_{23}$ --- the only relation between $d_{ij}$. We have det $D=-(-d_{14} d_{23}^2 - d_{23}^3 + 2 d_{13} d_{23} d_{24 }- d_{12} d_{24}^2 - d_{13}^2 d_{34} +
d_{12} d_{14} d_{34} + d_{12} d_{23} d_{34}) $, and obviously it is not a multiple of $d_{12}d_{34} -d_{13}d_{24}+d_{14}d_{23}$ --- a contradiction. $\square$
\medskip
{\bf 4. Computer evidence for Conjectures 9.7.} We consider a random affine space $Y$ defined over $\n Q$ of complementary dimension $d$, i.e. $d=$ codimension of $X_i=X(2,1,m,i)$ in $P^m$, and $Y\cap X_i$ is a finite set. Let $y_1,\dots,y_d$ be affine coordinates of $Y$, hence $a_0,\dots, a_m$ are their linear combinations. We substitute these linear combinations to $H_{jl}$ of (8.3.1), where $j, \ l$ run over the set of indices $S(k,i):=\{j=1,\dots,i-1, \ \ l=0,\dots,k-j\} \cup (0,0)$ (because we know that all $H_{0i}$ are proportional, see Proposition 9.12), and we denote the obtained polynomials in $y_1,\dots,y_d$ by $H_{jl}=H_{jl}(Y)$ as well. 

Definitions of the resultant of $n$ polynomials $P_1, \dots, P_n$ can be found for example in [GKZ]. We use the following inductive formula for its calculation. Let $P_1, \dots, P_n$ be polynomials in $n$ variables $x_1,\dots, x_n$, $P_j=\sum_I a_{j,I}x^I$ where $I$ is a multiindex and $a_{j,I}$ abstract coefficients, then
$R_{x_1,\dots, x_{n-1}}(P_1, \dots, P_n)$ --- the resultant of $P_1, \dots, P_n$ in variables $x_1,\dots, x_{n-1}$ --- can be calculated as follows:
$$R_{x_1,\dots, x_{n-1}}(P_1, \dots, P_n)=$$ $$GCD\{R_{x_{n-1}}[R_{x_1,\dots, x_{n-2}}(P_1, \dots, P_{n-1}), R_{x_1,\dots, x_{n-2}}(P_1, \dots, P_{n-2},P_n)], $$ $$ R_{x_{n-1}}[R_{x_1,\dots, x_{n-2}}(P_1, \dots, P_{n-1}), R_{x_1,\dots, x_{n-2}}(P_1, \dots, P_{n-3}, P_{n-1},P_n)]\} $$ where $R_{x_{n-1}}$ is the ordinary resultant in 1 variable of 2 polynomials. If $P_i$ are homogeneous of degrees $d_i$, then $R_{x_1,\dots, x_{n-1}}(P_1, \dots, P_n)$ is a polynomial in $x_n$ of degree $\prod d_i$. 

We consider all possible $d$-uples among the above polynomials $H_{jl}$. Namely, let $\al_i=(\al_{i1},\dots,\al_{id})$ be a subset of $S(k,i)$), and let $R_{\al_i}:=R_{y_1,\dots,y_{d-1}}(H_{\al_{i1}},...,H_{\al_{id}})$ be the corresponding resultant. We let $R:=GCD(R_{\al_1},\dots,R_{\al_k})$, where $\al_1,\dots,\al_k$ satisfy $\al_1\cup\dots\cup \al_k=S(k,i)$. We have $R\in \n Q[y_d]$. We can expect that for a generic $Y$ the set of $y_d$-coordinates of the points of $Y\cap X_i$ coincides with the set of roots of $R$. Again, for a generic $Y$ factorization of $R$ over $\n Q[y_d]$ corresponds to representation of $X_i$ as a sum of $\n Q$-irreducible divisors. 
\medskip

{\bf References}
\medskip
[A86] Anderson, Greg W. $t$-motives.  Duke Math. J.  53  (1986),  no. 2, 457 -- 502
\medskip
[A00] Anderson, Greg W. An elementary approach to $L$-functions mod $p$. J. Number theory 80 (2000), no. 2, 291 -- 303.
\medskip
[B02]  B\"ockle, Gebhard. Global $L$-functions over function fields. Math. Ann. 323 (2002), no. 4, 737 -- 795. 
\medskip
[B05] B\"ockle, Gebhard. Arithmetic over function fields: a cohomological approach. Number fields and function fields -- two parallel worlds, 1 -- 38, Progr. Math., 239, Birkhäuser Boston, Boston, MA, 2005.
\medskip
[BP] B\"ockle, Gebhard; Pink, Richard. Cohomological theory of crystals over function fields. EMS tracts in Mathematics, 9. European Mathematical Society (EMS), Z\"urich, 2009. viii+187 pp.
\medskip
[B12] B\"ockle, Gebhard. Cohomological theory of crystals over function fields and applications.
In "Arithmetic Geometry in Positive Characteristic". Advanced Courses in Mathematics. CRM Barcelona. Birkh\"auser Verlag, Basel (2012)
\medskip
[FP] Fulton, William; Pragacz, Piotr. Schubert varieties and degeneracy loci. LNM 1689. 
\medskip
[GKZ] Gelfand, Israel; Kapranov, Mikhail; Zelevinsky, Andrei. Discriminants, resultants, and multidimensional determinants. Birkh\"auser Verlag, Boston, Basel, Z\"urich, 1994
\medskip
[G] Goss, David. $L$-series of $t$-motives and Drinfeld modules. The arithmetic of function fields (Columbus, OH, 1991), 313 -- 402, Ohio State Univ. Math. Res. Inst. Publ., 2, de Gruyter, Berlin, 1992.
\medskip
[K] Krattenthaler, Christian. Advanced determinant calculus. 

arxiv.org/pdf/math/9902004v3.pdf
\medskip
[L] Lafforgue, Vincent. Valeurs sp\'eciales des fonctions $L$ en caract\'eristique $p$.  J. Number Theory  129  (2009),  no. 10, 2600 -- 2634
\medskip
[P] Pragacz, Piotr. A note on the elimination theory, Indagationes Math., 90, (1987), p. 215-221.
\medskip
[Sh] Shimura, Goro. Introduction to the arithmetic theory of automorphic functions. 1971. 
\medskip
[TW] Taguchi, Y.; Wan, D. $L$-functions of $\vf$-sheaves and Drinfeld modules. J. Amer. Math. Soc. 9 (1996), no. 3, 755 -- 781
\medskip
[T] Thakur, Dinesh S. On characteristic $p$ zeta functions. Compositio Math. 99 (1995), no. 3, 231 -- 247.
\medskip
[TVN] Tsfasman, Michael A.; Vl\u adu\c t, Serge; Nogin, Dmitry. Algebraic geometric codes: basic notions. 2007. 

\enddocument